\documentclass[reqno,12pt]{amsart}
\usepackage{amsmath, amssymb}
\usepackage{amsthm, amsbsy}
\usepackage{graphicx}
\usepackage{float}
\usepackage[ruled]{caption}

\theoremstyle{plain}
\newtheorem{theorem}{\indent\bf Theorem}[section]

\newtheorem{corollary}[theorem]{\indent\bf Corollary}

\theoremstyle{definition}

\begin{document}

\title[Irreducibility criteria]{Sch\" onemann-Eisenstein-Dumas-type irreducibility conditions that use arbitrarily many prime numbers}
\author[N.C. Bonciocat]{Nicolae Ciprian Bonciocat}
\address{Simion Stoilow Institute of Mathematics of the Romanian 
Academy, Research Unit 5,
P.O. Box 1-764, Bucharest 014700, Romania}
\email{Nicolae.Bonciocat@imar.ro}

\keywords{Irreducible polynomials, prime numbers, Newton polygon}
\subjclass[2000]{Primary 11R09; Secondary 11C08.}

\begin{abstract}
The famous irreducibility criteria of Sch\" onemann-Eisenstein and Dumas rely on information on the divisibility of the coefficients of a polynomial by a single prime number.
In this paper we provide several irreducibility criteria of Sch\" onemann-Eisenstein-Dumas-type for polynomials  
with integer coefficients, criteria that are given by some divisibility conditions for their coefficients with respect to arbitrarily many prime numbers.
A special attention will be paid to those irreducibility criteria that require information on the divisibility of the coefficients by two distinct prime numbers.
\end{abstract}
\maketitle

\section{Introduction} \label{se1}

The most famous irreducibility criterion is probably the one of Sch\" onemann and Eisenstein, that first appeared in a disguised form as a corollary of the following less known irreducibility criterion
of Sch\" onemann \cite{Schonemann} in 1846, and four years later in a paper of Eisenstein \cite{Eisenstein}.
\medskip

{\bf  Irreducibility criterion of Sch\" onemann} {\em  
Suppose that a polynomial $f(X)\in \mathbb{Z}[X]$ has the form $f(X)=\phi (X)^e+pM(X)$, where $p$ is a prime number, $\phi(X)$ is an irreducible polynomial modulo $p$, and $M(X)$ is
a polynomial relatively prime to $\phi (X)$ modulo $p$, with $\deg M<\deg f$. Then $f$ is irreducible over $\mathbb{Q}$.
}
\medskip

{\bf  Irreducibility criterion of Sch\" onemann-Eisenstein} {\em  
Let $f(X)=a_{0}+a_{1}X+\cdots +a_{n}X^{n}$ be a polynomial with integer coefficients, and let
$p$ be a prime number. If $p\nmid a_{n}$, $p\mid a_{i}$ for $i=0,\dots ,n-1$, and 
$p^{2}\nmid a_{0}$, then $f$ is irreducible over $\mathbb{Q}$.
}
\medskip

For an interesting review of the history of these results, and of some of the techniques used in their proof, as well as in the proof of some other important achievements of the $19$th century number theory, we refer the reader to \cite{Cox}.
Over the years many authors contributed to the development of the techniques used in the study of the irreducibility of polynomials, some of them generalizing in various ways these irreducibility criteria.
In 1895 K\" onigsberger \cite{Konigsberger}, and in 1896 Netto \cite{Netto} generalized the Sch\" onemann-Eisenstein criterion. Later, in 1905 K\" onigsberger's criterion was extended by Bauer \cite{Bauer}
and Perron \cite{Perron}.
A strengthened, more general version of Sch\" onemann-Eisenstein criterion was obtained by Dumas \cite{Dumas} in 1906, who had major contributions in the development of the Newton polygon method, one of the most powerful methods to study the irreducibility of polynomials.
\medskip

{\bf  Irreducibility criterion of Dumas} {\em  
Let $f(X)=a_{0}+a_{1}X+\cdots +a_{n}X^{n}$ be a polynomial with integer coefficients, and let
$p$ be a prime number. If
\smallskip

i) $\frac{\nu _{p}(a_{i})}{i}>\frac{\nu _{p}(a_{n})}{n}$ for $i=1,\dots ,n-1$,

ii) $\nu _{p}(a_{0})=0$,

iii) $\gcd(\nu _{p}(a_{n}),n)=1$,
\smallskip

\noindent then $f$ is irreducible over $\mathbb{Q}$.
}
\medskip

Here and henceforth, for an integer $n$ and a prime number $p$, $\nu _{p}(n)$ stands for the largest integer $i$ such that $p^{i}\mid n$ \ (by convention, $\nu _{p}(0)=\infty $).
\medskip

Further irreducibility criteria that rely on the use of Newton polygon method have been obtained by Kurschak \cite{Kurschak} in 1923, Ore \cite{Ore1}, \cite{Ore2}, \cite{Ore3}
in 1923 and 1924, and Rella \cite{Rella} in 1927. In 1938 MacLane \cite{MacLane} obtained a general result that includes all of these criteria as special cases. For some decades after MacLane's
paper, the interest in using Newton polygon method in factorization problems diminished, at least in the sense that the number of references to applications of this method decreased. 
In 1987 Coleman \cite{Coleman} used Newton polygons to study Galois groups of the exponential Taylor polynomials, and in 1995
Mott \cite{Mott} used Newton polygons to give new proofs of older results, and to obtain new irreducibility criteria as well. In recent years one of the most important applications of Newton polygon method was in the study
of the irreducibility of Bessel polynomials and Laguerre polynomials, and in the study of their Galois groups, and here we refer the reader to the work of Filaseta \cite{Filaseta1}, \cite{Filaseta2}, 
Filaseta and Lam \cite{FilasetaLam}, Filaseta and Trifonov \cite{FilasetaTrifonov}, Filaseta and Williams \cite{FilasetaWilliams}, Filaseta, Finch and Leidy \cite{Filaseta-Finch-Leidy}, 
Filaseta, Kidd and Trifonov \cite{Filaseta-Kidd-Trifonov}, Hajir \cite{Hajir1}, \cite{Hajir2}, \cite{Hajir3}, Hajir and Wong \cite{HajirWong}, and Sell \cite{Sell}. Other important applications of Newton polygons
have been obtained by Greve and Pauli \cite{GrevePauli} in the study of the splitting fields and Galois groups of Eisenstein polynomials.

More recent generalizations of the irreducibility criteria of Sch\" onemann and Sch\" onemann-Eisenstein, and results related to the study of generalized Sch\" onemann and Eisenstein-Dumas
polynomials have been obtained by Panaitopol and \c Stef\u anescu \cite{PanaitopolStefanescu}, Khanduja and Saha \cite{KhandujaSaha}, Brown \cite{Brown}, Bush and Hajir \cite{BushHajir}, Bishnoi and Khanduja \cite{BishnoiKhanduja}, 
Khanduja and Khassa \cite{KhandujaKhassa}, and Weintraub \cite{Weintraub}.

The aim of this paper is to provide some irreducibility criteria of Sch\" onemann-Eisenstein-Dumas-type for some classes of polynomials with integer coefficients, 
expressed by some divisibility conditions for their coefficients with respect to arbitrarily many prime numbers.

The main idea used to obtain the results in this paper is to gather information on the degrees of the non-constant factors of a given polynomial $f$ by studying its Newton polygons with respect to a family of
prime numbers that divide some of its coefficients, and then to search for some conditions that prevent these data to fit together, unless $f$ is irreducible. 
The idea to study the degrees of the non-constant factors of a given polynomial $f$ by studying its Newton polygons with respect to more than a single prime number goes back to Dumas, who proved the following result \cite[page 238, section 4]{Dumas}.
\medskip

{\em Si $p,p',\dots ,p^{(\nu )}$ sont des nombres premieres quelconques et si 
\[
\sum \mu _{i}s_{i},\ \sum \mu '_{i}s'_{i},\dots ,\sum \mu ^{(\gamma )}_{i}s^{(\gamma )}_{i},\dots ,\sum \mu ^{(\nu )}_{i}s^{(\nu )}_{i}
\]
repr\' esentent les sommes qui leur correspondent respectivement, lorsque pour chacun d'eux on construit le polygone correspondant \` a $f(x)$, seuls les polynomes, dont le degr\' e est susceptible
d'\^ etre egal, simultan\' ement, \` a une somme de la forme de chacune des expressions $\sum \mu ^{(\gamma )}_{i}s^{(\gamma )}_{i}$, pourront \^ etre diviseurs de $f(x)$.
}
\medskip

Here $\mu ^{(\gamma )}_{i}$ are non-negative integers, and $s^{(\gamma )}_{i}$ are widths of the segments in the Newton polygon of $f$ with respect to $p^{(\gamma )}$.

The results that we will prove in this paper rely on this fundamental theorem of Dumas, which can be easily turned into an irreducibility criterion in the following way.
\medskip



{\bf Theorem A.}\ {\em
Let $f(X)\in\mathbb{Z}[X]$ be a polynomial of degree $n$, let $k\geq 2$, and let
$p_{1},\dots ,p_{k}$ be pair-wise distinct prime numbers. For $i=1,\dots , k$ let us denote by $w_{i,1},\dots ,w_{i,n_{i}}$ the widths of all the segments of  the Newton polygon of $f$ with respect to $p_{i}$, and by $\mathcal{S}_{p_{i}}$ the set of all the integers in the interval  $(0,\lfloor \frac{n}{2}\rfloor ]$ that may be written as a linear combination of $w_{i,1},\dots ,w_{i,n_{i}}$ with coefficients $0$ or $1$. 
If $\mathcal{S}_{p_{1}}\cap\cdots \cap\mathcal{S}_{p_{k}}=\emptyset$,
then $f$ is irreducible over $\mathbb{Q}$.
}
\medskip

Note that in the statement of Theorem A, the integers $w_{i,1},\dots ,w_{i,n_{i}}$ are not necessarily distinct. In other words, Theorem A states that if for each $i=1,\dots , k$ we denote by $\mathcal{S}_{p_{i}}$ the set of all the possible lengths in the interval  $(0,\lfloor \frac{n}{2}\rfloor ]$ of the projections onto the $x$-axis of the (not necessarily connected) broken lines obtained by removing some of the segments of the Newton polygon of $f$ with respect to $p_{i}$, and if $\mathcal{S}_{p_{1}}\cap\cdots \cap\mathcal{S}_{p_{k}}=\emptyset$, then the polynomial $f(X)$ must be irreducible over $\mathbb{Q}$. Here the length of the projection onto the $x$-axis of a disconnected broken line, refers to the sum of the lengths of the projections onto the $x$-axis of its connected components. 

Our second result that relies on the information on the Newton polygons of $f$ with respect to different prime numbers is the following special case of Theorem A.
\medskip

{\bf Theorem B.}\ {\em 
Let $f(X)\in\mathbb{Z}[X]$ be a polynomial of degree $n$, let $k\geq 2$, let
$p_{1},\dots ,p_{k}$ be pair-wise distinct prime numbers, and for $i=1,\dots ,k$ let $d_{p_{i}}$ denote the greatest common divisor of the widths of the segments of  the Newton polygon of $f$ with respect to $p_{i}$. Then the degree of any factor of $f$ is divisible by $n/\gcd(\frac{n}{d_{p_{1}}},\dots ,\frac{n}{d_{p_{k}}})$. In particular, if $\gcd(\frac{n}{d_{p_{1}}},\dots ,\frac{n}{d_{p_{k}}})=1$,
then $f$ is irreducible over $\mathbb{Q}$.
}
\medskip

We note that Theorem B extends to the case of Newton polygons with arbitrarily many edges Proposition 3.10 in \cite{Mott}, that refers only to Newton polygons consisting of a single edge.
We also note here that Theorems A and B may be rephrased in the following effective form, by describing in an explicit form the sets $\mathcal{S}_{p_{1}},\dots ,\mathcal{S}_{p_{k}}$ and the numbers
$d_{p_{1}},\dots ,d_{p_{k}}$ too, in terms of the abscisae of the vertices of the Newton polygons, and of the $p_{i}$-adic valuations of the coefficients of $f$. 
\medskip

{\bf Theorem A'.}\ {\em 
Let $f(X)=a_{0}+a_{1}X+\cdots+a_{n}X^{n}\in\mathbb{Z}[X]$, $a_{0}a_{n}\neq 0$, let $k\geq 2$, and let
$p_{1},\dots ,p_{k}$ be pair-wise distinct prime numbers. For $i=1,\dots ,k$ let us denote the abscisae of the vertices in the Newton polygon of $f$ with respect to $p_{i}$ by the sequence of integers $0=j_{i,1}<j_{i,2}<\cdots<j_{i,r_{i}}=n$, and let
\begin{eqnarray*}
m_{i,l} & = & \gcd(\nu _{p_{i}}(a_{j_{i,l+1}})-\nu _{p_{i}}(a_{j_{i,l}}), j_{i,l+1}-j_{i,l}),\quad l=1,\dots, r_{i}-1,\\
x_{i,l} & = & \frac{j_{i,l+1}-j_{i,l}}{m_{i,l}}, \qquad \qquad \qquad \qquad \qquad \qquad \thinspace \thinspace l=1,\dots, r_{i}-1,\\
\mathcal{S}_{p_{i}} & = & \left\lbrace  \sum\limits _{l=1}^{r_{i}-1} n_{i,l}x_{i,l} : n_{i,l}\in \{ 0,1,\dots , m_{i,l} \} \right\rbrace \cap \left( 0,\left\lfloor \frac{n}{2}\right\rfloor \right] .
\end{eqnarray*}
If $\mathcal{S}_{p_{1}}\cap\cdots \cap\mathcal{S}_{p_{k}}=\emptyset$, then $f$ is irreducible over $\mathbb{Q}$.
}
\medskip

{\bf Theorem B'.}\ {\em 
Let $f(X)=a_{0}+a_{1}X+\cdots+a_{n}X^{n}\in\mathbb{Z}[X]$, $a_{0}a_{n}\neq 0$, let $k\geq 2$, and let
$p_{1},\dots ,p_{k}$ be pair-wise distinct prime numbers. For $i=1,\dots ,k$ let us denote the abscisae of the vertices in the Newton polygon of $f$ with respect to $p_{i}$ by the sequence of integers $0=j_{i,1}<j_{i,2}<\cdots<j_{i,r_{i}}=n$, and let
$d_{p_{i}}=\gcd(x_{i,1},\dots ,x_{i,r_{i}-1})$, with
\[
x_{i,l} = \frac{j_{i,l+1}-j_{i,l}}{\gcd(\nu _{p_{i}}(a_{j_{i,l+1}})-\nu _{p_{i}}(a_{j_{i,l}}), j_{i,l+1}-j_{i,l}  )}, \qquad l=1,\dots, r_{i}-1.
\]
Then the degree of any factor of $f$ is divisible by $n/\gcd(\frac{n}{d_{p_{1}}},\dots ,\frac{n}{d_{p_{k}}})$. In particular, if $\gcd(\frac{n}{d_{p_{1}}},\dots ,\frac{n}{d_{p_{k}}})=1$,
then $f$ is irreducible over $\mathbb{Q}$.
}
\medskip

Here we obviously understand that $d_{p_{i}}=x_{i,1}=\frac{n}{\gcd(\nu _{p_{i}}(a_{n})-\nu _{p_{i}}(a_{0}),\thinspace n)}$ if $r_{i}=2$ for some $i$.
\medskip

The proof of these results is unexpectedly simple, in contrast to their high level of generality. 
One might obviously expect that these results have a lot of corollaries, that mainly depend on the geometry of the Newton polygons of $f$ with respect to $p_{1},\dots ,p_{k}$.
To find a suitable notation for the possible corollaries of Theorem B, for instance, might be itself an extremely difficult task. Indeed, even if we restrict our attention to Newton polygons having the number of edges bounded
by a certain positive integer $m$, say, we are still left with the major task of finding all the relevant combinations of shapes for the Newton polygons that we will consider for each $p_{i}$, and the number of such possibilities will be quite overwhelming. A partial answer to this problem, at least for small values of $m$ and $k$, is to attach to each corollary a label consisting of a multi-index $(n_{1},n_{2},\dots ,n_{k})$ with
$n_{1}\leq n_{2}\leq \dots \leq n_{k}$, where $n_{i}$ is the number of edges in the Newton polygon of $f$ with respect to $p_{i}$, followed by an usual additional index, having no geometrical meaning. For instance, the first corollary for the case when $k=2$, $n_{1}=3$ and $n_{2}=4$ should be called Corollary (3,4).1, the second one being labelled as Corollary (3,4).2, and so on. A Corollary to such a result will then obviously receive an additional index, etc..  In this way, such a label will only contain the information on the numbers of edges, while the remaining information on the precise shape of the Newton polygons will be found only in the 
statement of that corollary. However, in this paper we will not use this unusual notation. We will content ourselves to state the main corollaries to Theorems A and B as theorems too, and their immediate consequences as corollaries, thus using an ordinary labelling.
Besides, once we state an irreducibility criterion, we will pay no attention at all to the corresponding irreducibility criterion obtained by considering the reciprocal $\overline{f}=X^{\deg f}f(1/X)$ of $f$
instead of $f$. The reason is that such a ``new'' irreducibility criterion will only use the coefficients in reverse order, and will actually bring no new important information. The main applications of Theorems A and B
that use information on the divisibility of the coefficients of $f$ with respect to a pair of distinct prime numbers $p$ and $q$, in case each one of the Newton polygons of $f$ with respect to $p$ and $q$ consists of at most two edges, will correspond to the $14$ relevant combinations of shapes of Newton polygons
displayed in Figure 1 below. By relevant we understand that none of two such combinations can be obtained by the other by considering the reverse order of the coefficients of $f$.
For instance, we don't consider to be relevant the combination of two Newton polygons each one consisting of a single edge with negative slope, since such a combination may be obtained by the one plotted in 
Figure 1.b by reversing the order of the coefficients of $f$, that is by considering the reciprocal of $f$ instead of $f$.
\medskip

\begin{center}

\setlength{\unitlength}{2mm}
\begin{picture}(58,7)
\linethickness{0.15mm}

{\small \put(6,3){$a$}}

\put(12,0){\line(0,1){3.2}}   
\thicklines

\put(0,4){\line(3,-1){12}}   
\linethickness{0.15mm}

\put(0,0){\circle{0.12}}
\put(0,0){\circle{0.08}}
\put(0.5,0.1){\circle{0.08}}
\put(1,0.2){\circle{0.08}}
\put(1.5,0.3){\circle{0.08}}
\put(2,0.4){\circle{0.08}}
\put(2.5,0.5){\circle{0.08}}
\put(3,0.6){\circle{0.08}}
\put(3.5,0.7){\circle{0.08}}
\put(4,0.8){\circle{0.08}}
\put(4.5,0.9){\circle{0.08}}
\put(5,1){\circle{0.08}}
\put(5.5,1.1){\circle{0.08}}
\put(6,1.2){\circle{0.08}}
\put(6.5,1.3){\circle{0.08}}
\put(7,1.4){\circle{0.08}}
\put(7.5,1.5){\circle{0.08}}
\put(8,1.6){\circle{0.08}}
\put(8.5,1.7){\circle{0.08}}
\put(9,1.8){\circle{0.08}}
\put(9.5,1.9){\circle{0.08}}
\put(10,2){\circle{0.08}}
\put(10.5,2.1){\circle{0.08}}
\put(11,2.2){\circle{0.08}}
\put(11.5,2.3){\circle{0.08}}
\put(12,2.4){\circle{0.08}}
\put(12,2.4){\circle{0.12}}

\put(0,0){\vector(1,0){13.5}}
\put(0,0){\vector(0,1){5.5}}

\put(0,4){\circle{0.12}}
\put(0,4){\circle{0.08}}

{\small \put(21,3){$b$}}

\setlength{\unitlength}{2mm}

\linethickness{0.15mm}

\put(27,0){\line(0,1){4.7}}   
\thicklines

\put(15,0){\line(3,1){12}}   
\linethickness{0.15mm}

\put(15,0){\circle{0.12}}
\put(15,0){\circle{0.08}}
\put(15.5,0.1){\circle{0.08}}
\put(16,0.2){\circle{0.08}}
\put(16.5,0.3){\circle{0.08}}
\put(17,0.4){\circle{0.08}}
\put(17.5,0.5){\circle{0.08}}
\put(18,0.6){\circle{0.08}}
\put(18.5,0.7){\circle{0.08}}
\put(19,0.8){\circle{0.08}}
\put(19.5,0.9){\circle{0.08}}
\put(20,1){\circle{0.08}}
\put(20.5,1.1){\circle{0.08}}
\put(21,1.2){\circle{0.08}}
\put(21.5,1.3){\circle{0.08}}
\put(22,1.4){\circle{0.08}}
\put(22.5,1.5){\circle{0.08}}
\put(23,1.6){\circle{0.08}}
\put(23.5,1.7){\circle{0.08}}
\put(24,1.8){\circle{0.08}}
\put(24.5,1.9){\circle{0.08}}
\put(25,2){\circle{0.08}}
\put(25.5,2.1){\circle{0.08}}
\put(26,2.2){\circle{0.08}}
\put(26.5,2.3){\circle{0.08}}
\put(27,2.4){\circle{0.08}}
\put(27,2.4){\circle{0.12}}

\put(15,0){\vector(1,0){13.5}}
\put(15,0){\vector(0,1){5.5}}

\put(27,4){\circle{0.12}}
\put(27,4){\circle{0.08}}

{\small \put(36,3){$c$}}

\setlength{\unitlength}{2mm}

\linethickness{0.15mm}

\put(42,0){\line(0,1){3.2}}   
\thicklines

\put(30,4){\line(5,-3){6.6666}}   

\put(36.6666,0){\line(3,1){5.3334}}   

\linethickness{0.15mm}

\put(30,0){\circle{0.12}}
\put(30,0){\circle{0.08}}
\put(30.5,0.1){\circle{0.08}}
\put(31,0.2){\circle{0.08}}
\put(31.5,0.3){\circle{0.08}}
\put(32,0.4){\circle{0.08}}
\put(32.5,0.5){\circle{0.08}}
\put(33,0.6){\circle{0.08}}
\put(33.5,0.7){\circle{0.08}}
\put(34,0.8){\circle{0.08}}
\put(34.5,0.9){\circle{0.08}}
\put(35,1){\circle{0.08}}
\put(35.5,1.1){\circle{0.08}}
\put(36,1.2){\circle{0.08}}
\put(36.5,1.3){\circle{0.08}}
\put(37,1.4){\circle{0.08}}
\put(37.5,1.5){\circle{0.08}}
\put(38,1.6){\circle{0.08}}
\put(38.5,1.7){\circle{0.08}}
\put(39,1.8){\circle{0.08}}
\put(39.5,1.9){\circle{0.08}}
\put(40,2){\circle{0.08}}
\put(40.5,2.1){\circle{0.08}}
\put(41,2.2){\circle{0.08}}
\put(41.5,2.3){\circle{0.08}}
\put(42,2.4){\circle{0.08}}
\put(42,2.4){\circle{0.12}}

\put(30,0){\vector(1,0){13.5}}
\put(30,0){\vector(0,1){5.5}}

\put(36.6666,0){\circle{0.12}}

\put(30,4){\circle{0.12}}
\put(30,4){\circle{0.08}}

\put(42,1.7778){\circle{0.12}}
\put(42,1.7778){\circle{0.08}}

{\small \put(51,3){$d$}}
\setlength{\unitlength}{2mm}

\linethickness{0.15mm}

\put(57,0){\line(0,1){3.3}}   

\thicklines

\put(45,0){\line(1,0){6.6666}}   
\put(51.6666,0){\line(3,1){5.3334}}   
\linethickness{0.15mm}

\put(45,0){\circle{0.12}}
\put(45,0){\circle{0.08}}
\put(45.5,0.1){\circle{0.08}}
\put(46,0.2){\circle{0.08}}
\put(46.5,0.3){\circle{0.08}}
\put(47,0.4){\circle{0.08}}
\put(47.5,0.5){\circle{0.08}}
\put(48,0.6){\circle{0.08}}
\put(48.5,0.7){\circle{0.08}}
\put(49,0.8){\circle{0.08}}
\put(49.5,0.9){\circle{0.08}}
\put(50,1){\circle{0.08}}
\put(50.5,1.1){\circle{0.08}}
\put(51,1.2){\circle{0.08}}
\put(51.5,1.3){\circle{0.08}}
\put(52,1.4){\circle{0.08}}
\put(52.5,1.5){\circle{0.08}}
\put(53,1.6){\circle{0.08}}
\put(53.5,1.7){\circle{0.08}}
\put(54,1.8){\circle{0.08}}
\put(54.5,1.9){\circle{0.08}}
\put(55,2){\circle{0.08}}
\put(55.5,2.1){\circle{0.08}}
\put(56,2.2){\circle{0.08}}
\put(56.5,2.3){\circle{0.08}}
\put(57,2.4){\circle{0.08}}
\put(57,2.4){\circle{0.12}}

\put(45,0){\vector(1,0){13.5}}
\put(45,0){\vector(0,1){5.5}}

\put(51.6666,0){\circle{0.12}}

\put(57,1.7778){\circle{0.12}}
\put(57,1.7778){\circle{0.08}}

\end{picture}
\end{center}

\begin{center}

\setlength{\unitlength}{2mm}
\begin{picture}(58,7)

{\small \put(6,3){$e$}}
\linethickness{0.15mm}

\put(12,0){\line(0,1){3.3}}   

\thicklines
\put(0,4){\line(5,-3){6.6666}}   

\put(6.6666,0){\line(1,0){5.3334}}   
\linethickness{0.15mm}

\put(0,0){\circle{0.12}}
\put(0,0){\circle{0.08}}
\put(0.5,0.1){\circle{0.08}}
\put(1,0.2){\circle{0.08}}
\put(1.5,0.3){\circle{0.08}}
\put(2,0.4){\circle{0.08}}
\put(2.5,0.5){\circle{0.08}}
\put(3,0.6){\circle{0.08}}
\put(3.5,0.7){\circle{0.08}}
\put(4,0.8){\circle{0.08}}
\put(4.5,0.9){\circle{0.08}}
\put(5,1){\circle{0.08}}
\put(5.5,1.1){\circle{0.08}}
\put(6,1.2){\circle{0.08}}
\put(6.5,1.3){\circle{0.08}}
\put(7,1.4){\circle{0.08}}
\put(7.5,1.5){\circle{0.08}}
\put(8,1.6){\circle{0.08}}
\put(8.5,1.7){\circle{0.08}}
\put(9,1.8){\circle{0.08}}
\put(9.5,1.9){\circle{0.08}}
\put(10,2){\circle{0.08}}
\put(10.5,2.1){\circle{0.08}}
\put(11,2.2){\circle{0.08}}
\put(11.5,2.3){\circle{0.08}}
\put(12,2.4){\circle{0.08}}
\put(12,2.4){\circle{0.12}}

\put(0,0){\vector(1,0){13.5}}
\put(0,0){\vector(0,1){5.5}}

\put(6.6666,0){\circle{0.12}}

\put(0,4){\circle{0.12}}
\put(0,4){\circle{0.08}}

\put(12,0){\circle{0.12}}
\put(12,0){\circle{0.08}}

{\small \put(21,3.2){$f$}}
\setlength{\unitlength}{2mm}

\linethickness{0.15mm}

\put(27,0){\line(0,1){5.4}}  

\thicklines

\put(15,0){\line(6,1){5}}   
\put(20,0.83333){\line(2,1){7}}   
\linethickness{0.15mm}

\put(15,0){\circle{0.12}}
\put(15,0){\circle{0.08}}
\put(15.5,0.2){\circle{0.08}}
\put(16,0.4){\circle{0.08}}
\put(16.5,0.6){\circle{0.08}}
\put(17,0.8){\circle{0.08}}
\put(17.5,1){\circle{0.08}}
\put(18,1.2){\circle{0.08}}
\put(18.5,1.4){\circle{0.08}}
\put(19,1.6){\circle{0.08}}
\put(19.5,1.8){\circle{0.08}}
\put(20,2){\circle{0.08}}
\put(20.5,2.2){\circle{0.08}}
\put(21,2.4){\circle{0.08}}
\put(21.5,2.6){\circle{0.08}}
\put(22,2.8){\circle{0.08}}
\put(22.5,3){\circle{0.08}}
\put(23,3.2){\circle{0.08}}
\put(23.5,3.4){\circle{0.08}}
\put(24,3.6){\circle{0.08}}
\put(24.5,3.8){\circle{0.08}}
\put(25,4){\circle{0.08}}
\put(25.5,4.2){\circle{0.08}}
\put(26,4.4){\circle{0.08}}
\put(26.5,4.6){\circle{0.08}}
\put(27,4.8){\circle{0.08}}
\put(27,4.8){\circle{0.12}}

\put(15,0){\vector(1,0){13.5}}
\put(15,0){\vector(0,1){5.5}}

\put(20,0.8333){\circle{0.12}}
\put(20,0.8333){\circle{0.08}}

\put(27,4.333){\circle{0.12}}
\put(27,4.333){\circle{0.08}}

{\small \put(36,3){$g$}}
\setlength{\unitlength}{2mm}

\linethickness{0.15mm}

\put(42,0){\line(0,1){3.3}}   

\thicklines
\put(30,4){\line(1,-1){2}}   

\put(32,2){\line(5,-1){10}}   
\linethickness{0.15mm}

\put(30,0){\circle{0.12}}
\put(30,0){\circle{0.08}}
\put(30.5,0.1){\circle{0.08}}
\put(31,0.2){\circle{0.08}}
\put(31.5,0.3){\circle{0.08}}
\put(32,0.4){\circle{0.08}}
\put(32.5,0.5){\circle{0.08}}
\put(33,0.6){\circle{0.08}}
\put(33.5,0.7){\circle{0.08}}
\put(34,0.8){\circle{0.08}}
\put(34.5,0.9){\circle{0.08}}
\put(35,1){\circle{0.08}}
\put(35.5,1.1){\circle{0.08}}
\put(36,1.2){\circle{0.08}}
\put(36.5,1.3){\circle{0.08}}
\put(37,1.4){\circle{0.08}}
\put(37.5,1.5){\circle{0.08}}
\put(38,1.6){\circle{0.08}}
\put(38.5,1.7){\circle{0.08}}
\put(39,1.8){\circle{0.08}}
\put(39.5,1.9){\circle{0.08}}
\put(40,2){\circle{0.08}}
\put(40.5,2.1){\circle{0.08}}
\put(41,2.2){\circle{0.08}}
\put(41.5,2.3){\circle{0.08}}
\put(42,2.4){\circle{0.08}}
\put(42,2.4){\circle{0.12}}

\put(30,0){\vector(1,0){13.5}}
\put(30,0){\vector(0,1){5.5}}

\put(32,2){\circle{0.12}}
\put(32,2){\circle{0.08}}

\put(30,4){\circle{0.12}}
\put(30,4){\circle{0.08}}

\put(42,0){\circle{0.12}}
\put(42,0){\circle{0.08}}

{\small \put(51,3){$h$}}
\setlength{\unitlength}{2mm}

\linethickness{0.15mm}

\put(57,0){\line(0,1){3.3}}   
\thicklines

\put(45,4){\line(5,-3){6.6666}}   

\put(51.6666,0){\line(3,1){5.3334}}   

\linethickness{0.15mm}

\put(45,1.5){\circle{0.12}}
\put(45,1.5){\circle{0.08}}
\put(45.5,1.25){\circle{0.08}}
\put(46,1){\circle{0.08}}
\put(46.5,0.75){\circle{0.08}}
\put(47,0.5){\circle{0.08}}
\put(47.5,0.25){\circle{0.08}}
\put(48,0){\circle{0.12}}
\put(48,0){\circle{0.08}}

\put(48.5,0.15){\circle{0.08}}
\put(49,0.3){\circle{0.08}}
\put(49.5,0.45){\circle{0.08}}
\put(50,0.6){\circle{0.08}}
\put(50.5,0.75){\circle{0.08}}
\put(51,0.9){\circle{0.08}}
\put(51.5,1.05){\circle{0.08}}
\put(52,1.2){\circle{0.08}}
\put(52.5,1.35){\circle{0.08}}
\put(53,1.5){\circle{0.08}}
\put(53.5,1.65){\circle{0.08}}
\put(54,1.8){\circle{0.08}}
\put(54.5,1.95){\circle{0.08}}
\put(55,2.1){\circle{0.08}}
\put(55.5,2.25){\circle{0.08}}
\put(56,2.4){\circle{0.08}}
\put(56.5,2.55){\circle{0.08}}
\put(57,2.7){\circle{0.08}}

\put(45,0){\vector(1,0){13.5}}
\put(45,0){\vector(0,1){5.5}}

\put(51.6666,0){\circle{0.12}}

\put(45,4){\circle{0.12}}
\put(45,4){\circle{0.08}}

\put(57,1.7778){\circle{0.12}}
\put(57,1.7778){\circle{0.08}}

\end{picture}
\end{center}

\begin{center}

\setlength{\unitlength}{2mm}
\begin{picture}(58,7)
{\small \put(6,3){$i$}}
\linethickness{0.15mm}

\put(12,0){\line(0,1){4.5}}   

\thicklines
\put(0,0){\line(1,0){4}}   

\put(4,0){\line(2,1){8}}   

\linethickness{0.15mm}

\put(0,4){\circle{0.12}}
\put(0,4){\circle{0.08}}
\put(0.5,3.75){\circle{0.08}}
\put(1,3.5){\circle{0.08}}
\put(1.5,3.25){\circle{0.08}}
\put(2,3){\circle{0.08}}
\put(2.5,2.75){\circle{0.08}}
\put(3,2.5){\circle{0.08}}

\put(3.5,2.25){\circle{0.08}}
\put(4,2){\circle{0.08}}
\put(4.5,1.75){\circle{0.08}}
\put(5,1.5){\circle{0.08}}
\put(5.5,1.25){\circle{0.08}}
\put(6,1){\circle{0.08}}
\put(6.5,0.75){\circle{0.08}}
\put(7,0.5){\circle{0.08}}
\put(7.5,0.25){\circle{0.08}}
\put(8,0){\circle{0.08}}

\put(8.5,0.25){\circle{0.08}}
\put(9,0.5){\circle{0.08}}
\put(9.5,0.75){\circle{0.08}}
\put(10,1){\circle{0.08}}
\put(10.5,1.25){\circle{0.08}}
\put(11,1.5){\circle{0.08}}
\put(11.5,1.75){\circle{0.08}}
\put(12,2){\circle{0.08}}
\put(12,2){\circle{0.12}}

\put(0,0){\vector(1,0){13.5}}
\put(0,0){\vector(0,1){5.5}}

\put(0,0){\circle{0.12}}
\put(0,0){\circle{0.08}}

\put(4,0){\circle{0.12}}
\put(4,0){\circle{0.08}}

\put(0,4){\circle{0.12}}
\put(0,4){\circle{0.08}}

\put(12,4){\circle{0.12}}
\put(12,4){\circle{0.08}}

{\small \put(21,3){$j$}}
\setlength{\unitlength}{2mm}

\linethickness{0.15mm}

\put(27,0){\line(0,1){5}}   

\thicklines
\put(15,0){\line(5,1){5}}   

\put(20,1){\line(2,1){7}}   

\linethickness{0.15mm}

\put(15,4){\circle{0.12}}
\put(15,4){\circle{0.08}}
\put(15.5,3.75){\circle{0.08}}
\put(16,3.5){\circle{0.08}}
\put(16.5,3.25){\circle{0.08}}
\put(17,3){\circle{0.08}}
\put(17.5,2.75){\circle{0.08}}
\put(18,2.5){\circle{0.08}}

\put(18.5,2.25){\circle{0.08}}
\put(19,2){\circle{0.08}}
\put(19.5,1.75){\circle{0.08}}
\put(20,1.5){\circle{0.08}}
\put(20.5,1.25){\circle{0.08}}
\put(21,1){\circle{0.08}}
\put(21.5,0.75){\circle{0.08}}
\put(22,0.5){\circle{0.08}}
\put(22.5,0.25){\circle{0.08}}
\put(23,0){\circle{0.08}}
\put(23,0){\circle{0.12}}

\put(23.5,0.25){\circle{0.08}}
\put(24,0.5){\circle{0.08}}
\put(24.5,0.75){\circle{0.08}}
\put(25,1){\circle{0.08}}
\put(25.5,1.25){\circle{0.08}}
\put(26,1.5){\circle{0.08}}
\put(26.5,1.75){\circle{0.08}}
\put(27,2){\circle{0.08}}
\put(27,2){\circle{0.12}}

\put(15,0){\vector(1,0){13.5}}
\put(15,0){\vector(0,1){5.5}}

\put(15,0){\circle{0.12}}
\put(15,0){\circle{0.08}}

\put(20,1){\circle{0.12}}
\put(20,1){\circle{0.08}}

\put(15,4){\circle{0.12}}
\put(15,4){\circle{0.08}}

\put(27,4.5){\circle{0.12}}
\put(27,4.5){\circle{0.08}}

{\small \put(36,3){$k$}}

\setlength{\unitlength}{2mm}

\linethickness{0.15mm}

\put(42,0){\line(0,1){4.7}}   

\thicklines

\put(30,0){\line(1,0){4}}   

\put(34,0){\line(2,1){8}}   

\linethickness{0.15mm}

\put(34,0){\circle{0.08}}
\put(34,0){\circle{0.12}}

\put(42,4){\circle{0.12}}
\put(42,4){\circle{0.08}}

\put(30,0){\circle{0.12}}
\put(30,0){\circle{0.08}}
\put(30.5,0.05){\circle{0.08}}
\put(31,0.1){\circle{0.08}}
\put(31.5,0.15){\circle{0.08}}
\put(32,0.2){\circle{0.08}}
\put(32.5,0.25){\circle{0.08}}
\put(33,0.3){\circle{0.08}}
\put(33.5,0.35){\circle{0.08}}
\put(34,0.4){\circle{0.08}}
\put(34.5,0.45){\circle{0.08}}
\put(35,0.5){\circle{0.08}}
\put(35.5,0.55){\circle{0.08}}
\put(36,0.6){\circle{0.08}}
\put(36.5,0.65){\circle{0.08}}
\put(37,0.7){\circle{0.08}}
\put(37.5,0.75){\circle{0.08}}
\put(38,0.8){\circle{0.08}}
\put(38,0.8){\circle{0.12}}
\put(38.5,1){\circle{0.08}}
\put(39,1.2){\circle{0.08}}
\put(39.5,1.4){\circle{0.08}}
\put(40,1.6){\circle{0.08}}
\put(40.5,1.8){\circle{0.08}}
\put(41,2){\circle{0.08}}
\put(41.5,2.2){\circle{0.08}}
\put(42,2.4){\circle{0.08}}
\put(42,2.4){\circle{0.12}}

\put(30,0){\vector(1,0){13.5}}
\put(30,0){\vector(0,1){5.5}}

{\small \put(51,3){$l$}}

\setlength{\unitlength}{2mm}

\linethickness{0.15mm}

\put(57,0){\line(0,1){4.5}}   

\thicklines
\put(45,0){\line(1,0){4}}   

\put(49,0){\line(2,1){8}}   

\linethickness{0.15mm}

\put(45,4){\circle{0.12}}
\put(45,4){\circle{0.08}}
\put(45.5,3.75){\circle{0.08}}
\put(46,3.5){\circle{0.08}}
\put(46.5,3.25){\circle{0.08}}
\put(47,3){\circle{0.08}}
\put(47.5,2.75){\circle{0.08}}
\put(48,2.5){\circle{0.08}}

\put(48.5,2.25){\circle{0.08}}
\put(49,2){\circle{0.08}}
\put(49.5,1.75){\circle{0.08}}
\put(50,1.5){\circle{0.08}}
\put(50.5,1.25){\circle{0.08}}
\put(51,1){\circle{0.08}}
\put(51.5,0.75){\circle{0.08}}

\put(52,0.5){\circle{0.12}}
\put(52,0.5){\circle{0.08}}
\put(52.5,0.45){\circle{0.08}}
\put(53,0.4){\circle{0.08}}

\put(53.5,0.35){\circle{0.08}}
\put(54,0.3){\circle{0.08}}
\put(54.5,0.25){\circle{0.08}}
\put(55,0.2){\circle{0.08}}
\put(55.5,0.15){\circle{0.08}}
\put(56,0.1){\circle{0.08}}
\put(56.5,0.05){\circle{0.08}}
\put(57,0){\circle{0.08}}
\put(57,0){\circle{0.12}}

\put(45,0){\vector(1,0){13.5}}
\put(45,0){\vector(0,1){5.5}}

\put(45,0){\circle{0.12}}
\put(45,0){\circle{0.08}}

\put(49,0){\circle{0.12}}
\put(49,0){\circle{0.08}}

\put(45,4){\circle{0.12}}
\put(45,4){\circle{0.08}}

\put(57,4){\circle{0.12}}
\put(57,4){\circle{0.08}}

\end{picture}
\end{center}

\begin{center}

\setlength{\unitlength}{2mm}
\begin{picture}(58,7)

{\small \put(21,3){$m$}}
\linethickness{0.15mm}

\put(27,0){\line(0,1){5}}   
\thicklines

\put(15,0){\line(5,1){6}}   

\put(21,1.2){\line(2,1){6}}   
\linethickness{0.15mm}

\put(21,1.2){\circle{0.08}}
\put(21,1.2){\circle{0.12}}
\put(27,4.2){\circle{0.12}}

\put(15,0){\circle{0.12}}
\put(15,0){\circle{0.08}}
\put(15.5,0.05){\circle{0.08}}
\put(16,0.1){\circle{0.08}}
\put(16.5,0.15){\circle{0.08}}
\put(17,0.2){\circle{0.08}}
\put(17.5,0.25){\circle{0.08}}
\put(18,0.3){\circle{0.08}}
\put(18.5,0.35){\circle{0.08}}
\put(19,0.4){\circle{0.08}}
\put(19.5,0.45){\circle{0.08}}
\put(20,0.5){\circle{0.08}}
\put(20.5,0.55){\circle{0.08}}
\put(21,0.6){\circle{0.08}}
\put(21.5,0.65){\circle{0.08}}
\put(22,0.7){\circle{0.08}}
\put(22.5,0.75){\circle{0.08}}
\put(23,0.8){\circle{0.08}}
\put(23,0.8){\circle{0.12}}
\put(23.5,1){\circle{0.08}}
\put(24,1.2){\circle{0.08}}
\put(24.5,1.4){\circle{0.08}}
\put(25,1.6){\circle{0.08}}
\put(25.5,1.8){\circle{0.08}}
\put(26,2){\circle{0.08}}
\put(26.5,2.2){\circle{0.08}}
\put(27,2.4){\circle{0.08}}
\put(27,2.4){\circle{0.12}}

\put(15,0){\vector(1,0){13.5}}
\put(15,0){\vector(0,1){5.5}}

{\small \put(36,3){$n$}}

\setlength{\unitlength}{2mm}

\linethickness{0.15mm}

\put(42,0){\line(0,1){4}}   
\thicklines

\put(30,4){\line(2,-1){6}}   
\put(36,1){\line(6,-1){6}}   
\linethickness{0.15mm}

\put(30,4){\circle{0.08}}
\put(30,4){\circle{0.12}}
\put(36,1){\circle{0.08}}
\put(36,1){\circle{0.12}}
\put(42,0){\circle{0.08}}
\put(42,0){\circle{0.12}}

\put(30,0){\circle{0.12}}
\put(30,0){\circle{0.08}}
\put(30.5,0.05){\circle{0.08}}
\put(31,0.1){\circle{0.08}}
\put(31.5,0.15){\circle{0.08}}
\put(32,0.2){\circle{0.08}}
\put(32.5,0.25){\circle{0.08}}
\put(33,0.3){\circle{0.08}}
\put(33.5,0.35){\circle{0.08}}
\put(34,0.4){\circle{0.08}}
\put(34.5,0.45){\circle{0.08}}
\put(35,0.5){\circle{0.08}}
\put(35.5,0.55){\circle{0.08}}
\put(36,0.6){\circle{0.08}}
\put(36.5,0.65){\circle{0.08}}
\put(37,0.7){\circle{0.08}}
\put(37.5,0.75){\circle{0.08}}
\put(38,0.8){\circle{0.08}}
\put(38.5,0.85){\circle{0.08}}
\put(39,0.9){\circle{0.08}}
\put(39,0.9){\circle{0.12}}
\put(39.5,1.2){\circle{0.08}}
\put(40,1.6){\circle{0.08}}
\put(40.5,2){\circle{0.08}}
\put(41,2.4){\circle{0.08}}
\put(41.5,2.8){\circle{0.08}}
\put(42,3.2){\circle{0.08}}
\put(42,3.2){\circle{0.12}}

\put(30,0){\vector(1,0){13.5}}
\put(30,0){\vector(0,1){5.5}}

\end{picture}
\end{center}
\medskip

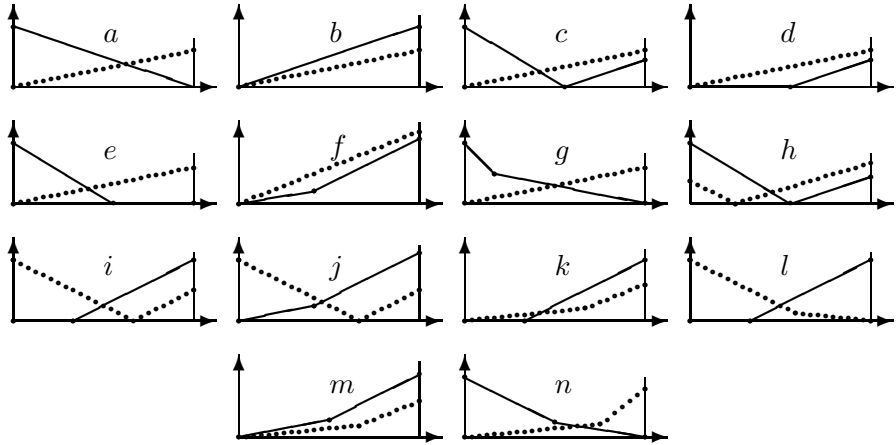
{\em {\small  {\bf Figure 1.}  Relevant combinations of shapes of the Newton polygons of $f$ with respect to a pair of prime numbers, each one having at most two edges.
}
}
\medskip

For most of the cases displayed in Figure 1, the corresponding results that we will consider will be special cases of Theorem A, and will consist  of a main theorem, followed by a corollary depending on some integer parameters $k_{1}, k_{2}, \dots $, and also by the simplest corollary, obtained in general by letting all or some of these parameters to be equal to 1. Some of the 14 cases in Figure 1 will have even more corresponding results, and this will be the case when the information on the divisibility of the leading coefficient of $f$, or of its free term, with respect to one of the primes will be irrelevant. Moreover, some immediate applications of Theorem B that use arbitrarily many Newton polygons
having one or two edges will be also provided.

The reader might lose his patience reading the statements of the irreducibility criteria that are immediate applications of Theorems A and B, since the 
number of such applications is quite discouraging large, even for Newton polygons that consist of only two edges. The reason is that such statements may look very much alike, and it might
be difficult to spot some significant features that make a particular result easier to remember. Despite this, some of these irreducibility conditions are sufficiently simple, or have enough symmetry that allow them to deserve further attention. That's why we advise the reader who is not familiar with Newton polygon method to focus on the simplest corollaries of Theorems A and B that correspond to the 14 cases plotted in Figure 1.
We might gain further insight on how ``wild" the irreducibility problem is, by imagining how difficult would be to describe all the relevant combinations of shapes of only two Newton polygons having each one arbitrarily many edges.

The idea to simultaneously use information on the Newton polygons of $f$ with respect to a family of prime numbers
proves to be quite rewarding. Indeed, if the irreducibility criteria of Sch\" onemann-Eisenstein and Dumas might seem to be in some sense some ``isolated" results of this type, the reader will see that they may be included as the simplest criteria in an infinite family of irreducibility criteria expressed by means of some divisibility conditions for the coefficients of a polynomial with respect to arbitrarily many prime numbers.

We will include here as an application of Theorem B the following result, which is also a special case of Proposition 3.10 in \cite{Mott}, result that refers to the case when each one of the Newton polygons of $f$ with respect to $p_{1},\dots ,p_{k}$ consists of a single edge, having either positive, or negative slope, covering both cases plotted in Figure 1.a and Figure 1.b.


\begin{theorem}\label{teorema1} 
Let $f(X)=a_{0}+a_{1}X+\cdots +a_{n}X^{n}\in \mathbb{Z}[X]$, $a_{0}a_{n}\neq 0$, let $k\geq 2$ and let
$p_{1},\dots ,p_{k}$ be pair-wise distinct prime numbers. Assume that for each $i=1,\dots ,k$ we have
\[
\nu _{p_{i}}(a_{j})\geq \frac{n-j}{n}\cdot \nu _{p_{i}}(a_{0})+\frac{j}{n}\cdot \nu _{p_{i}}(a_{n})\qquad  for\ j=1,\dots ,n-1, 
\]
where exactly one of the integers $\nu _{p_{i}}(a_{0})$ and $\nu _{p_{i}}(a_{n})$ is zero, and the non-zero one is denoted by $\alpha _{i}$.
If $\gcd(\alpha_{1}, n),\dots ,\gcd(\alpha _{k}, n)$ are relatively prime, then $f$ is irreducible over $\mathbb{Q}$.
\end{theorem}
Note that the condition that $\gcd(\alpha_{1}, n),\dots ,\gcd(\alpha _{k}, n)$ are relatively prime is satisfied if, for instance, 
$\alpha_{1},\dots ,\alpha _{k}$ are relatively prime (not necessarily pair-wise relatively prime), but this is obviously a more restrictive condition on the $\alpha _{i}$'s.

In particular, we obtain the following corollary.
\begin{corollary}
\label{corolarul1}
Let $f(X)=a_{0}+a_{1}X+\cdots +a_{n}X^{n}\in\mathbb{Z}[X]$, $a_{0}a_{n}\neq 0$, let $k\geq 2$, let
$p_{1},\dots ,p_{k}$ be pair-wise distinct prime numbers, and let $m_{1},\dots ,m_{k}$ be positive integers. Assume that for each $j=1,\dots ,k$ we have either
\smallskip

i) $p_{j}^{m_{j}}\mid a_{i}$ \ for \ $i=1,\dots ,n$, \  $p_{j}^{m_{j}+1}\nmid a_{n}$, \ and \ $p_{j}\nmid a_{0}$, 

\noindent or
 
ii) $p_{j}^{m_{j}}\mid a_{i}$ \ for \ $i=0,\dots ,n-1$, \  $p_{j}^{m_{j}+1}\nmid a_{0}$, \ and \ $p_{j}\nmid a_{n}$.
\smallskip

\noindent If $\gcd(m_{1},n),\dots ,\gcd(m_{k},n)$ are relatively prime, then $f$ is irreducible over $\mathbb{Q}$.
\end{corollary}

The following result refers to the case when the Newton polygon of $f$ with respect to a prime number $p$ consists of two segments whose slopes have different sign, while
the Newton polygon of $f$ with respect to another prime number, say $q$, consists of a single edge with positive slope, as in Figure 1.c.

\begin{theorem}\label{teorema8}
Let $f(X)=a_{0}+a_{1}X+\cdots +a_{n}X^{n}\in \mathbb{Z}[X]$, $a_{0}a_{n}\neq 0$. If
there exist two distinct prime numbers $p$ and $q$, and an index $j\in \{1,\dots ,n-1\} $ which is not a multiple of $n/\gcd(\nu _{q}(a_{n}),n)$, such that
\smallskip

i) \ $\nu _{p}(a_{0})>0$, \ \ $\frac{\nu _{p}(a_{i})}{j-i}>\frac{\nu _{p}(a_{0})}{j}$ \ \ for \ $0<i<j$,  \ \ $\nu _{p}(a_{j})=0$,

\quad \ $\frac{\nu _{p}(a_{i})}{i-j}>\frac{\nu _{p}(a_{n})}{n-j}$ \ \ for \ $j<i<n$, \ \ and \ \ $\nu _{p}(a_{n})>0$,

ii) $\nu _{q}(a_{0})=0$, \ \ $\frac{\nu _{q}(a_{i})}{i}\geq\frac{\nu _{q}(a_{n})}{n}$ \ \ for \ $i\geq 1$, \ \ and \ \ $\nu _{q}(a_{n})>0$,

iii) $\gcd(\nu _{p}(a_{0}),j)=\gcd(\nu _{p}(a_{n}),n-j)=1$,
\smallskip

\noindent then $f$ is irreducible over $\mathbb{Q}$.
\end{theorem}
In particular, we obtain the following corollary.

\begin{corollary}
\label{corolarul12}
Let $f(X)=a_{0}+a_{1}X+\cdots +a_{n}X^{n}\in \mathbb{Z}[X]$, $a_{0}a_{n}\neq 0$. If there exist two distinct prime numbers $p$ and $q$, three positive integers $k_{1},k_{2},k_{3}$, and an index $j\in \{1,\dots ,n-1\} $ which is not a multiple of $n/\gcd(k_{3},n)$, such that
\smallskip

i) $p^{k_{1}}\mid a_{i}$\quad for \ $i<j$,\quad  $p^{k_{1}+1}\nmid a_{0}$,\quad  $p\nmid a_{j}$,\quad  $p^{k_{2}}\mid a_{i}$\quad for \ $i>j$,\quad  $p^{k_{2}+1}\nmid a_{n}$

ii) $q\nmid a_{0}$,\quad  $q^{k_{3}}\mid a_{i}$\quad for \ $i>0$,\quad  $q^{k_{3}+1}\nmid a_{n}$,

iii)  $\gcd(k_{1},j)=\gcd(k_{2},n-j)=1$,
\smallskip

\noindent then $f$ is irreducible over $\mathbb{Q}$.
\end{corollary}
The simplest irreducibility conditions of this kind are obtained by taking $k_{1}=k_{2}=1$:
\begin{corollary}
\label{corolarul13}
Let $f(X)=a_{0}+a_{1}X+\cdots +a_{n}X^{n}\in \mathbb{Z}[X]$, $a_{0}a_{n}\neq 0$. If there exist two distinct prime numbers $p$ and $q$, a positive integer $k$, and an index $j\in \{1,\dots ,n-1\} $ which is not a multiple of $n/\gcd(k,n)$, such that
\smallskip

i) $p\mid a_{i}$\quad  for \  $i\neq j$,\quad $p\nmid a_{j}$,\quad $p^{2}\nmid a_{0}$, \quad $p^{2}\nmid a_{n}$,

ii) $q\nmid a_{0}$,\quad $q^{k}\mid a_{i}$\quad for \ $i=1,\dots ,n$,\quad  $q^{k+1}\nmid a_{n}$,

\smallskip

\noindent then $f$ is irreducible over $\mathbb{Q}$.
\end{corollary}

A special case of Theorem \ref{teorema8}, namely the case when $j=1$, deserves further attention, since it requires no information at all on the divisibility of $a_{0}$ by $p$. 

\begin{theorem}\label{teorema9}
Let $f(X)=a_{0}+a_{1}X+\cdots +a_{n}X^{n}\in \mathbb{Z}[X]$, $a_{0}a_{n}\neq 0$. If
there exist two distinct prime numbers $p$ and $q$ such that
\smallskip

i) \ $\nu _{p}(a_{1})=0$, \ \ $\nu _{p}(a_{n})>0$, \ \  and \ \  $\frac{\nu _{p}(a_{i})}{i-1}>\frac{\nu _{p}(a_{n})}{n-1}$\quad for \ $1<i<n$,

ii) \thinspace $\nu _{q}(a_{0})=0$, \ \ $\frac{\nu _{q}(a_{i})}{i}\geq\frac{\nu _{q}(a_{n})}{n}$ \ \ for \ $i\geq 1$, \ \ and \ \ $n\nmid \nu _{q}(a_{n})$,

iii) $\gcd(\nu _{p}(a_{n}),n-1)=1$,
\smallskip

\noindent then $f$ is irreducible over $\mathbb{Q}$.
\end{theorem}
In particular, we obtain the following corollary.

\begin{corollary}
\label{corolarul14}
Let $f(X)=a_{0}+a_{1}X+\cdots +a_{n}X^{n}\in \mathbb{Z}[X]$, $a_{0}a_{n}\neq 0$. If there exist two distinct prime numbers $p$ and $q$ and two positive integers $k_{1},k_{2}$ such that
\smallskip

i) \ \thinspace $p\nmid a_{1}$,\quad  $p^{k_{1}}\mid a_{i}$\quad for \ $i>1$,\quad  $p^{k_{1}+1}\nmid a_{n}$,

ii) \thinspace $q\nmid a_{0}$,\quad  $q^{k_{2}}\mid a_{i}$\quad for \ $i>0$,\quad  $q^{k_{2}+1}\nmid a_{n}$,\quad and\quad $n\nmid k_{2}$,

iii)  $\gcd(k_{1},n-1)=1$,
\smallskip

\noindent then $f$ is irreducible over $\mathbb{Q}$.
\end{corollary}
Here the simplest irreducibility conditions of this kind are obtained by taking $k_{1}=1$:
\begin{corollary}
\label{corolarul15}
Let $f(X)=a_{0}+a_{1}X+\cdots +a_{n}X^{n}\in \mathbb{Z}[X]$, $a_{0}a_{n}\neq 0$. If there exist two distinct prime numbers $p$ and $q$ and a positive integer $k$ such that
\smallskip

i) \ $p\nmid a_{1}$,\quad  $p\mid a_{i}$\quad for \ $i>1$,\quad  $p^{2}\nmid a_{n}$,

ii) $q\nmid a_{0}$,\quad  $q^{k}\mid a_{i}$\quad for \ $i>0$,\quad  $q^{k+1}\nmid a_{n}$,\quad and\quad $n\nmid k$,

\smallskip

\noindent then $f$ is irreducible over $\mathbb{Q}$.
\end{corollary}

Another special case of Theorem \ref{teorema8}, when $j=n-1$, deserves further attention too, since it will require no information at all on the divisibility of $a_{n}$ by $p$. In this case we have the following results.

\begin{theorem}\label{teorema10}
Let $f(X)=a_{0}+a_{1}X+\cdots +a_{n}X^{n}\in \mathbb{Z}[X]$, $a_{0}a_{n}\neq 0$. If
there exist two distinct prime numbers $p$ and $q$ such that
\smallskip

i) \ $\nu _{p}(a_{0})>0$, \ \ $\frac{\nu _{p}(a_{i})}{n-1-i}>\frac{\nu _{p}(a_{0})}{n-1}$ \ \ for \ $0<i<n-1$, \ \ and \ \ $\nu _{p}(a_{n-1})=0$,

ii) $\nu _{q}(a_{0})=0$, \ \ $\frac{\nu _{q}(a_{i})}{i}\geq\frac{\nu _{q}(a_{n})}{n}$ \ \ for \ $i\geq 1$, \ \ and \ \ $n\nmid \nu _{q}(a_{n})$,

iii) $\gcd(\nu _{p}(a_{0}),n-1)=1$,
\smallskip

\noindent then $f$ is irreducible over $\mathbb{Q}$.
\end{theorem}

\begin{corollary}
\label{corolarul16}
Let $f(X)=a_{0}+a_{1}X+\cdots +a_{n}X^{n}\in \mathbb{Z}[X]$, $a_{0}a_{n}\neq 0$. If there exist two distinct prime numbers $p$ and $q$, and two positive integers $k_{1},k_{2}$ such that
\smallskip

i) \ \thinspace $p^{k_{1}}\mid a_{i}$ \ \ for \ $i<n-1$, \ \  $p^{k_{1}+1}\nmid a_{0}$, \ \ and \ \  $p\nmid a_{n-1}$,

ii) \thinspace $q\nmid a_{0}$, \ \  $q^{k_{2}}\mid a_{i}$ \ \ for \ $i>0$, \ \  $q^{k_{2}+1}\nmid a_{n}$, \ \ and \ \ $n\nmid k_{2}$,

iii)  $\gcd(k_{1},n-1)=1$,
\smallskip

\noindent then $f$ is irreducible over $\mathbb{Q}$.
\end{corollary}
Again, the simplest irreducibility conditions of this kind are obtained by taking $k_{1}=1$:
\begin{corollary}
\label{corolarul17}
Let $f(X)=a_{0}+a_{1}X+\cdots +a_{n}X^{n}\in \mathbb{Z}[X]$, $a_{0}a_{n}\neq 0$. If there exist two distinct prime numbers $p$ and $q$ and a positive integer $k$ such that
\smallskip

i) \ $p\mid a_{i}$ \ \  for \  $i<n-1$, \ \ $p\nmid a_{n-1}$, \ \ $p^{2}\nmid a_{0}$,

ii) $q\nmid a_{0}$, \ \ $q^{k}\mid a_{i}$ \ \ for \ $i>0$, \ \  $q^{k+1}\nmid a_{n}$, \ \ and \ \ $n\nmid k$,

\smallskip

\noindent then $f$ is irreducible over $\mathbb{Q}$.
\end{corollary}

We will consider now the case when the Newton polygon of $f$ with respect to $q$
consists of a single edge with positive slope, while the Newton polygon of $f$ with respect to $p$ consists of two edges, one of which lies on the $x$-axis, while the other has positive slope, as in Figure 1.d. Our first result in this case is the following.

\begin{theorem}\label{teorema11}
Let $f(X)=a_{0}+a_{1}X+\cdots +a_{n}X^{n}\in \mathbb{Z}[X]$, $a_{0}a_{n}\neq 0$. If
there exist two distinct prime numbers $p$ and $q$ and an index $j<\frac{n}{ \gcd(\nu _{q}(a_{n}),n) }$ such that
\smallskip

i) \ \thinspace $\nu _{p}(a_{i})=0$ \ \ for \ $i\leq j$, \ \ $\frac{\nu _{p}(a_{i})}{i-j}>\frac{\nu _{p}(a_{n})}{n-j}$ \ \ for \ $j<i<n$, \ \ and \ \ $\nu _{p}(a_{n})>0$,

ii) \ $\nu _{q}(a_{0})=0$, \ \  $\frac{\nu _{q}(a_{i})}{i}\geq\frac{\nu _{q}(a_{n})}{n}$ \ \ for \ $i>0$, \ \ and \ \ $\nu _{q}(a_{n})>0$,

iii) $\gcd(\nu _{p}(a_{n}),n-j)=1$,
\smallskip

\noindent then $f$ is irreducible over $\mathbb{Q}$.
\end{theorem}
In particular, we have the following corollary.
\begin{corollary}\label{corolarul18}
Let $f(X)=a_{0}+a_{1}X+\cdots +a_{n}X^{n}\in \mathbb{Z}[X]$, $a_{0}a_{n}\neq 0$. If
there exist two distinct prime numbers $p$ and $q$, two positive integers $k_{1},k_{2}$, and an index $j<\frac{n}{ \gcd(k_{2},n) }$ such that
\smallskip

i) \ \thinspace  $p\nmid a_{i}$ \ \ for \ $i\leq j$, \ \ $p^{k_{1}}\mid a_{i}$ \ \ for \ $i>j$, \ \ and \ \  $p^{k_{1}+1}\nmid a_{n}$,

ii) \ $q\nmid a_{0}$, \ \  $q^{k_{2}}\mid a_{i}$ \ \ for \ $i>0$, \ \ and \ \ $q^{k_{2}+1}\nmid a_{n}$,

iii) $\gcd(k_{1},n-j)=1$,
\smallskip

\noindent then $f$ is irreducible over $\mathbb{Q}$.
\end{corollary}
The simplest conditions of this type are obtained for $k_{1}=1$.
\begin{corollary}\label{corolarul19}
Let $f(X)=a_{0}+a_{1}X+\cdots +a_{n}X^{n}\in \mathbb{Z}[X]$, $a_{0}a_{n}\neq 0$. If
there exist two distinct prime numbers $p$ and $q$, a positive integer $k$, and an index $j<\frac{n}{ \gcd(k,n) }$ such that
\smallskip

i) \ $p\nmid a_{i}$ \ \ for \ $i\leq j$, \ \ $p\mid a_{i}$ \ \ for \ $i>j$, \ \ and \ \  $p^{2}\nmid a_{n}$,

ii) $q\nmid a_{0}$, \ \  $q^{k}\mid a_{i}$ \ \ for \ $i>0$, \ \ and \ \ $q^{k+1}\nmid a_{n}$,
\smallskip

\noindent then $f$ is irreducible over $\mathbb{Q}$.
\end{corollary}

The following case is the one when the Newton polygon of $f$ with respect to $q$
consists of a single edge having positive slope, while the Newton polygon of $f$ with respect to $p$ consists of two edges, one of which lies on the $x$-axis, and the other has negative slope, as in Figure 1.e. 
In this case we have the following results.

\begin{theorem}\label{teorema12}
Let $f(X)=a_{0}+a_{1}X+\cdots +a_{n}X^{n}\in \mathbb{Z}[X]$, $a_{0}a_{n}\neq 0$. If
there exist two distinct prime numbers $p$ and $q$ and an index $j>n-\frac{n}{ \gcd(\nu _{q}(a_{n}),n) }$ such that
\smallskip

i) \ \thinspace $\nu _{p}(a_{0})>0$, \ \ $\frac{\nu _{p}(a_{i})}{j-i}>\frac{\nu _{p}(a_{0})}{j}$ \ \ for \ $0<i<j$, \ \ and \ \ $\nu _{p}(a_{i})=0$ \ \ for \ $i\geq j$,

ii) \ $\nu _{q}(a_{0})=0$, \ \ $\frac{\nu _{q}(a_{i})}{i}\geq\frac{\nu _{q}(a_{n})}{n}$ \ \ for \ $i>0$, \ \ and \ \ $\nu _{q}(a_{n})>0$,

iii) $\gcd(\nu _{p}(a_{0}),j)=1$,
\smallskip

\noindent then $f$ is irreducible over $\mathbb{Q}$.
\end{theorem}
\begin{corollary}\label{corolarul20}
Let $f(X)=a_{0}+a_{1}X+\cdots +a_{n}X^{n}\in \mathbb{Z}[X]$. If
there exist two distinct prime numbers $p$ and $q$, two positive integers $k_{1},k_{2}$, and an index $j>n-\frac{n}{ \gcd(k_{2},n) }$ such that
\smallskip

i) \ \thinspace $p\nmid a_{i}$ \ \ for \ $i\geq j$, \ \ $p^{k_{1}}\mid a_{i}$ \ \ for \ $i<j$, \ \ and \ \  $p^{k_{1}+1}\nmid a_{0}$,

ii) \ $q\nmid a_{0}$, \ \  $q^{k_{2}}\mid a_{i}$ \ \ for \ $i>0$, \ \ and \ \ $q^{k_{2}+1}\nmid a_{n}$,

iii) $\gcd(k_{1},j)=1$,
\smallskip

\noindent then $f$ is irreducible over $\mathbb{Q}$.
\end{corollary}
\begin{corollary}\label{corolarul21}
Let $f(X)=a_{0}+a_{1}X+\cdots +a_{n}X^{n}\in \mathbb{Z}[X]$. If
there exist two distinct prime numbers $p$ and $q$, a positive integer $k$, and an index $j>n-\frac{n}{ \gcd(k,n) }$  such that
\smallskip

i) \thinspace $p\nmid a_{i}$ \ \ for \ $i\geq j$, \ \ $p\mid a_{i}$ \ \ for \ $i<j$, \ \ and \ \  $p^{2}\nmid a_{0}$,

ii) $q\nmid a_{0}$, \ \  $q^{k}\mid a_{i}$ \ \ for \ $i>0$, \ \ and \ \ $q^{k+1}\nmid a_{n}$,
\smallskip

\noindent then $f$ is irreducible over $\mathbb{Q}$.
\end{corollary}

Our next result refers to the case when the Newton polygon of $f$ with respect to $p$ consists of two segments with positive different slopes, while the Newton polygon of $f$ with respect to $q$ consists of a 
single edge with positive slope too, as in Figure 1.f.

\begin{theorem}\label{teorema13}
Let $f(X)=a_{0}+a_{1}X+\cdots +a_{n}X^{n}\in \mathbb{Z}[X]$, $a_{0}a_{n}\neq 0$. If
there exist two distinct prime numbers $p$ and $q$ and an index $j<\frac{n}{ \gcd(\nu _{q}(a_{n}),n) }$ such that
\smallskip

i) \thinspace   $\nu _{p}(a_{0}) = 0$, \ \ $\frac{\nu _{p}(a_{i})}{i} > \frac{\nu _{p}(a_{j})}{j} $ \ \ for \  $0<i<j$, \ \  $\nu _{p}(a_{j})>0$,

\quad \thinspace $\nu _{p}(a_{i}) > \frac{n-i}{n-j} \nu _{p}(a_{j})+\frac{i-j}{n-j} \nu _{p}(a_{n})$\ \ for\ \  $j<i<n$, \ \ $\frac{\nu _{p}(a_{n})}{n}>\frac{\nu _{p}(a_{j})}{j}$, 
\smallskip

ii) \ $\nu _{q}(a_{0})=0$, \ \ $\frac{\nu _{q}(a_{i})}{i}\geq\frac{\nu _{q}(a_{n})}{n}$ \ \ for \ $i>0$, \ \ and \ \ $\nu _{q}(a_{n})>0$,

iii) $\gcd(\nu _{p}(a_{j}),j)=\gcd(\nu _{p}(a_{n})-\nu _{p}(a_{j}),n-j)=1$,
\smallskip

\noindent then $f$ is irreducible over $\mathbb{Q}$.
\end{theorem}

\begin{corollary}\label{corolarul22}
Let $f(X)=a_{0}+a_{1}X+\cdots +a_{n}X^{n}\in \mathbb{Z}[X]$. If
there exist two distinct prime numbers $p$ and $q$, three positive integers $k_{1}, k_{2}, k_{3}$ and an index $j<\frac{n}{ \gcd(k_{3},n) }$ such that
\smallskip

i) \ \thinspace   $p\nmid a_{0}$, \  $p^{k_{1}}\mid a_{i}$ \ \ for \  $0<i\leq j$, \  $p^{k_{1}+1}\nmid a_{j}$, \ \ $p^{k_{2}}\mid a_{i}$ \ \ for \ $i>j$, \  $p^{k_{2}+1}\nmid a_{n}$, \  $\frac{k_{1}}{j}<\frac{k_{2}}{n}$,

ii) \ $q\nmid a_{0}$, \   $q^{k_{3}}\mid a_{i}$ \  for \ $i>0$,\  and \  $q^{k_{3}+1}\nmid a_{n}$,

iii) $\gcd(k_{1},j)=\gcd(k_{2}-k_{1},n-j)=1$,
\smallskip

\noindent then $f$ is irreducible over $\mathbb{Q}$.
\end{corollary}

By letting $k_{1}=1$ and by writing $k_{1}$ for $k_{2}$ and $k_{2}$ for $k_{3}$, one obtains

\begin{corollary}\label{corolarul23}
Let $f(X)=a_{0}+a_{1}X+\cdots +a_{n}X^{n}\in \mathbb{Z}[X]$. If
there exist two distinct prime numbers $p$ and $q$, two positive integers $k_{1}, k_{2}$ and an index $j<\frac{n}{ \gcd(k_{2},n) }$ such that
\smallskip

i) \ \thinspace   $p\nmid a_{0}$, \  $p\mid a_{i}$ \  for \ $0<i\leq j$, \ $p^{2}\nmid a_{j}$, \  $p^{k_{1}}\mid a_{i}$ \  for \ $i>j$, \  $p^{k_{1}+1}\nmid a_{n}$, \  $k_{1}>\frac{n}{j}$,

ii) \ $q\nmid a_{0}$, \ \  $q^{k_{2}}\mid a_{i}$ \ \ for \ $i>0$, \ \ and \ \ $q^{k_{2}+1}\nmid a_{n}$,

iii) $\gcd(k_{1}-1,n-j)=1$,
\smallskip

\noindent then $f$ is irreducible over $\mathbb{Q}$.
\end{corollary}

The following result refers to the case when the Newton polygon of $f$ with respect to $p$ consists of two segments with negative different slopes, while the Newton polygon of $f$ with respect to $q$ consists of a single edge with positive slope, as in Figure 1.g.

\begin{theorem}\label{teorema14}
Let $f(X)=a_{0}+a_{1}X+\cdots +a_{n}X^{n}\in \mathbb{Z}[X]$, $a_{0}a_{n}\neq 0$. If
there exist two distinct prime numbers $p$ and $q$ and an index $j<\frac{n}{ \gcd(\nu _{q}(a_{n}),n) }$ such that
\smallskip

i) \ $\frac{\nu _{p}(a_{0})}{n}>\frac{\nu _{p}(a_{j})}{n-j}$, \ \ $\nu _{p}(a_{i}) > \frac{j-i}{j} \nu _{p}(a_{0})+\frac{i}{j} \nu _{p}(a_{j})$\ \ for\ \  $0<i<j$, \ \ $\nu _{p}(a_{j})>0$, 

\quad \thinspace \thinspace $\frac{\nu _{p}(a_{i})}{n-i}  > \frac{\nu _{p}(a_{j})}{n-j} $\ \ for\ \ $j<i<n$, \ \ and \ \ $\nu _{p}(a_{n}) = 0$,

\smallskip

ii) \ $\nu _{q}(a_{0})=0$, \ \ $\frac{\nu _{q}(a_{i})}{i}\geq \frac{\nu _{q}(a_{n})}{n}$ \ \ for \ $i>0$, \ \ and \ \ $\nu _{q}(a_{n})>0$,

iii) $\gcd(\nu _{p}(a_{j}),n-j)=\gcd(\nu _{p}(a_{0})-\nu _{p}(a_{j}),j)=1$,
\smallskip

\noindent then $f$ is irreducible over $\mathbb{Q}$.
\end{theorem}

\begin{corollary}\label{corolarul24}
Let $f(X)=a_{0}+a_{1}X+\cdots +a_{n}X^{n}\in \mathbb{Z}[X]$. If
there exist two distinct prime numbers $p$ and $q$, three positive integers $k_{1}, k_{2}, k_{3}$ and an index $j<\frac{n}{ \gcd(k_{3},n) }$ such that
\smallskip

i) \ \thinspace   $\frac{k_{1}}{n}>\frac{k_{2}}{n-j}$, \  $p^{k_{1}}\mid a_{i}$ \   for \  $i<j$, \ $p^{k_{1}+1}\nmid a_{0}$, \ $p^{k_{2}}\mid a_{i}$ \  for \  $j\leq i<n$, \  $p^{k_{2}+1}\nmid a_{j}$, \  $p\nmid a_{n}$, 

ii) \ $q\nmid a_{0}$, \   $q^{k_{3}}\mid a_{i}$ \  for \ $i>0$,\  and \  $q^{k_{3}+1}\nmid a_{n}$,

iii) $\gcd(k_{1}-k_{2},j)=\gcd(k_{2},n-j)=1$,
\smallskip

\noindent then $f$ is irreducible over $\mathbb{Q}$.
\end{corollary}

The simplest irreducibility conditions of this type are obtained by letting $k_{2}=1$ and by writing $k_{2}$ for $k_{3}$.

\begin{corollary}\label{corolarul25}
Let $f(X)=a_{0}+a_{1}X+\cdots +a_{n}X^{n}\in \mathbb{Z}[X]$. If
there exist two distinct prime numbers $p$ and $q$, two positive integers $k_{1}, k_{2}$ and an index $j<\frac{n}{ \gcd(k_{2},n) }$ such that
\smallskip

i) \ \thinspace   $k_{1}>\frac{n}{n-j}$, \  $p^{k_{1}}\mid a_{i}$ \  for \ $i<j$, \  $p^{k_{1}+1}\nmid a_{0}$, \  $p\mid a_{i}$ \  for \ $j\leq i<n$, \  $p^{2}\nmid a_{j}$, \  $p\nmid a_{n}$,

ii) \ $q\nmid a_{0}$, \ $q^{k_{2}}\mid a_{i}$ \  for \ $i>0$, \  and \  $q^{k_{2}+1}\nmid a_{n}$,

iii) $\gcd(k_{1}-1,j)=1$,
\smallskip

\noindent then $f$ is irreducible over $\mathbb{Q}$.
\end{corollary}

In the case when each one of the Newton polygons of $f$ with respect to two prime numbers $p$ and $q$ consists of two segments, whose slopes have different sign, as in Figure 1.h,
Theorem A leads us to the following irreducibility conditions.

\begin{theorem}\label{teorema2}
Let $f(X)=a_{0}+a_{1}X+\cdots +a_{n}X^{n}\in \mathbb{Z}[X]$, $a_{0}a_{n}\neq 0$. If there exist two distinct indices $j_{1},j_{2}\in \{1,\dots ,n-1\} $ such that $j_{1}+j_{2}\neq n$, 
and two distinct prime numbers $p$ and $q$ such that $\nu _{p}(a_{j_{1}})=\nu _{q}(a_{j_{2}})=0$, $\nu _{p}(a_{0})\nu _{q}(a_{0})\nu _{p}(a_{n})\nu _{q}(a_{n})\neq 0$, and
\smallskip

i) \ $\frac{\nu _{p}(a_{i})}{j_{1}-i}>\frac{\nu _{p}(a_{0})}{j_{1}}$ \ for \ $0<i<j_{1}$ \quad and \quad  $\frac{\nu _{p}(a_{i})}{i-j_{1}}>\frac{\nu _{p}(a_{n})}{n-j_{1}}$ \ for \ $j_{1}<i<n$,

ii)  $\frac{\nu _{q}(a_{i})}{j_{2}-i}>\frac{\nu _{q}(a_{0})}{j_{2}}$ \ for \ $0<i<j_{2}$ \quad and \quad  $\frac{\nu _{q}(a_{i})}{i-j_{2}}>\frac{\nu _{q}(a_{n})}{n-j_{2}}$ \ for \ $j_{2}<i<n$,

iii) $\gcd(\nu _{p}(a_{0}),j_{1})=\gcd(\nu _{p}(a_{n}),n-j_{1})=\gcd(\nu _{q}(a_{0}),j_{2})=\gcd(\nu _{q}(a_{n}),n-j_{2})=1$,
\smallskip

\noindent then $f$ is irreducible over $\mathbb{Q}$.
\end{theorem}
In particular, we obtain the following corollaries.
\begin{corollary}\label{corolarul2}
Let $f(X)=a_{0}+a_{1}X+\cdots +a_{n}X^{n}\in \mathbb{Z}[X]$, $a_{0}a_{n}\neq 0$. If there exist two distinct indices $j_{1},j_{2}\in \{1,\dots ,n-1\} $ such that $j_{1}+j_{2}\neq n$, two distinct prime numbers $p$ and $q$
and four positive integers $k_{1},k_{2},k_{3}, k_{4}$ such that $p\nmid a_{j_{1}}$, $q\nmid a_{j_{2}}$ and
\smallskip

i) \ $p^{k_{1}}\mid a_{i}$ \ for \ $i<j_{1}$, $p^{k_{1}+1}\nmid a_{0}$\quad  and\quad  $p^{k_{2}}\mid a_{i}$ \ for \ $i>j_{1}$, $p^{k_{2}+1}\nmid a_{n}$,

ii) \thinspace $q^{k_{3}}\mid a_{i}$ \ for \ $i<j_{2}$, $q^{k_{3}+1}\nmid a_{0}$\quad  and\quad  $q^{k_{4}}\mid a_{i}$ \ for \ $i>j_{2}$, $q^{k_{4}+1}\nmid a_{n}$,

iii) $\gcd(k_{1},j_{1})=\gcd(k_{2},n-j_{1})=\gcd(k_{3},j_{2})=\gcd(k_{4},n-j_{2})=1$,
\smallskip

\noindent then $f$ is irreducible over $\mathbb{Q}$.
\end{corollary}

Here the simplest conditions of this type are obtained by letting all the four parameters $k_{1}, k_{2}, k_{3}, k_{4}$ to be equal to $1$.

\begin{corollary}\label{corolarul3}
Let $f(X)=a_{0}+a_{1}X+\cdots +a_{n}X^{n}\in \mathbb{Z}[X]$, $a_{0}a_{n}\neq 0$. If there exist two distinct indices 
$j_{1},j_{2}\in \{1,\dots ,n-1\} $ such that $j_{1}+j_{2}\neq n$, and two distinct prime numbers $p$ and $q$ such that
\smallskip

i) \ $p\mid a_{i}$ \ for \ $i\neq j_{1}$, \ $p\nmid a_{j_{1}}$, \ $p^{2}\nmid a_{0}$ \ and \ $p^{2}\nmid a_{n}$,

ii) $q\mid a_{i}$ \ for \ $i\neq j_{2}$, \ $q\nmid a_{j_{2}}$, \ $q^{2}\nmid a_{0}$ \ and \ $q^{2}\nmid a_{n}$,

\smallskip

\noindent then $f$ is irreducible over $\mathbb{Q}$.
\end{corollary}

Here we notice that one can easily drop the condition $j_{1}+j_{2}\neq n$ by adding a condition similar to i) and ii) for a third prime, and a third corresponding index, as in the following result.

\begin{corollary}\label{corolarul4}
Let $f(X)=a_{0}+a_{1}X+\cdots +a_{n}X^{n}\in \mathbb{Z}[X]$. If
there exist three distinct indices $j_{1},j_{2},j_{3}\in \{1,\dots ,n-1\} $, and three distinct prime numbers $p$, $q$ and $r$ such that 
\smallskip

i) \ \thinspace $p\mid a_{i}$ \ for \ $i\neq j_{1}$, \ $p\nmid a_{j_{1}}$, \ $p^{2}\nmid a_{0}$ \ and \ $p^{2}\nmid a_{n}$,

ii) \thinspace $q\mid a_{i}$ \ for \ $i\neq j_{2}$, \ $q\nmid a_{j_{2}}$, \ $q^{2}\nmid a_{0}$ \ and \ $q^{2}\nmid a_{n}$,

iii) $r\mid a_{i}$ \ for \ $i\neq j_{3}$, \ $r\nmid a_{j_{3}}$, \ $r^{2}\nmid a_{0}$ \ and \ $r^{2}\nmid a_{n}$,

\smallskip

\noindent then $f$ is irreducible over $\mathbb{Q}$.
\end{corollary}
Obviously, in a similar way one may remove the condition $j_{1}+j_{2}\neq n$  in the statements of Theorem \ref{teorema2} and Corollary \ref{corolarul2} too, by adding an additional condition for a third prime number 
and the corresponding third index, as in Corollary \ref{corolarul4}. However, in practice it is easier to test a condition like $j_{1}+j_{2}\neq n$, rather than a number of divisibilities corresponding to a third prime.

We remark that when one of the indices $j_{1}$ and $j_{2}$ in Theorem \ref{teorema2} is equal to $1$, while the other is different from $0$, $1$, $n-1$, and $n$, then no information is required on the divisibility
of $a_{0}$ with respect to one of the two prime numbers $p$ and $q$. Indeed, in this case we obtain from Theorem \ref{teorema2} the following result.

\begin{theorem}\label{teorema3}
Let $f(X)=a_{0}+a_{1}X+\cdots +a_{n}X^{n}\in \mathbb{Z}[X]$, $a_{0}a_{n}\neq 0$. If
there exists an index $j\in \{2,\dots ,n-2\} $ and two distinct prime numbers $p$ and $q$ such that 
\smallskip

i) \ $\nu _{p}(a_{1})=\nu _{q}(a_{j})=0$, $\nu _{p}(a_{n})\nu _{q}(a_{0})\nu _{q}(a_{n})\neq 0$,

ii) \thinspace $\frac{\nu _{p}(a_{i})}{i-1}>\frac{\nu _{p}(a_{n})}{n-1}$ \ for \ $i=2,\dots ,n-1$,

iii) $\frac{\nu _{q}(a_{i})}{j-i}>\frac{\nu _{q}(a_{0})}{j}$ \ for \ $i=1,\dots ,j-1$ \ \ and \ \  $\frac{\nu _{q}(a_{i})}{i-j}>\frac{\nu _{q}(a_{n})}{n-j}$ \ for \ $i=j+1,\dots ,n-1$,

iv) $\gcd(\nu _{p}(a_{n}),n-1)=\gcd(\nu _{q}(a_{0}),j)=\gcd(\nu _{q}(a_{n}),n-j)=1$,
\smallskip

\noindent then $f$ is irreducible over $\mathbb{Q}$.
\end{theorem}
In particular, we obtain the following corollaries.
\begin{corollary}\label{corolarul5}
Let $f(X)=a_{0}+a_{1}X+\cdots +a_{n}X^{n}\in \mathbb{Z}[X]$, $a_{0}a_{n}\neq 0$. If
there exists an index $j\in \{2,\dots ,n-2\} $, two distinct prime numbers $p$ and $q$, and three positive integers $k_{1},k_{2},k_{3}$ such that 
\smallskip

i) \ $p\nmid a_{1}$, \ $p^{k_{1}}\mid a_{i}$ \ for \ $i>1$, \ $p^{k_{1}+1}\nmid a_{n}$,

ii) \thinspace $q^{k_{2}}\mid a_{i}$ \ for \ $i<j$, \ $q^{k_{2}+1}\nmid a_{0}$, \ $q\nmid a_{j}$, \ and \ $q^{k_{3}}\mid a_{i}$ \ for \ $i>j$, \ $q^{k_{3}+1}\nmid a_{n}$,

iii) $\gcd(k_{1},n-1)=\gcd(k_{2},j)=\gcd(k_{3},n-j)=1$,
\smallskip

\noindent then $f$ is irreducible over $\mathbb{Q}$.
\end{corollary}
By letting $k_{1}=k_{2}=k_{3}=1$, we obtain the following result.

\begin{corollary}\label{corolarul6}
Let $f(X)=a_{0}+a_{1}X+\cdots +a_{n}X^{n}\in \mathbb{Z}[X]$, $a_{0}a_{n}\neq 0$. If
there exists an index $j\in \{2,\dots ,n-2\} $ and two distinct prime numbers $p$ and $q$ such that 
\smallskip

i) \ $p\mid a_{i}$ \ for \ $i>1$, \ $p\nmid a_{1}$, \ and \ $p^{2}\nmid a_{n}$,

ii) $q\mid a_{i}$ \ for \ $i\neq j$, \ $q\nmid a_{j}$, \ $q^{2}\nmid a_{0}$ \ and \ $q^{2}\nmid a_{n}$,

\smallskip

\noindent then $f$ is irreducible over $\mathbb{Q}$.
\end{corollary}
One may easily prove similar results that will use no information on the divisibility
of $a_{n}$ with respect to one of the two prime numbers $p$ and $q$, by considering the reciprocal of $f$ instead of $f$.

Another consequence of Theorem B is the following result that refers to the case when each one of the Newton polygons of $f$ with respect to $p_{1},\dots ,p_{k}$ consists of exactly two edges, whose slopes have different sign.

\begin{theorem}\label{teorema4} 
Let $f(X)=a_{0}+a_{1}X+\cdots +a_{n}X^{n}\in \mathbb{Z}[X]$, $a_{0}a_{n}\neq 0$, let $k\geq 2$ and let
$p_{1},\dots ,p_{k}$ be pair-wise distinct prime numbers. Assume that there exist indices $j_{1},\dots ,j_{k}\in \{1,\dots ,n-1\} $ such that for each $i=1,\dots ,k$ we have $\nu _{p_{i}}(a_{j_{i}})=0$, \ $\nu _{p_{i}}(a_{0})\nu _{p_{i}}(a_{n})\neq 0$, and
\smallskip

$\frac{\nu _{p_{i}}(a_{l})}{j_{i}-l} \geq \frac{\nu _{p_{i}}(a_{0})}{j_{i}}$ \ for \ $l<j_{i}$, \ and \  $\frac{\nu _{p_{i}}(a_{l})}{l-j_{i}} \geq \frac{\nu _{p_{i}}(a_{n})}{n-j_{i}}$ \ for \ $l>j_{i}$.
\smallskip

\noindent Let $d_{p_{i}}=\gcd \left( \frac{j_{i}}{\gcd (\nu _{p_{i}}(a_{0}),j_{i})}, \frac{n-j_{i}}{\gcd (\nu _{p_{i}}(a_{n}),n-j_{i})}   \right) $, for $i=1,\dots ,k.$
If \ $\gcd(\frac{n}{d_{p_{1}}},\dots ,\frac{n}{d_{p_{k}}})=1$, then $f$ is irreducible over $\mathbb{Q}$.
\end{theorem}

The following results refer to the case when the Newton polygon of $f$ with respect to $p$ consists of an edge situated on the $x$-axis, followed by a segment with positive slope,
while the Newton polygon of $f$ with respect to $q$ consists of two segments whose slopes have different sign, as in Figure 1.i.

\begin{theorem}\label{teorema15} 
Let $f(X)=a_{0}+a_{1}X+\cdots +a_{n}X^{n}\in \mathbb{Z}[X]$, $a_{0}a_{n}\neq 0$. If there exist two indices $j_{1},j_{2}$ with $j_{1}<{\rm min}(j_{2},n-j_{2})$, 
and two distinct prime numbers $p$ and $q$ such that 
\smallskip

i) \ $\nu _{p}(a_{i})=0$ \ for \ $i\leq j_{1}$, \  $\frac{\nu _{p}(a_{i})}{i-j_{1}}>\frac{\nu _{p}(a_{n})}{n-j_{1}}$ \ for \ $j_{1}<i<n$, \ and \ $\nu _{p}(a_{n})>0$,

ii)  $\nu _{q}(a_{0})>0$, \ $\frac{\nu _{q}(a_{i})}{j_{2}-i}>\frac{\nu _{q}(a_{0})}{j_{2}}$ \ for \ $0<i<j_{2}$, \ $\nu _{q}(a_{j_{2}})=0$, 

\quad \ and \ $\frac{\nu _{q}(a_{i})}{i-j_{2}}>\frac{\nu _{q}(a_{n})}{n-j_{2}}$ \ for \ $j_{2}<i<n$, \
$\nu _{q}(a_{n})>0$,

iii) $\gcd(\nu _{p}(a_{n}),n-j_{1})=\gcd(\nu _{q}(a_{0}),j_{2})=\gcd(\nu _{q}(a_{n}),n-j_{2})=1$,
\smallskip

\noindent then $f$ is irreducible over $\mathbb{Q}$.
\end{theorem}
\begin{corollary}\label{corolarul26}
Let $f(X)=a_{0}+a_{1}X+\cdots +a_{n}X^{n}\in \mathbb{Z}[X]$, $a_{0}a_{n}\neq 0$. If
there exist two distinct prime numbers $p$ and $q$, two indices $j_{1},j_{2}$ with $j_{1}<{\rm min}(j_{2},n-j_{2})$, and three positive integers $k_{1}, k_{2}, k_{3}$ such that
\smallskip

i) \ \thinspace   $p\nmid a_{i}$ \ \ for \ $i\leq j_{1}$, \  $p^{k_{1}}\mid a_{i}$\  \ for\  $i>j_{1}$,\ \ $p^{k_{1}+1}\nmid a_{n}$, 

ii) \ $q^{k_{2}}\mid a_{i}$ \  for \ $i<j_{2}$,\  \ $q^{k_{2}+1}\nmid a_{0}$, \ $q\nmid a_{j_{2}}$, \ \ $q^{k_{3}}\mid a_{i}$ \  for \ $i>j_{2}$,\  \ $q^{k_{3}+1}\nmid a_{n}$, 

iii) $\gcd(k_{1},n-j_{1})=\gcd(k_{2},j_{2})=\gcd(k_{3},n-j_{2})=1$,
\smallskip

\noindent then $f$ is irreducible over $\mathbb{Q}$.
\end{corollary}

\begin{corollary}\label{corolarul27}
Let $f(X)=a_{0}+a_{1}X+\cdots +a_{n}X^{n}\in \mathbb{Z}[X]$, $a_{0}a_{n}\neq 0$. If
there exist two distinct prime numbers $p$ and $q$ and two indices $j_{1},j_{2}$ with $j_{1}<{\rm min}(j_{2},n-j_{2})$ such that
\smallskip

i) \ \thinspace   $p\nmid a_{i}$ \  for \ $i\leq j_{1}$, \ \ $p\mid a_{i}$\  \ for\  $i>j_{1}$,\ \ and \ \ $p^{2}\nmid a_{n}$, 

ii) \ $q\mid a_{i}$ \  for \ $i\neq j_{2}$, \ \ $q\nmid a_{j_{2}}$, \ \ $q^{2}\nmid a_{0}$,\ \ and \ \ $q^{2}\nmid a_{n}$, 

\smallskip

\noindent then $f$ is irreducible over $\mathbb{Q}$.
\end{corollary}

Our next result refers to the case when the Newton polygon of $f$ with respect to $p$ consists of two segments with positive different slopes, while the Newton polygon of $f$ with respect to $q$ consists of two segments whose slopes have different sign, as in Figure 1.j.

\begin{theorem}\label{teorema16}
Let $f(X)=a_{0}+a_{1}X+\cdots +a_{n}X^{n}\in\mathbb{Z}[X]$, $a_{0}a_{n}\neq 0$. If
there exist two distinct primes $p,q$ and two distinct indices $j_{1},j_{2}\in \{ 1,\dots ,n-1 \} $ with $j_{1}+j_{2}\neq n$ such that
\smallskip

i) \thinspace   $\nu _{p}(a_{0}) = 0$, \ \ $\frac{\nu _{p}(a_{i})}{i} > \frac{\nu _{p}(a_{j_{1}})}{j_{1}} $ \ \ for \  $0<i<j_{1}$, \ \  $\nu _{p}(a_{j_{1}})>0$,

\quad \thinspace $\nu _{p}(a_{i}) > \frac{n-i}{n-j_{1}} \nu _{p}(a_{j_{1}})+\frac{i-j_{1}}{n-j_{1}} \nu _{p}(a_{n})$\ \ for\ \  $j_{1}<i<n$, \ \ $\frac{\nu _{p}(a_{n})}{n}>\frac{\nu _{p}(a_{j_{1}})}{j_{1}}$, 

\smallskip

ii)  $\nu _{q}(a_{0})>0$, \ $\frac{\nu _{q}(a_{i})}{j_{2}-i}>\frac{\nu _{q}(a_{0})}{j_{2}}$ \ for \ $0<i<j_{2}$, \ $\nu _{q}(a_{j_{2}})=0$, 

\quad \ and \ $\frac{\nu _{q}(a_{i})}{i-j_{2}}>\frac{\nu _{q}(a_{n})}{n-j_{2}}$ \ for \ $j_{2}<i<n$, \
$\nu _{q}(a_{n})>0$,

iii) $\gcd(\nu _{p}(a_{j_{1}}),j_{1})=\gcd(\nu _{p}(a_{n})-\nu _{p}(a_{j_{1}}),n-j_{1})=1$,

\quad \ \ $\gcd(\nu _{q}(a_{0}),j_{2})=\gcd(\nu _{q}(a_{n}),n-j_{2})=1$,
\smallskip

\noindent then $f$ is irreducible over $\mathbb{Q}$.
\end{theorem}

In particular one obtains the following result.

\begin{corollary}\label{corolarul28}
Let $f(X)=a_{0}+a_{1}X+\cdots +a_{n}X^{n}\in \mathbb{Z}[X]$, $a_{0}a_{n}\neq 0$. If
there exist two distinct prime numbers $p$ and $q$, two distinct indices $j_{1},j_{2}\in \{ 1,\dots ,n-1 \} $ with $j_{1}+j_{2}\neq n$, and four positive integers $k_{1}, k_{2}, k_{3}, k_{4}$ such that
\smallskip

i) \ \thinspace   $p\nmid a_{0}$, \  $p^{k_{1}}\mid a_{i}$\  for\  $0<i\leq j_{1}$,\ $p^{k_{1}+1}\nmid a_{j_{1}}$, $p^{k_{2}}\mid a_{i}$\  for\  $i>j_{1}$, $p^{k_{2}+1}\nmid a_{n}$, $k_{1}<\frac{j_{1}}{n}k_{2}$,

ii) \ $q^{k_{3}}\mid a_{i}$ \  for \ $i<j_{2}$,\  \ $q^{k_{3}+1}\nmid a_{0}$, \ $q\nmid a_{j_{2}}$, \ \ $q^{k_{4}}\mid a_{i}$ \  for \ $i>j_{2}$,\  \ $q^{k_{4}+1}\nmid a_{n}$, 

iii) $\gcd(k_{1},j_{1})=\gcd(k_{2}-k_{1},n-j_{1})=\gcd(k_{3},j_{2})=\gcd(k_{4},n-j_{2})=1$,
\smallskip

\noindent then $f$ is irreducible over $\mathbb{Q}$.
\end{corollary}

Here the simplest irreducibility conditions of this type are obtained by letting $k_{1}=k_{3}=k_{4}=1$ and by writing $k$ for $k_{2}$.

\begin{corollary}\label{corolarul29}
Let $f(X)=a_{0}+a_{1}X+\cdots +a_{n}X^{n}\in \mathbb{Z}[X]$, $a_{0}a_{n}\neq 0$. If
there exist two distinct prime numbers $p$ and $q$, two distinct indices $j_{1},j_{2}\in \{ 1,\dots ,n-1 \} $ with $j_{1}+j_{2}\neq n$, and a positive integer $k>\frac{n}{j_{1}}$ such that
\smallskip

i) \ \thinspace   $p\nmid a_{0}$, \ $p\mid a_{i}$ \ for \  $0<i\leq j_{1}$, \ \ $p^{2}\nmid a_{j_{1}}$, \ \ $p^{k}\mid a_{i}$ \ for \  $i>j_{1}$, \ $p^{k+1}\nmid a_{n}$, 

ii) \ $q\mid a_{i}$ \  for \ $i\neq j_{2}$,\  \ $q\nmid a_{j_{2}}$, \ \ $q^{2}\mid a_{0}$ \  and \ $q^{2}\nmid a_{n}$, 

iii) $\gcd(k-1,n-j_{1})=1$,
\smallskip

\noindent then $f$ is irreducible over $\mathbb{Q}$.
\end{corollary}

The following results refer to the case when the Newton polygon of $f$ with respect to $p$ consists of an edge situated on the $x$-axis, followed by a segment with positive slope, while the Newton polygon of $f$ with respect to $q$ consists of two segments with positive different slopes, as in Figure 1.k.

\begin{theorem}\label{teorema17} 
Let $f(X)=a_{0}+a_{1}X+\cdots +a_{n}X^{n}\in \mathbb{Z}[X]$, $a_{0}a_{n}\neq 0$. If there exist two distinct prime numbers $p,q$
and two distinct indices $j_{1},j_{2}$ such that  $j_{1}<{\rm min}(j_{2},n-j_{2})$ and
\smallskip

i) \ $\nu _{p}(a_{i})=0$ \ for \ $i\leq j_{1}$, \  $\frac{\nu _{p}(a_{i})}{i-j_{1}}>\frac{\nu _{p}(a_{n})}{n-j_{1}}$ \ for \ $j_{1}<i<n$, \ and \ $\nu _{p}(a_{n})>0$,

ii) \thinspace   $\nu _{q}(a_{0}) = 0$, \ \ $\frac{\nu _{q}(a_{i})}{i} > \frac{\nu _{q}(a_{j_{2}})}{j_{2}} $ \ \ for \  $0<i<j_{2}$, \ \  $\nu _{q}(a_{j_{2}})>0$,

\quad \ $\nu _{q}(a_{i}) > \frac{n-i}{n-j_{2}} \nu _{q}(a_{j_{2}})+\frac{i-j_{2}}{n-j_{2}} \nu _{q}(a_{n})$\ \ for\ \  $j_{2}<i<n$, \ \ $\frac{\nu _{q}(a_{n})}{n}>\frac{\nu _{q}(a_{j_{2}})}{j_{2}}$, 

iii) $\gcd(\nu _{p}(a_{n}),n-j_{1})=gcd(\nu _{q}(a_{j_{2}}),j_{2})=\gcd(\nu _{q}(a_{n})-\nu _{q}(a_{j_{2}}),n-j_{2})=1$,

\smallskip

\noindent then $f$ is irreducible over $\mathbb{Q}$.
\end{theorem}

\begin{corollary}\label{corolarul30}
Let $f(X)=a_{0}+a_{1}X+\cdots +a_{n}X^{n}\in \mathbb{Z}[X]$, $a_{0}a_{n}\neq 0$. If there exist two distinct prime numbers $p,q$,
two distinct indices $j_{1},j_{2}$ with  $j_{1}<{\rm min}(j_{2},n-j_{2})$, and three positive integers $k_{1},k_{2},k_{3}$ such that $k_{2}<\frac{j_{2}}{n}k_{3}$ and
\smallskip

i) \thinspace $p\nmid a_{i}$ \ for \ $i\leq j_{1}$, \   $p^{k_{1}}\mid a_{i}$ \  for \  $i>j_{1}$, \ $p^{k_{1}+1}\nmid a_{n}$,

ii) $q\nmid a_{0}$, \   $q^{k_{2}}\mid a_{i}$ \ for\ $0<i\leq j_{2}$,  \ $q^{k_{2}+1}\nmid a_{j_{2}}$,
\ $q^{k_{3}}\mid a_{i}$ \  for \  $i>j_{2}$, \ $q^{k_{3}+1}\nmid a_{n}$,

iii) $\gcd(k_{1},n-j_{1})=\gcd(k_{2},j_{2})=\gcd(k_{3}-k_{2},n-j_{2})=1$,

\smallskip

\noindent then $f$ is irreducible over $\mathbb{Q}$.
\end{corollary}

By letting $k_{1}=k_{2}=1$, and writing $k$ for $k_{3}$ we obtain the following result.

\begin{corollary}\label{corolarul31}
Let $f(X)=a_{0}+a_{1}X+\cdots +a_{n}X^{n}\in \mathbb{Z}[X]$, $a_{0}a_{n}\neq 0$. If there exist two distinct prime numbers $p,q$,
two distinct indices $j_{1},j_{2}$ with $j_{1}<{\rm min}(j_{2},n-j_{2})$, and a positive integer $k>\frac{n}{j_{2}}$ such that
\smallskip

i) \thinspace $p\nmid a_{i}$ \ for \ $i\leq j_{1}$, \ $p\mid a_{i}$ \ for\ $i>j_{1}$,  \ $p^{2}\nmid a_{n}$,

ii) $q\nmid a_{0}$, \   $q\mid a_{i}$ \ for\ $0<i\leq j_{2}$,  \ $q^{2}\nmid a_{j_{2}}$,
 \ $q^{k}\mid a_{i}$ \ for\  $i>j_{2}$,  \ $q^{k+1}\nmid a_{n}$,

iii) $\gcd(k-1,n-j_{2})=1$,

\smallskip

\noindent then $f$ is irreducible over $\mathbb{Q}$.
\end{corollary}

The following results refer to the case when the Newton polygon of $f$ with respect to $p$ consists of an edge situated on the $x$-axis, followed by a segment  with positive slope, while the Newton polygon of $f$ with respect to $q$ consists of two segments with negative different slopes, as in Figure 1.l.

\begin{theorem}\label{teorema18} 
Let $f(X)=a_{0}+a_{1}X+\cdots +a_{n}X^{n}\in \mathbb{Z}[X]$, $a_{0}a_{n}\neq 0$. If there exist two distinct prime numbers $p,q$
and two distinct indices $j_{1},j_{2}$ such that  $j_{1}<{\rm min}(j_{2},n-j_{2})$ and
\smallskip

i) \ \thinspace $\nu _{p}(a_{i})=0$ \ for \ $i\leq j_{1}$, \  $\frac{\nu _{p}(a_{i})}{i-j_{1}}>\frac{\nu _{p}(a_{n})}{n-j_{1}}$ \ for \ $j_{1}<i<n$, \ and \ $\nu _{p}(a_{n})>0$,

ii) \ $\frac{\nu _{q}(a_{0})}{n}>\frac{\nu _{q}(a_{j_{2}})}{n-j_{2}}$, \ \ $\nu _{q}(a_{i}) > \frac{j_{2}-i}{j_{2}} \nu _{q}(a_{0})+\frac{i}{j_{2}} \nu _{q}(a_{j_{2}})$\ \ for\ \  $0<i<j_{2}$, \ \ $\nu _{q}(a_{j_{2}})>0$, 

\quad \thinspace \thinspace \thinspace $\frac{\nu _{q}(a_{i})}{n-i}  > \frac{\nu _{q}(a_{j_{2}})}{n-j_{2}} $\ \ for\ \ $j_{2}<i<n$, \ \ and \ \ $\nu _{q}(a_{n}) = 0$,

iii) $\gcd(\nu _{p}(a_{n}),n-j_{1})=\gcd(\nu _{q}(a_{0})-\nu _{q}(a_{j_{2}}),j_{2})=\gcd(\nu _{q}(a_{j_{2}}),n-j_{2})=1$,

\smallskip

\noindent then $f$ is irreducible over $\mathbb{Q}$.
\end{theorem}

In particular, one obtains the following corollary.

\begin{corollary}\label{corolarul32}
Let $f(X)=a_{0}+a_{1}X+\cdots +a_{n}X^{n}\in \mathbb{Z}[X]$, $a_{0}a_{n}\neq 0$. If there exist two distinct prime numbers $p,q$,
two distinct indices $j_{1},j_{2}$ with  $j_{1}<{\rm min}(j_{2},n-j_{2})$, and three positive integers $k_{1},k_{2},k_{3}$ such that $k_{2}>\frac{n}{n-j_{2}}k_{3}$ and   
\smallskip

i) \thinspace $p\nmid a_{i}$ \ for \ $i\leq j_{1}$, \   $p^{k_{1}}\mid a_{i}$ \ for \ $i>j_{1}$, \ $p^{k_{1}+1}\nmid a_{n}$,

ii) $q^{k_{2}}\mid a_{i}$ \ for \  $i<j_{2}$,  \ $q^{k_{2}+1}\nmid a_{0}$,
\ $q^{k_{3}}\mid a_{i}$ \ for \  $j_{2}\leq i<n$, \ $q^{k_{3}+1}\nmid a_{j_{2}}$, \ $q\nmid a_{n}$,

iii) $\gcd(k_{1},n-j_{1})=\gcd(k_{2}-k_{3},j_{2})=\gcd(k_{3},n-j_{2})=1$,

\smallskip

\noindent then $f$ is irreducible over $\mathbb{Q}$.
\end{corollary}

By letting $k_{1}=k_{3}=1$ and by writing $k$ for $k_{2}$ one obtains the following result.

\begin{corollary}\label{corolarul33}
Let $f(X)=a_{0}+a_{1}X+\cdots +a_{n}X^{n}\in \mathbb{Z}[X]$, $a_{0}a_{n}\neq 0$. If there exist two distinct prime numbers $p,q$,
two distinct indices $j_{1},j_{2}$ with  $j_{1}<{\rm min}(j_{2},n-j_{2})$, and a positive integer $k$ such that $k>\frac{n}{n-j_{2}}$ and   
\smallskip

i) \thinspace $p\nmid a_{i}$ \ for \ $i\leq j_{1}$, \ $p\mid a_{i}$ \ for \ $i>j_{1}$,  \ $p^{2}\nmid a_{n}$,

ii) $q^{k}\mid a_{i}$ \ for \  $i<j_{2}$, \ $q^{k+1}\nmid a_{0}$,
\ $q\mid a_{i}$ \ for \  $j_{2}\leq i<n$, \ $q^{2}\nmid a_{j_{2}}$,  \ $q\nmid a_{n}$,

iii) $\gcd(k-1,j_{2})=1$,

\smallskip

\noindent then $f$ is irreducible over $\mathbb{Q}$.
\end{corollary}

For the case when each one of the Newton polygons of $f$ with respect to $p$ and $q$ consists of two segments with positive 
slopes, as in Figure 1.m, we have the following result.

\begin{theorem}\label{teorema5} 
Let $f(X)=a_{0}+a_{1}X+\cdots +a_{n}X^{n}\in \mathbb{Z}[X]$, $a_{0}a_{n}\neq 0$. If there exist two distinct prime numbers $p$ and $q$
and two distinct indices $j_{1},j_{2}\in \{1,\dots ,n-1\} $ such that $j_{1}+j_{2}\neq n$ and
\smallskip

i) \ \thinspace   $\nu _{p}(a_{0}) = 0$, \ \ $\frac{\nu _{p}(a_{i})}{i} > \frac{\nu _{p}(a_{j_{1}})}{j_{2}} $ \ \ for \  $0<i<j_{1}$, \ \  $\nu _{p}(a_{j_{1}})>0$,

\quad \ $\nu _{p}(a_{i}) > \frac{n-i}{n-j_{1}} \nu _{p}(a_{j_{1}})+\frac{i-j_{1}}{n-j_{1}} \nu _{p}(a_{n})$\ \ for\ \  $j_{1}<i<n$, \ \ $\frac{\nu _{p}(a_{n})}{n}>\frac{\nu _{p}(a_{j_{1}})}{j_{1}}$,

ii) \thinspace   $\nu _{q}(a_{0}) = 0$, \ \ $\frac{\nu _{q}(a_{i})}{i} > \frac{\nu _{q}(a_{j_{2}})}{j_{2}} $ \ \ for \  $0<i<j_{2}$, \ \  $\nu _{q}(a_{j_{2}})>0$,

\quad \ $\nu _{q}(a_{i}) > \frac{n-i}{n-j_{2}} \nu _{q}(a_{j_{2}})+\frac{i-j_{2}}{n-j_{2}} \nu _{q}(a_{n})$\ \ for\ \  $j_{2}<i<n$, \ \ $\frac{\nu _{q}(a_{n})}{n}>\frac{\nu _{q}(a_{j_{2}})}{j_{2}}$, 

iii) $\gcd(\nu _{p}(a_{j_{1}}),j_{1})=\gcd(\nu _{p}(a_{n})-\nu _{p}(a_{j_{1}}),n-j_{1})=1$ and 

\quad \ \ $\gcd(\nu _{q}(a_{j_{2}}),j_{2})=\gcd(\nu _{q}(a_{n})-\nu _{q}(a_{j_{2}}),n-j_{2})=1$,

\smallskip

\noindent then $f$ is irreducible over $\mathbb{Q}$.
\end{theorem}

\begin{corollary}\label{corolarul7}
Let $f(X)=a_{0}+a_{1}X+\cdots +a_{n}X^{n}\in \mathbb{Z}[X]$, $a_{0}a_{n}\neq 0$. If there exist two distinct prime numbers $p$ and $q$,
two distinct indices $j_{1},j_{2}\in \{1,\dots ,n-1\} $ such that $j_{1}+j_{2}\neq n$, and four positive integers $k_{1},k_{2},k_{3},k_{4}$ such that $k_{1}<\frac{j_{1}}{n}k_{2}$, $k_{3}<\frac{j_{2}}{n}k_{4}$, and
\smallskip

i) \thinspace $p\nmid a_{0}$, \ \  $p^{k_{1}}\mid a_{i}$\ \ for\ \ $0<i\leq j_{1}$, \ \ $p^{k_{1}+1}\nmid a_{j_{1}}$,
\ \ $p^{k_{2}}\mid a_{i}$\ \ for\ \  $i>j_{1}$, \ \ $p^{k_{2}+1}\nmid a_{n}$,

ii) $q\nmid a_{0}$, \ \  $q^{k_{3}}\mid a_{i}$\ \ for\ \ $0<i\leq j_{2}$, \ \ $q^{k_{3}+1}\nmid a_{j_{2}}$,
\ \ $q^{k_{4}}\mid a_{i}$\ \ for\ \  $i>j_{2}$, \ \ $q^{k_{4}+1}\nmid a_{n}$,

iii) $\gcd(k_{1},j_{1})=\gcd(k_{2}-k_{1},n-j_{1})=\gcd(k_{3},j_{2})=\gcd(k_{4}-k_{3},n-j_{2})=1$,

\smallskip

\noindent then $f$ is irreducible over $\mathbb{Q}$.
\end{corollary}

By letting $k_{1}=k_{3}=1$, and writing $k_{1}$ for $k_{2}$ and $k_{2}$ for $k_{4}$ we obtain the following result.

\begin{corollary}\label{corolarul8}
Let $f(X)=a_{0}+a_{1}X+\cdots +a_{n}X^{n}\in \mathbb{Z}[X]$, $a_{0}a_{n}\neq 0$. If there exist two distinct prime numbers $p$ and $q$,
two distinct indices $j_{1},j_{2}\in \{1,\dots ,n-1\} $ such that $j_{1}+j_{2}\neq n$, and two positive integers $k_{1},k_{2}$ such that $k_{1}>\frac{n}{j_{1}}$, $k_{2}>\frac{n}{j_{2}}$, and
\smallskip

i) \thinspace \thinspace $p\nmid a_{0}$, \ \ $p\mid a_{i}$\ \ for\ \ $0<i\leq j_{1}$, \ \ $p^{2}\nmid a_{j_{1}}$,
\ \ $p^{k_{1}}\mid a_{i}$\ \ for\ \  $i>j_{1}$, \ \ $p^{k_{1}+1}\nmid a_{n}$,

ii) $q\nmid a_{0}$, \ \  $q\mid a_{i}$\ \ for\ \ $0<i\leq j_{2}$, \ \ $q^{2}\nmid a_{j_{2}}$,
\ \ $q^{k_{2}}\mid a_{i}$\ \ for\ \  $i>j_{2}$, \ \ $q^{k_{2}+1}\nmid a_{n}$,

iii) $\gcd(k_{1}-1,n-j_{1})=\gcd(k_{2}-1,n-j_{2})=1$,

\smallskip

\noindent then $f$ is irreducible over $\mathbb{Q}$.
\end{corollary}

\begin{corollary}\label{corolarul9}
Let $f(X)=a_{0}+a_{1}X+\cdots +a_{n}X^{n}\in \mathbb{Z}[X]$, $a_{0}a_{n}\neq 0$. If there exist three distinct prime numbers $p$, $q$ and $r$,
three distinct indices $j_{1},j_{2},j_{3}\in \{1,\dots ,n-1\} $ and three positive integers $k_{1},k_{2},k_{3}$ such that $k_{1}>\frac{n}{j_{1}}$, $k_{2}>\frac{n}{j_{2}}$, $k_{3}>\frac{n}{j_{3}}$, and
\smallskip

i) \thinspace \ $p\nmid a_{0}$, \ \  $p\mid a_{i}$\ \ for\ \ $0<i\leq j_{1}$, \ \ $p^{2}\nmid a_{j_{1}}$,
\ \ $p^{k_{1}}\mid a_{i}$\ \ for\ \  $i>j_{1}$, \ \ $p^{k_{1}+1}\nmid a_{n}$,

ii) \thinspace \thinspace $q\nmid a_{0}$, \ \  $q\mid a_{i}$\ \ for\ \ $0<i\leq j_{2}$, \ \ $q^{2}\nmid a_{j_{2}}$,
\ \ $q^{k_{2}}\mid a_{i}$\ \ for\ \  $i>j_{2}$, \ \ $q^{k_{2}+1}\nmid a_{n}$,

iii) $r\nmid a_{0}$, \ \  $r\mid a_{i}$\ \ for\ \ $0<i\leq j_{2}$, \ \ $r^{2}\nmid a_{j_{2}}$,
\ \ $r^{k_{3}}\mid a_{i}$\ \ for\ \  $i>j_{3}$, \ \ $r^{k_{3}+1}\nmid a_{n}$,

iv) $\gcd(k_{1}-1,n-j_{1})=\gcd(k_{2}-1,n-j_{2})=\gcd(k_{3}-1,n-j_{3})=1$,

\smallskip

\noindent then $f$ is irreducible over $\mathbb{Q}$.
\end{corollary}

The following results refer to the last case, when each one of the Newton polygons of $f$ with respect to $p$ and $q$ consists of two segments, whose intersection point lies above the $x$-axis, and one of the two Newton polygons has segments with negative different slopes, while the other has segments with positive different slopes, as in Figure 1.n. 

\begin{theorem}\label{teorema6} 
Let $f(X)=a_{0}+a_{1}X+\cdots +a_{n}X^{n}\in \mathbb{Z}[X]$, $a_{0}a_{n}\neq 0$. If there exist two distinct prime numbers $p$ and $q$
and two distinct indices $j_{1},j_{2}\in \{1,\dots ,n-1\} $ such that $j_{1}+j_{2}\neq n$ and
\smallskip

i) \ $\frac{\nu _{p}(a_{0})}{n}>\frac{\nu _{p}(a_{j_{1}})}{n-j_{1}}$, \ \ $\nu _{p}(a_{i}) > \frac{j_{1}-i}{j_{1}} \nu _{p}(a_{0})+\frac{i}{j_{1}} \nu _{p}(a_{j_{1}})$\ \ for\ \  $0<i<j_{1}$, \ \ $\nu _{p}(a_{j_{1}})>0$, 

\quad \thinspace \thinspace $\frac{\nu _{p}(a_{i})}{n-i}  > \frac{\nu _{p}(a_{j_{1}})}{n-j_{1}} $ \ \ for \ \ $j_{1}<i<n$, \ \ and \ \ $\nu _{p}(a_{n}) = 0$,

ii) \thinspace   $\nu _{q}(a_{0}) = 0$, \ \ $\frac{\nu _{q}(a_{i})}{i} > \frac{\nu _{q}(a_{j_{2}})}{j_{2}} $ \ \ for \  $0<i<j_{2}$, \ \  $\nu _{q}(a_{j_{2}})>0$,

\quad \ $\nu _{q}(a_{i}) > \frac{n-i}{n-j_{2}} \nu _{q}(a_{j_{2}})+\frac{i-j_{2}}{n-j_{2}} \nu _{q}(a_{n})$\ \ for\ \  $j_{2}<i<n$, \ \ $\frac{\nu _{q}(a_{n})}{n}>\frac{\nu _{q}(a_{j_{2}})}{j_{2}}$, 

iii) $\gcd(\nu _{p}(a_{0})-\nu _{p}(a_{j_{1}}),j_{1})=\gcd(\nu _{p}(a_{j_{1}}),n-j_{1})=1$ and 

\quad \ \ $\gcd(\nu _{q}(a_{j_{2}}),j_{2})=\gcd(\nu _{q}(a_{n})-\nu _{q}(a_{j_{2}}),n-j_{2})=1$,

\smallskip

\noindent then $f$ is irreducible over $\mathbb{Q}$.
\end{theorem}
In particular we have the following corollaries.

\begin{corollary}\label{corolarul10}
Let $f(X)=a_{0}+a_{1}X+\cdots +a_{n}X^{n}\in \mathbb{Z}[X]$, $a_{0}a_{n}\neq 0$. If there exist two distinct prime numbers $p$ and $q$,
two distinct indices $j_{1},j_{2}\in \{1,\dots ,n-1\} $ such that $j_{1}+j_{2}\neq n$, and four positive integers $k_{1},k_{2},k_{3},k_{4}$ such that $k_{2}<\frac{n-j_{1}}{n}k_{1}$, $k_{3}<\frac{j_{2}}{n}k_{4}$, and
\smallskip

i) \thinspace $p^{k_{1}}\mid a_{i}$\ \ for\ \ $i<j_{1}$, \ \ $p^{k_{1}+1}\nmid a_{0}$,
\ \ $p^{k_{2}}\mid a_{i}$\ \ for\ \  $j_{1}\leq i<n$, \ \ $p^{k_{2}+1}\nmid a_{j_{1}}$, \ \ $p\nmid a_{n}$

ii) $q\nmid a_{0}$, \ \  $q^{k_{3}}\mid a_{i}$\ \ for\ \ $0<i\leq j_{2}$, \ \ $q^{k_{3}+1}\nmid a_{j_{2}}$,
\ \ $q^{k_{4}}\mid a_{i}$\ \ for\ \  $i>j_{2}$, \ \ $q^{k_{4}+1}\nmid a_{n}$,

iii) $\gcd(k_{1}-k_{2},j_{1})=\gcd(k_{2},n-j_{1})=\gcd(k_{3},j_{2})=\gcd(k_{4}-k_{3},n-j_{2})=1$,

\smallskip

\noindent then $f$ is irreducible over $\mathbb{Q}$.
\end{corollary}

By letting $k_{2}=k_{3}=1$, and writing $k_{2}$ for $k_{4}$ we obtain the following result.

\begin{corollary}\label{corolarul11}
Let $f(X)=a_{0}+a_{1}X+\cdots +a_{n}X^{n}\in \mathbb{Z}[X]$, $a_{0}a_{n}\neq 0$. If there exist two distinct prime numbers $p$ and $q$,
two distinct indices $j_{1},j_{2}\in \{1,\dots ,n-1\} $ such that $j_{1}+j_{2}\neq n$, and two positive integers $k_{1},k_{2}$ such that $k_{1}>\frac{n}{n-j_{1}}$, $k_{2}>\frac{n}{j_{2}}$, and
\smallskip

i) \thinspace $p^{k_{1}}\mid a_{i}$\ \ for\ \ $i<j_{1}$, \ \ $p^{k_{1}+1}\nmid a_{0}$,
\ \ $p\mid a_{i}$\ \ for\ \  $j_{1}\leq i<n$, \ \ $p^{2}\nmid a_{j_{1}}$, \ \ $p\nmid a_{n}$,

ii) $q\nmid a_{0}$, \ \  $q\mid a_{i}$\ \ for\ \ $0<i\leq j_{2}$, \ \ $q^{2}\nmid a_{j_{2}}$,
\ \ $q^{k_{2}}\mid a_{i}$\ \ for\ \  $i>j_{2}$, \ \ $q^{k_{2}+1}\nmid a_{n}$,

iii) $\gcd(k_{1}-1,j_{1})=\gcd(k_{2}-1,n-j_{2})=1$,

\smallskip

\noindent then $f$ is irreducible over $\mathbb{Q}$.
\end{corollary}

Our last application of Theorem B is the following result that refers to the case when each one of the Newton polygons of $f$ with respect to $p_{1},\dots ,p_{k}$ consists of two edges, whose intersection point lies above the $x$-axis.

\begin{theorem}\label{teorema7} 
Let $f(X)=a_{0}+a_{1}X+\cdots +a_{n}X^{n}\in \mathbb{Z}[X]$, $a_{0}a_{n}\neq 0$, let $k\geq 2$ and let
$p_{1},\dots ,p_{k}$ be pair-wise distinct prime numbers. Assume that there exist indices $j_{1},\dots ,j_{k}\in \{1,\dots ,n-1\} $ such that for each $i=1,\dots ,k$ we have either
\smallskip

i) \thinspace   $\nu _{p_{i}}(a_{0}) = 0$, \ \ $\frac{\nu _{p_{i}}(a_{l})}{l} > \frac{\nu _{p_{i}}(a_{j_{i}})}{j_{i}} $ \ \ for \  $0<l<j_{i}$, \ \  $\nu _{p_{i}}(a_{j_{i}})>0$,

\quad \ $\nu _{p_{i}}(a_{l}) > \frac{n-l}{n-j_{i}} \nu _{p_{i}}(a_{j_{i}})+\frac{l-j_{i}}{n-j_{i}} \nu _{p_{i}}(a_{n})$\ \ for\ \  $j_{i}<l<n$, \ \ $\frac{\nu _{p_{i}}(a_{n})}{n}>\frac{\nu _{p_{i}}(a_{j_{i}})}{j_{i}}$,

\noindent or

ii) \thinspace $\frac{\nu _{p_{i}}(a_{0})}{n}>\frac{\nu _{p_{i}}(a_{j_{i}})}{n-j_{i}}$, \ \ $\nu _{p_{i}}(a_{l}) > \frac{j_{i}-l}{j_{i}} \nu _{p_{i}}(a_{0})+\frac{l}{j_{i}} \nu _{p_{i}}(a_{j_{i}})$ \ \ for\ \  $0<l<j_{i}$, 

\quad \   $\nu _{p_{i}}(a_{j_{i}})>0$, \ \ $\frac{\nu _{p_{i}}(a_{l})}{n-l}  > \frac{\nu _{p_{i}}(a_{j_{i}})}{n-j_{i}} $ \ \ for \ \ $j_{i}<l<n$, \ \ and \ \ $\nu _{p_{i}}(a_{n}) = 0$,

\smallskip

\noindent Let $d_{p_{i}}=\gcd \left( \frac{j_{i}}{\gcd (\nu _{p_{i}}(a_{j_{i}})-\nu _{p_{i}}(a_{0}),j_{i})}, \frac{n-j_{i}}{\gcd (\nu _{p_{i}}(a_{n})-\nu _{p_{i}}(a_{j_{i}}),n-j_{i})}   \right) $ for $i=1,\dots ,k.$

\noindent If \ $\gcd(\frac{n}{d_{p_{1}}},\dots ,\frac{n}{d_{p_{k}}})=1$, then $f$ is irreducible over $\mathbb{Q}$.
\end{theorem}
The idea to obtain irreducibility criteria by using divisibility conditions for the coefficients 
of a polynomial with respect to more than a single prime number goes back more than a century. For instance, without using Newton polygons, Perron \cite{Perron} obtained several irreducibility criteria for some classes of monic polynomials, like for instance monic polynomials whose coefficients, except for the leading one, are all divisible by the same prime numbers $p_{1},\dots ,p_{k}$. Other results of Perron refer to polynomials of the form $X^{n}+a_{n-1}X^{n-1}+\cdots +a_{1}X+a_{0}$, where all the $a_{i}$'s except for $a_{n-1}$ are divisible by $p_{1},\dots ,p_{k}$, while $a_{n-1}$ is prime to $p_{1}\cdots p_{k}$. In the language of Newton polygons, Perron's results admit simpler proofs, as shown by Mott \cite{Mott}, who proved these results in a more general setting. 

The results in this paper too may be obviously stated in a more general context. For instance, instead of polynomials with integer coefficients,
one may consider polynomials with coefficients in a unique factorization domain $R$ with quotient field $K$, and instead of rational primes, some nonassociated primes of $R$. 
Moreover, instead of polynomials with integer coefficients,
one may consider, for instance, polynomials $f(X)$ with coefficients in a field $K$ endowed with a nonarchimedean valuation $v$, where $v(K^{*})$ is a subgroup of $\mathbb{Z}$. The nonarchimedean valuation $v$ determines 
a metric topology on $K$ and, in turn, a complete field $\widehat{K}$ that will contain an isomorphic copy of $K$. Besides, $\widehat{K}$ is endowed with a nonarchimedean valuation $\hat{v}$ defined on 
$\widehat{K}$, extending the valuation $v$ on $K$, and the two valuations $v$ and $\hat{v}$ will have the same value group, and will have isomorphic residue fields (for a proof we refer the reader to 
Reiner \cite{Reiner} and Bourbaki \cite{Bourbaki}). In this case one may study the irreducibility of $f$ in ${K}[X]$ or in $\widehat{K}[X]$, and may as well consider different valuations on $K$. 
With respect to a discrete valuation $v$ on $K$, the Newton polygon of a polynomial $f(X)=a_{0}+a_{1}X+\cdots +a_{n}X^{n}\in K[X]$ is the lower convex hull of the points $(0,v(a_{0})), (1,v(a_{1})),\dots , (n,v(a_{n}))$. Thus, all the points $(i,v(a_{i}))$ lie on or above the line segments that form the edges of the Newton polygon of $f$ with respect to $v$.

For instance, in the case when the field $K$ containing the coefficients of $f$ is endowed with several discrete valuations, one may state Theorems A and B as follows.
\medskip

{\bf Theorem A".}\ {\em
Let $v_{1},\dots ,v_{k}$ be discrete valuations of a field $K$, and let $f(X)\in K[X]$ be a polynomial of degree $n$. For $i=1,\dots , k$ denote by $w_{i,1},\dots ,w_{i,n_{i}}$ the widths of all the segments of  the Newton polygon of  $f$ with respect to $v_{i}$, and by $\mathcal{S}_{v_{i}}$ the set of all the integers in the interval  $(0,\lfloor \frac{n}{2}\rfloor ]$ that may be written as a linear combination of $w_{i,1},\dots ,w_{i,n_{i}}$
with coefficients $0$ or $1$. 
If $\mathcal{S}_{v_{1}}\cap\cdots \cap\mathcal{S}_{v_{k}}=\emptyset$,
then $f$ is irreducible over $K$.
}
\medskip

{\bf Theorem B".}\ {\em Let $v_{1},\dots ,v_{k}$ be discrete valuations of a field $K$, and let $f(X)\in K[X]$ be a polynomial of degree $n$. For $i=1,\dots ,k$ let $d_{v_{i}}$ denote the greatest common divisor of the widths of the segments of  the Newton polygon of $f$ with respect to $v_{i}$. Then the degree of any factor of $f$ is divisible by $n/\gcd(\frac{n}{d_{v_{1}}},\dots ,\frac{n}{d_{v_{k}}})$. In particular, if $\gcd(\frac{n}{d_{v_{1}}},\dots ,\frac{n}{d_{v_{k}}})=1$,
then $f$ is irreducible over $K$.
}
\medskip

One may obviously state similar results for valuations of arbitrary rank over a field $K$. Such results, as well as results for multivariate polynomials that require the study of their Newton polytopes, will be provided in a forthcoming paper. For several interesting results on the irreducibility of polynomials over valued fields, of a different nature, we refer the reader to Zaharescu \cite{Zaharescu}, who extended some results of Krasner \cite{Krasner} by investigating the irreducibility of polynomials over valued fields in the presence of a secondary valuation, used to compensate the lack of completeness of the field that contains their coefficients.

The proofs of our main results are presented in Section 2 below. We will also give several examples in the last section of the paper.

\section{Proof of the main results} \label{se2} We will first recall some facts about Newton polygons (see for instance \cite{Prasolov}), that will be required in the proof of our results. 

Let $p$ be a fixed prime number, and let $f(X)=\sum _{i=0}^{n}a_{i}X^{i}$ be a polynomial with integer coefficients, $a_{0}a_{n}\neq 0$. Let us represent now the non-zero coefficients of $f$
in the form $a_{i}=\alpha _{i}p^{\beta _{i}}$ where $\alpha _{i}$ is an integer not divisible by $p$, and let us assign to each non-zero coefficient $\alpha _{i}p^{\beta _{i}}$
a point in the plane with integer coordinates $(i, \beta _{i})$. The Newton polygon of $f$ corresponding to the prime $p$ (sometimes called the Newton diagram of $f$) is constructed from these points as follows. 
Let $A_{0}=(0,\beta _{0})$ and let $A_{1}=(i_{1},\beta _{i_{1}})$, where $i_{1}$ is the largest integer for which there are no points $(i,\beta _{i})$ below the segment $A_{0}A_{1}$.
Next, let $A_{2}=(i_{2},\beta _{i_{2}})$, where $i_{2}$ is the largest integer for which there are no points $(i,\beta_{i})$ below the segment $A_{1}A_{2}$, and so on (see Figure 2).
The very last segment that we will draw will be $A_{m-1}A_{m}$, say, where $A_{m}=(n,\beta _{n})$. We note that the broken line constructed so far is the lower convex hull of the points $(i, \beta _{i}),\ i=0,\dots ,n$. Now, if some segments of the broken line $A_{0}...A_{m}$ pass through points in the plane
that have integer coordinates, then such points in the plane will be also considered  as vertices of the broken line. In this way, to the vertices $A_{0},\dots ,A_{m}$ plotted in the first phase,
we might need to add a number of $s\geq 0$ more vertices. 

\begin{center}
\hspace{1.1cm}
\setlength{\unitlength}{8mm}
\begin{picture}(12,7)
\linethickness{0.15mm}
\multiput(0,0)(1,0){12}{\line(0,1){5}}

\multiput(0,0)(0,1){6}{\line(1,0){11}}

\put(0,0){\vector(1,0){12}}
\put(0,0){\vector(0,1){6}}

\thicklines

\put(0,4){\line(1,-2){1}}   

\put(1,2){\line(2,-1){4}}

\put(5,0){\line(1,0){3}}

\put(8,0){\line(1,1){2}}

\put(10,2){\line(1,3){1}}

\put(0,4){\circle{0.08}}
\put(1,2){\circle{0.08}}
\put(3,1){\circle{0.08}}
\put(5,0){\circle{0.08}}
\put(7,0){\circle{0.08}}
\put(8,0){\circle{0.08}}
\put(10,2){\circle{0.08}}
\put(11,5){\circle{0.08}}

\put(0,4){\circle{0.12}}
\put(1,2){\circle{0.12}}
\put(5,0){\circle{0.12}}
\put(8,0){\circle{0.12}}
\put(10,2){\circle{0.12}}
\put(11,5){\circle{0.12}}

\put(2,2){\circle{0.08}}
\put(4,1){\circle{0.08}}
\put(6,1){\circle{0.08}}
\put(9,4){\circle{0.08}}

{\small 
\put(0.1,4.2){$A_{0}=B_{0}$}

\put(0.1,1.2){$A_{1}=B_{1}$}

\put(2.4,0.55){$B_{2}$}

\put(3.4,-0.45){$A_{2}=B_{3}$}

\put(5.7,-0.45){$B_{4}$}

\put(6.7,-0.45){$B_{5}$}

\put(7.6,-0.45){$A_{3}=B_{6}$}

\put(9.1,0.55){$B_{7}$}

\put(8.1,2.2){$A_{4}=B_{8}$}

\put(9.1,5.2){$A_{5}=B_{9}$}
}
\end{picture}
\end{center}
\bigskip

{\em {\small 
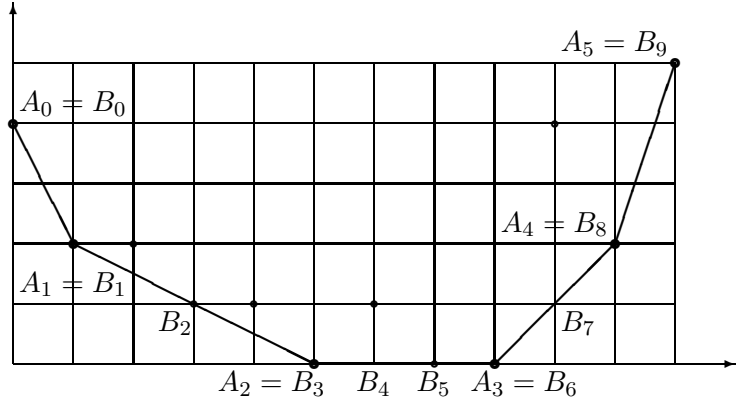
{\bf Figure 2.} The Newton polygon of $f(X)=2^{4}+2^{2}X-2^{2}X^{2}+2X^{3}-2X^{4}+X^{5}+2X^{6}-X^{7}-X^{8}+2^{4}X^{9}+2^{2}X^{10}+2^{5}X^{11}$ with respect to $p=2$.}}
\medskip

The {\it Newton polygon of} $f$ (with respect to the prime number $p$) is the resulting broken line $B_{0}...B_{m+s}$ obtained after relabelling 
all these points from left to the right, as they appear in this construction (here $B_{0}=A_{0}$ and $B_{m+s}=A_{m}$). With this notation, $A_{j}$ are called the {\it vertices} of the Newton polygon, $A_{j}A_{j+1}$ and $B_{i}B_{i+1}$ are called {\it edges} and {\it segments}
of the Newton polygon, respectively, and the vectors {$\overrightarrow{B_{i}B_{i+1}}$ are called the {\it vectors of the segments} of the Newton polygon. Therefore, a segment $B_{i}B_{i+1}$ of the Newton
polygon contains no points with integer coordinates other than its end-points $B_{i}$ and $B_{i+1}$. The collection of the vectors of the segments of the Newton polygon, taking each vector with its multiplicity,
that is as many times as it appears, is called {\it the system of vectors} for the Newton polygon. The length of the projection onto the $x$-axis of a segment (edge) of the Newton polygon is usually called the {\it width} of that segment (edge), while the projection of a segment (edge) onto the $y$-axis determines the {\it height} of that segment (edge). Usually the heights may be positive, zero, or negative, while the widths are always
taken to be positive.

An immediate consequence of the way the edges are constructed is the fact that the sequence of their slopes, when measured from left to the right, is a strictly increasing sequence, so we are not allowed to have two edges with the same slope, while the sequence of the slopes of the segments, measured also from left to the right, is an increasing, but not necessarily strictly increasing sequence.
We note here that the Newton polygon of $f$ may be alternatively defined by labelling the coefficients of $f$ in reverse order, but the useful information will be practically the same, since if $f(0)\neq 0$, then the canonical decompositions of $f$ and its reciprocal $\overline{f}=X^{\deg f}f(1/X)$ have the same type, that is the degrees of their irreducible factors along with their corresponding multiplicities coincide. 

With these definitions, we have the following celebrated result of Dumas \cite{Dumas}, which plays a central role in the theory of Newton polygons.
\medskip

{\bf Theorem (Dumas).} {\em Let $f=gh$ where $f, g$ and $h$ are non-constant polynomials with integer coefficients, and let $p$ be a prime number.
Then the system of vectors of the segments for the Newton polygon of $f$ with respect to $p$ is the union of the systems
of vectors of the segments for the Newton polygons of $g$ and $h$ with respect to $p$.
}
\medskip

In other words, Dumas' Theorem says that the edges in the Newton polygon of the product of two non-constant polynomials $g,h$ with respect to a prime number $p$ are formed by constructing a polygonal path composed by translates of all the edges that appear in the Newton polygons of $g$ and $h$ with respect to $p$, using exactly one translate for each edge, in such a way as to form a polygonal path with the slopes of the edges increasing.
\medskip

Note that the irreducibility criterion of Dumas follows immediately by Dumas' Theorem, since the Newton polygon of $f$ with respect to $p$
consists of only one segment joining the points $A_{0}=(0,0)$ and $A_{1}=(n,\nu _{p}(a_{n}))$. Indeed, according to condition i), all the points $(1,\nu _{p}(a_{1})),\dots ,(n-1,\nu _{p}(a_{n-1}))$ are situated
above $A_{0}A_{1}$, and moreover, according to condition iii), $A_{0}A_{1}$ is a segment, since it contains no other points with integer coordinates other than $A_{0}$ and $A_{1}$. 
To prove this, we see that the equation of the line passing through $A_{0}$ and $A_{1}$ is given by
\[
y=\frac{\nu _{p}(a_{n})}{n}x,
\]
so for each $i=1,\dots ,n-1$, the $y$-coordinate of a point situated on $A_{0}A_{1}$ is never a rational integer, since $n$ and $\nu _{p}(a_{n})$ are coprime.
For an alternative proof, we see that if $A_{0}A_{1}$ would contain another point with integer coordinates, say $(i,m)$ with $0<i<n$ and $0<m<\nu _{p}(a_{n})$, then by Thales' Theorem on similar triangles we would deduce that $mn=i\nu _{p}(a_{n})$, which since $gcd(\nu _{p}(a_{n}),n)=1$, would force $\nu _{p}(a_{n})$ to be a divisor of $m$, a contradiction. 

Therefore $f$ must be irreducible, for otherwise,
if $f=gh$ with $g,h\in \mathbb{Z}[X]$, $\deg g\geq 1$ and $\deg h\geq 1$, each one of the Newton polygons of its two alleged factors $g$ and $h$ with respect to the prime number $p$ would contain at least one segment, thus producing by Dumas' Theorem at least two segments in the Newton polygon of $f$.
\medskip

{\it Proof of Theorem A.} Let $f(X)$ and $p_{1},\dots ,p_{k}$ be as in the statement of our theorem, 
and let us denote by $t$ the least of the degrees of the non-constant irreducible factors of $f$.
Obviously $1\leq t\leq n$. If we assume now that $f$ is reducible over $\mathbb{Q}$, and hence by Gauss Lemma over $\mathbb{Z}$, then we must in fact have $t\leq \lfloor \frac{n}{2}\rfloor $.
Now let us look at the Newton polygon of $f$ with respect to $p_{i}$, for some fixed $i\in \{1,\dots ,k\}$. In view of Dumas' Theorem, $t$ must be the sum of some of the integers $w_{i,1},\dots , w_{i,n_{i}}$
(in such a sum, each $w_{i,l}$ enters at most once, but two such integers entering the sum, say $w_{i,l_{1}}$ and $w_{i,l_{2}}$, may be obviously equal).
Since this must happen for each $i=1,\dots ,k$, this shows that $t$ must belong to each $\mathcal{S}_{p_{i}}$, $i=1,\dots ,k$, which obviously can not hold, since according to our hypotheses, $\mathcal{S}_{p_{1}}\cap\cdots \cap\mathcal{S}_{p_{k}}=\emptyset$. Therefore the least of the degrees of the irreducible factors of $f$ must exceed $\lfloor \frac{n}{2}\rfloor $, and hence $f$ must be irreducible over $\mathbb{Q}$. This completes the proof of the theorem. $\square $

\medskip

{\it Proof of Theorem A'.}\  Let us first note that a sequence of indices $0=j_{1}<j_{2}<\cdots<j_{r}=n$ is the sequence of the abscisae of the vertices in the Newton polygon of $f$ with respect to a prime number $p$ if and only if they satisfy the following two conditions:
\smallskip

{\it i)} the sequence of slopes $u_{1},\dots ,u_{r-1}$ is strictly increasing, where
\[
u_{i}=\frac{\nu _{p}(a_{j_{i+1}})-\nu _{p}(a_{j_{i}})}{j_{i+1}-j_{i}},\qquad i=1,\dots ,r-1,
\]

{\it ii)} for each $i=1,\dots ,r-1$ and for all $k$ with $j_{i}<k<j_{i+1}$ we have
\[
\nu _{p}(a_{k})\geq \frac{j_{i+1}-k}{j_{i+1}-j_{i}}\cdot \nu _{p}(a_{j_{i}})+\frac{k-j_{i}}{j_{i+1}-j_{i}}\cdot \nu _{p}(a_{j_{i+1}}). 
\]

Indeed, condition {\it ii)} shows that there are no points $(k,\nu _{p}(a_{k}))$ with $j_{i}<k<j_{i+1}$ situated below the line segment joining the points  $P_{i}=(j_{i},\nu _{p}(a_{j_{i}}))$ and $P_{i+1}=(j_{i+1},\nu _{p}(a_{j_{i+1}}))$, while our assumption that the sequence of slopes $u_{1},\dots ,u_{r-1}$ is strictly increasing shows that for each $i=1,\dots ,r-1$, the integer $j_{i+1}$ is in fact the largest integer with the property that no point $(k,\nu _{p}(a_{k}))$ with $j_{i}<k<j_{i+1}$ is situated below the line segment $P_{i}P_{i+1}$, exactly as in the construction of the edges of a Newton polygon.
Therefore, with respect to our prime numbers $p_{1},\dots ,p_{k}$ the following conditions must be satisfied:
\smallskip

{\it i')} for each $i=1,\dots ,k$ the sequences $u_{i,1},\dots ,u_{i,r_{i}-1}$ are strictly increasing, where
\[
u_{i,l}=\frac{\nu _{p_{i}}(a_{j_{i,l+1}})-\nu _{p_{i}}(a_{j_{i,l}})}{j_{i,l+1}-j_{i,l}}, \qquad l=1,\dots ,r_{i}-1,
\]

{\it ii')} for each $i=1,\dots k$, for all $s=1,\dots ,r_{i}-1$, and for all $l$ with $j_{i,s}<l<j_{i,s+1}$ we have
\[
\nu _{p_{i}}(a_{l})\geq \frac{j_{i,s+1}-l}{j_{i,s+1}-j_{i,s}}\cdot \nu _{p_{i}}(a_{j_{i,s}})+\frac{l-j_{i,s}}{j_{i,s+1}-j_{i,s}}\cdot \nu _{p_{i}}(a_{j_{i,s+1}}).
\]
In order to prove Theorem A', all that remains to do is to notice that according to Thales' Theorem on similar triangles, the points $(x,y)$ with integer coordinates situated on 
the edge joining the vertices $P_{i,l}=(j_{i,l},\nu _{p_{i}}(a_{j_{i,l}}))$ and $P_{i,l+1}=(j_{i,l+1},\nu _{p_{i}}(a_{j_{i,l+1}}))$ in the Newton polygon of $f$ with respect to $p_{i}$ are precisely 
\[
(x_{t},y_{t})=(t\cdot x_{i,l},t\cdot y_{i,l}), \quad \makebox{with} \quad t\in\{ 0,1,\dots ,m_{i,l} \} ,
\]
where $x_{i,l}$ and $y_{i,l}$ are the coprime integers defined by the equations
\begin{eqnarray*}
j_{i,l+1}-j_{i,l} & = & m_{i,l}\cdot x_{i,l} \quad \makebox{and} \\
\nu _{p_{i}}(a_{j_{i,l+1}})-\nu _{p_{i}}(a_{j_{i,l}}) & = & m_{i,l}\cdot y_{i,l}.
\end{eqnarray*}
This shows that all the projections onto the $x$-axis of the segments situated on the edge $P_{i,l}P_{i,l+1}$ have length equal to
\[
x_{i,l}=\frac{j_{i,l+1}-j_{i,l}}{\gcd(\nu _{p_{i}}(a_{j_{i,l+1}})-\nu _{p_{i}}(a_{j_{i,l}}),j_{i,l+1}-j_{i,l}) },
\]
and moreover, the fact that the edge $P_{i,l}P_{i,l+1}$ is itself a segment of the Newton polygon, if and only if $\nu _{p_{i}}(a_{j_{i,l+1}})-\nu _{p_{i}}(a_{j_{i,l}})$ and $j_{i,l+1}-j_{i,l}$ are coprime.

According to our notation, we may now see that the Newton polygon of $f$ with respect to $p_{i}$ has precisely $m_{i,1}+\cdots +m_{i,r_{i}-1}$ segments, of which the first $m_{i,1}$ of them have width equal to
$x_{i,1}$, the following $m_{i,2}$ of them have width equal to $x_{i,2}$, and so on. Therefore, the sequence of integers $w_{i,1},\dots ,w_{i,n_{i}}$ in the statement of Theorem A will consist of precisely $m_{i,1}$ terms
equal to $x_{i,1}$, $m_{i,2}$ terms equal to $x_{i,2}$,\dots , and $m_{i,r_{i}-1}$ terms equal to $x_{i,r_{i}-1}$, and hence $n_{i}=m_{i,1}+\cdots +m_{i,r_{i}-1}$ for $i=1,\dots ,k$. In this way the sets 
$\mathcal{S}_{p_{1}},\dots ,\mathcal{S}_{p_{k}}$ may be computed explicitly by
\[
\mathcal{S}_{p_{i}} = \left\lbrace  \sum\limits _{l=1}^{r_{i}-1} n_{i,l}x_{i,l} : n_{i,l}\in \{ 0,1,\dots , m_{i,l} \} \right\rbrace \cap \left( 0,\left\lfloor \frac{n}{2}\right\rfloor \right] ,
\]
as in the statement of Theorem A'. $\square $
\medskip

{\it Proof of Theorem B.} We have to prove the fact that the degree of every non-constant factor of $f$ is divisible by $n/\gcd(\frac{n}{d_{p_{1}}},\dots ,\frac{n}{d_{p_{k}}})$.
Let us denote by $t$ the degree of a non-constant factor of $f$, and let us look at the Newton polygon of $f$ with respect to $p_{i}$, for some fixed $i\in \{1,\dots ,k\}$. In view of Dumas' Theorem,
$t$ must be the sum of some of the lengths of the segments in that Newton polygon, hence it must be a multiple of $d_{p_{i}}$. Since this must obviously hold for each $i=1,\dots, k$, we
see that $t$ must be a common multiple of $d_{p_{1}},\dots ,d_{p_{k}}$, so $t$ must be divisible by ${\rm lcm}(d_{p_{1}},\dots ,d_{p_{k}})$. The conclusion follows now by using the fact that if $d_{1},\dots ,d_{k}$ are divisors of a positive integer $n$, then
\begin{equation}\label{simplu}
{\rm lcm}(d_{1},\dots ,d_{k})=\frac{n}{\gcd(\frac{n}{d_{1}},\dots ,\frac{n}{d_{k}})}.
\end{equation}

\noindent For a proof of (\ref{simplu}), we refer the reader to \cite{Mott}, where this identity is stated as Lemma 3.11. $\square $

\medskip

{\bf Remarks.}\ i) We note that instead of using the condition $\gcd(\frac{n}{d_{p_{1}}},\dots ,\frac{n}{d_{p_{k}}})=1$, one may as well use the condition ${\rm lcm}(d_{p_{1}},\dots ,d_{p_{k}})=n$.

ii) If the Newton polygon of $f$ with respect to some $p_{i}$ consists of a single segment, then $d_{p_{i}}=n$ and hence $\gcd(\frac{n}{d_{p_{1}}},\dots ,\frac{n}{d_{p_{k}}})=1$, but in fact, in this case the irreducibility of $f$ follows directly from the Theorem of Dumas, so our result becomes useful if each one of the 
Newton polygons of $f$ with respect to $p_{1},\dots ,p_{k}$ has at least two segments, that is when $d_{p_{i}}<n$, for each $i=1,\dots ,k$.

iii) For an alternative proof of the irreducibility of $f$, that makes no use of the equality in (\ref{simplu}), let us assume to the contrary that $f$ has a non-constant
factor whose degree $t$ is strictly smaller than $n$.
Since $t$ must be a common multiple of $d_{p_{1}},\dots ,d_{p_{k}}$, we have
\[
t=\alpha _{1}\cdot d_{p_{1}}=\cdots =\alpha _{k}\cdot d_{p_{k}},
\]
for some positive integers $\alpha _{i}$, $i=1,\dots ,k$. Therefore we have
\begin{equation}\label{tab2}
\frac{n}{d_{p_{i}}}=\frac{n\alpha _{i}}{t}\qquad  \makebox{for} \ i=1,\dots ,k.
\end{equation}
If we write $\frac{n}{t}=\frac{n'}{t'}$ with $\gcd(n',t')=1$, then we must obviously have
\begin{equation}\label{t'}
t'<n'.
\end{equation}
On the other hand, in view of (\ref{tab2}) we see that
\begin{equation}\label{util}
1=\gcd\left( \frac{n'\alpha _{1}}{t'},\dots ,\frac{n'\alpha _{k}}{d_{p_{k}}}\right),
\end{equation}
and since $\frac{n'\alpha _{1}}{t'},\dots ,\frac{n'\alpha _{k}}{t'}$ must all be integers, and $t'$ and $n'$ are coprime, we deduce that $t'\mid \alpha _{i}$ for each $i=1,\dots ,k$, 
so by (\ref{util}) we must have $n'=1$. This obviously can not hold, since in view of (\ref{t'}), this would force $t'<1$, a contradiction.
Therefore we must have $t=n$.
\medskip

{\it Proof of Theorem B'.}\ In order to prove Theorem B', all that we have to do is to notice that since all the projections onto the $x$-axis of the segments situated on the edge $P_{i,l}P_{i,l+1}$ 
joining the vertices $P_{i,l}=(j_{i,l},\nu _{p_{i}}(a_{j_{i,l}}))$ and $P_{i,l+1}=(j_{i,l+1},\nu _{p_{i}}(a_{j_{i,l+1}}))$ in the Newton polygon of $f$ with respect to $p_{i}$,
have length equal to
\[
x_{i,l}=\frac{j_{i,l+1}-j_{i,l}}{\gcd(\nu _{p_{i}}(a_{j_{i,l+1}})-\nu _{p_{i}}(a_{j_{i,l}}),j_{i,l+1}-j_{i,l}) },
\]
then for a fixed, arbitrary $i\in \{ 1,\dots ,k\} $ we actually have $d_{p_{i}}=x_{i,1}=\frac{n}{\gcd(\nu _{p_{i}}(a_{n})-\nu _{p_{i}}(a_{0}),\thinspace n )}$ if $r_{i}=2$, and  $d_{p_{i}}=\gcd(x_{i,1},\dots ,x_{i,r_{i}-1})$ if $r_{i}>2$. 
$\square $
\medskip

{\it Proof of Theorem \ref{teorema1}.} Let $f(X)$ and $p_{1},\dots ,p_{k}$ be as in the statement of our theorem. First of all, we note that if at least one of the integers $\alpha _{i}$ is coprime to $n$, then 
the inequality in the statement of the theorem will in fact be a strict inequality, namely either $\nu _{p_{i}}(a_{j})>\frac{n-j}{n}\cdot \nu _{p_{i}}(a_{0})$ or 
$\nu _{p_{i}}(a_{j})>\frac{j}{n}\cdot \nu _{p_{i}}(a_{n})$ for $j=1,\dots ,n-1$, and the irreducibility
of $f$ will follow directly by the irreducibility criterion of Dumas applied to $f$ or to $\overline{f}$. Therefore we may obviously assume that none of the $\alpha _{i}$'s is coprime to $n$,
so the inequality in the statement reads either $\nu _{p_{i}}(a_{j})\geq \frac{n-j}{n}\cdot \nu _{p_{i}}(a_{0})$ or 
$\nu _{p_{i}}(a_{j})\geq \frac{j}{n}\cdot \nu _{p_{i}}(a_{n})$ for $j=1,\dots ,n-1$. The Newton polygons of $f$ with respect to $p_{1},\dots ,p_{k}$ have one of the two possible shapes, like in Figure 3, where they have been plotted together. One applies now Theorem B', with $r_{i}=2$ for each $i$. 

\begin{center}

\setlength{\unitlength}{4.65mm}
\begin{picture}(12,6)
\linethickness{0.075mm}

\put(12,0){\line(0,1){4}}   
\thicklines

\put(0,0){\line(4,1){12}}   

\put(0,2.4){\circle{0.04}}  
\put(0,2.4){\circle{0.12}}
\put(0.25,2.35){\circle{0.04}} 
\put(0.5,2.3){\circle{0.04}}   
\put(0.75,2.25){\circle{0.04}}  
\put(1,2.2){\circle{0.04}}
\put(1.25,2.15){\circle{0.04}}  
\put(1.5,2.1){\circle{0.04}}
\put(1.75,2.05){\circle{0.04}}  
\put(2,2){\circle{0.04}}
\put(2.25,1.95){\circle{0.04}}  
\put(2.5,1.9){\circle{0.04}}
\put(2.75,1.85){\circle{0.04}}  
\put(3,1.8){\circle{0.04}}
\put(3.25,1.75){\circle{0.04}}  
\put(3.5,1.7){\circle{0.04}}
\put(3.75,1.65){\circle{0.04}}  
\put(4,1.6){\circle{0.04}}
\put(4.25,1.55){\circle{0.04}}  
\put(4.5,1.5){\circle{0.04}}
\put(4.75,1.45){\circle{0.04}}  
\put(5,1.4){\circle{0.04}}
\put(5.25,1.35){\circle{0.04}}  
\put(5.5,1.3){\circle{0.04}}
\put(5.75,1.25){\circle{0.04}}  
\put(6,1.2){\circle{0.04}}
\put(6.25,1.15){\circle{0.04}}  
\put(6.5,1.1){\circle{0.04}}
\put(6.75,1.05){\circle{0.04}}  
\put(7,1){\circle{0.04}}
\put(7.25,0.95){\circle{0.04}}  
\put(7.5,0.9){\circle{0.04}}
\put(7.75,0.85){\circle{0.04}}  
\put(8,0.8){\circle{0.04}}  
\put(8.25,0.75){\circle{0.04}}
\put(8.5,0.7){\circle{0.04}}  
\put(8.75,0.65){\circle{0.04}}
\put(9,0.6){\circle{0.04}}  
\put(9.25,0.55){\circle{0.04}}
\put(9.5,0.5){\circle{0.04}}  
\put(9.75,0.45){\circle{0.04}}
\put(10,0.4){\circle{0.04}}  
\put(10.25,0.35){\circle{0.04}}
\put(10.5,0.3){\circle{0.04}}  
\put(10.75,0.25){\circle{0.04}}
\put(11,0.2){\circle{0.04}}  
\put(11.25,0.15){\circle{0.04}}
\put(11.5,0.1){\circle{0.04}}  
\put(11.75,0.05){\circle{0.04}}
\put(12,0){\circle{0.04}}

\linethickness{0.15mm}
\put(0,0){\vector(1,0){14}}
\put(0,0){\vector(0,1){5.5}}

\put(0,0){\circle{0.08}}
\put(0,0){\circle{0.12}}
\put(12,3){\circle{0.08}}
\put(12,3){\circle{0.12}}
\put(0,2.4){\circle{0.08}}
\put(0,2.4){\circle{0.12}}
\put(0,0){\circle{0.08}}
\put(0,0){\circle{0.12}}

{\small
\put(12.2,3){$A_{1}=(n,\nu _{p_{i}}(a_{n}))$}

\put(-6.4,2.4){$A_{0}'=(n,\nu _{p_{j}}(a_{0}))$}

\put(-2,-0.85){$A_{0}=(0,0)$}

\put(12,-0.85){$A_{1}'=(n,0)$}
}

\end{picture}
\end{center}
\bigskip

{\em {\small 
{\bf Figure 3.} Newton polygons consisting of a single edge.}}
\bigskip

{\it Proof of Theorem \ref{teorema8}.} \ The conditions in the statement show that the Newton polygon of $f$ with respect to $p$ consists of two segments whose slopes have different sign, while the Newton polygon of $f$ with respect $q$ consists of a single edge (composed by $\gcd(\nu _{q}(a_{n}),n)$ segments), like in Figure 4 below.

\begin{center}

\setlength{\unitlength}{4.65mm}
\begin{picture}(12,6)
\linethickness{0.075mm}

\put(12,0){\line(0,1){3.5}}   
\thicklines

\put(0,4){\line(5,-3){6.6666}}   

\put(6.6666,0){\line(3,1){5.3334}}

\put(0,0){\circle{0.12}}
\put(0,0){\circle{0.08}}
\put(0.25,0.05){\circle{0.08}}
\put(0.5,0.1){\circle{0.08}}
\put(0.75,0.15){\circle{0.08}}
\put(1,0.2){\circle{0.08}}
\put(1.25,0.25){\circle{0.08}}
\put(1.5,0.3){\circle{0.08}}
\put(1.75,0.35){\circle{0.08}}
\put(2,0.4){\circle{0.08}}
\put(2.25,0.45){\circle{0.08}}
\put(2.5,0.5){\circle{0.08}}
\put(2.75,0.55){\circle{0.08}}
\put(3,0.6){\circle{0.08}}
\put(3.25,0.65){\circle{0.08}}
\put(3.5,0.7){\circle{0.08}}
\put(3.75,0.75){\circle{0.08}}
\put(4,0.8){\circle{0.08}}
\put(4.25,0.85){\circle{0.08}}
\put(4.5,0.9){\circle{0.08}}
\put(4.75,0.95){\circle{0.08}}
\put(5,1){\circle{0.08}}
\put(5.25,1.05){\circle{0.08}}
\put(5.5,1.1){\circle{0.08}}
\put(5.75,1.15){\circle{0.08}}
\put(6,1.2){\circle{0.08}}
\put(6.25,1.25){\circle{0.08}}
\put(6.5,1.3){\circle{0.08}}
\put(6.75,1.35){\circle{0.08}}
\put(7,1.4){\circle{0.08}}
\put(7.25,1.45){\circle{0.08}}
\put(7.5,1.5){\circle{0.08}}
\put(7.75,1.55){\circle{0.08}}
\put(8,1.6){\circle{0.08}}
\put(8.25,1.65){\circle{0.08}}
\put(8.5,1.7){\circle{0.08}}
\put(8.75,1.75){\circle{0.08}}
\put(9,1.8){\circle{0.08}}
\put(9.25,1.85){\circle{0.08}}
\put(9.5,1.9){\circle{0.08}}
\put(9.75,1.95){\circle{0.08}}
\put(10,2){\circle{0.08}}
\put(10.25,2.05){\circle{0.08}}
\put(10.5,2.1){\circle{0.08}}
\put(10.75,2.15){\circle{0.08}}
\put(11,2.2){\circle{0.08}}
\put(11.25,2.25){\circle{0.08}}
\put(11.5,2.3){\circle{0.08}}
\put(11.75,2.35){\circle{0.08}}
\put(12,2.4){\circle{0.08}}
\put(12,2.4){\circle{0.12}}

\linethickness{0.15mm}

\put(0,0){\vector(1,0){14}}
\put(0,0){\vector(0,1){5.5}}

\put(6.6666,0){\circle{0.12}}

\put(0,4){\circle{0.12}}
\put(0,4){\circle{0.08}}

\put(12,1.7778){\circle{0.12}}
\put(12,1.7778){\circle{0.08}}

{\small 
\put(-6,4){$A_{0}=(0,\nu _{p}(a_{0}))$}

\put(5,-0.85){$A_{1}=(j,0)$}

\put(12.2,1.5){$A_{2}=(n,\nu _{p}(a_{n}))$}

\put(12.2,2.4){$A_{1}'=(n,\nu _{q}(a_{n}))$}

\put(-1.8,-0.85){$A_{0}'=(0,0)$}

\put(11.85,-0.85){$(n,0)$}
}

\end{picture}
\end{center}
\bigskip

{\em {\small 
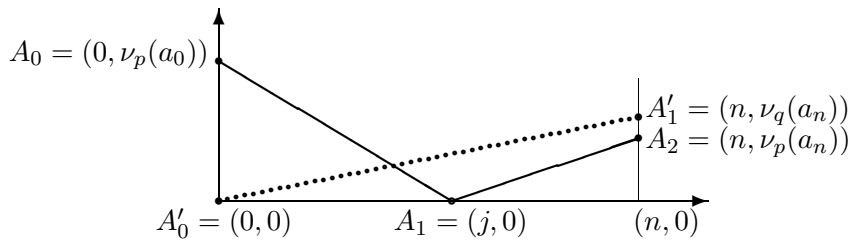
{\bf Figure 4.} The Newton polygons of $f$ with respect to $p$ and $q$, one consisting of two segments whose slopes have different sign, and the other one consisting of a single edge with positive slope.}}
\medskip

If we use now the notation from the statement of Theorem A, we see that 
\[
\mathcal{S}_{p}=\{ \min(j,n-j) \} \quad  \makebox{and} \quad \mathcal{S}_{q}=\left\lbrace  i\cdot \frac{n}{\gcd(\nu _{q}(a_{n}),n)}:i\in \mathbb{N} \right\rbrace \cap  \left( 0,\left\lfloor \frac{n}{2}\right\rfloor \right]  , 
\]
hence $\mathcal{S}_{p}\cap \mathcal{S}_{q}=\emptyset$, since $j$ is not a multiple of $\frac{n}{\gcd(\nu _{q}(a_{n}),n)}$.
\medskip

{\it Proof of Theorem \ref{teorema9}}. \ The point in this case is that the left-most edge of the Newton polygon of $f$ with respect to $p$ joins the 
points $A_{0}=(0,\nu _{p}(a_{0}))$ and $A_{1}=(1,0)$, and hence
it must be a segment in the Newton polygon, since it contains no points with integer coordinates other than $A_{0}$ and $A_{1}$, no matter what value $\nu _{p}(a_{0})$ takes. So again, the Newton polygon of 
$f$ with respect to $p$ consists of exactly two segments, while the Newton polygon of $f$ with respect to $q$ consists of one edge composed by $\gcd(\nu _{q}(a_{n}),n)$ segments. Here 
$\mathcal{S}_{p}=\{ 1 \}$ and $\mathcal{S}_{p}\cap \mathcal{S}_{q}=\emptyset$, since $n\nmid \nu _{q}(a_{n})$ and hence $\frac{n}{\gcd(\nu _{q}(a_{n}),n)}\neq 1$.
\medskip

{\it Proof of Theorem \ref{teorema10}}. \ In this case the right-most edge of the Newton polygon of $f$ with respect to $p$ joins the 
points $A_{1}=(n-1,\nu _{p}(a_{n-1}))$ and $A_{2}=(n,\nu _{p}(a_{n}))$, so
it must be a segment in the Newton polygon, since it contains no points with integer coordinates other than $A_{1}$ and $A_{2}$, no matter what value $\nu _{p}(a_{n})$ takes.

Here again, the Newton polygon of 
$f$ with respect to $p$ will consist of exactly two segments, while the Newton polygon of $f$ with respect to $q$ consists of one edge composed by $\gcd(\nu _{q}(a_{n}),n)$ segments. The rest of the proof is
similar to that of Theorem \ref{teorema9} and will be omitted.
\medskip

{\it Proof of Theorem \ref{teorema11}}. \ First of all we note that condition $j<\frac{n}{\gcd(\nu _{q}(a_{n}),n)}$ prevents $j$ to be equal to $n$, so the trivial case when $p$ divides none of the coefficients of $f$ is avoided. Now, we see that if $j=0$, then the Newton polygon of $f$ with respect to $p$ will consist of a single segment, and the irreducibility of $f$ will follow by directly using the irreducibility criterion of Dumas, so we may assume that $j\geq 1$. The conditions in the statement of the theorem show in this case that the Newton polygon of $f$ with respect to $p$ consists of two edges, the first one being composed of $j$ segments that lie on the $x$-axis, and the second one being in fact a segment 
joining the points $(j,0)$ and $(n,\nu _{p}(a_{n}))$, like in Figure 5 below.

\begin{center}
\setlength{\unitlength}{4.65mm}
\begin{picture}(12,6)
\linethickness{0.075mm}

\put(12,0){\line(0,1){3.5}}   

\thicklines

\put(0,0){\line(1,0){6.6666}}   
\put(6.6666,0){\line(3,1){5.3334}}   
\linethickness{0.15mm}

\put(0,0){\circle{0.12}}
\put(0,0){\circle{0.08}}
\put(0.25,0.05){\circle{0.08}}
\put(0.5,0.1){\circle{0.08}}
\put(0.75,0.15){\circle{0.08}}
\put(1,0.2){\circle{0.08}}
\put(1.25,0.25){\circle{0.08}}
\put(1.5,0.3){\circle{0.08}}
\put(1.75,0.35){\circle{0.08}}
\put(2,0.4){\circle{0.08}}
\put(2.25,0.45){\circle{0.08}}
\put(2.5,0.5){\circle{0.08}}
\put(2.75,0.55){\circle{0.08}}
\put(3,0.6){\circle{0.08}}
\put(3.25,0.65){\circle{0.08}}
\put(3.5,0.7){\circle{0.08}}
\put(3.75,0.75){\circle{0.08}}
\put(4,0.8){\circle{0.08}}
\put(4.25,0.85){\circle{0.08}}
\put(4.5,0.9){\circle{0.08}}
\put(4.75,0.95){\circle{0.08}}
\put(5,1){\circle{0.08}}
\put(5.25,1.05){\circle{0.08}}
\put(5.5,1.1){\circle{0.08}}
\put(5.75,1.15){\circle{0.08}}
\put(6,1.2){\circle{0.08}}
\put(6.25,1.25){\circle{0.08}}
\put(6.5,1.3){\circle{0.08}}
\put(6.75,1.35){\circle{0.08}}
\put(7,1.4){\circle{0.08}}
\put(7.25,1.45){\circle{0.08}}
\put(7.5,1.5){\circle{0.08}}
\put(7.75,1.55){\circle{0.08}}
\put(8,1.6){\circle{0.08}}
\put(8.25,1.65){\circle{0.08}}
\put(8.5,1.7){\circle{0.08}}
\put(8.75,1.75){\circle{0.08}}
\put(9,1.8){\circle{0.08}}
\put(9.25,1.85){\circle{0.08}}
\put(9.5,1.9){\circle{0.08}}
\put(9.75,1.95){\circle{0.08}}
\put(10,2){\circle{0.08}}
\put(10.25,2.05){\circle{0.08}}
\put(10.5,2.1){\circle{0.08}}
\put(10.75,2.15){\circle{0.08}}
\put(11,2.2){\circle{0.08}}
\put(11.25,2.25){\circle{0.08}}
\put(11.5,2.3){\circle{0.08}}
\put(11.75,2.35){\circle{0.08}}
\put(12,2.4){\circle{0.08}}
\put(12,2.4){\circle{0.12}}

\put(0,0){\vector(1,0){14}}
\put(0,0){\vector(0,1){5.5}}

\put(6.6666,0){\circle{0.12}}

\put(12,1.7778){\circle{0.12}}
\put(12,1.7778){\circle{0.08}}

{\small 

\put(5,-0.85){$A_{1}=(j,0)$}

\put(12.3,1.5){$A_{2}=(n,\nu _{p}(a_{n}))$}

\put(12.3,2.4){$A_{1}'=(n,\nu _{q}(a_{n}))$}

\put(-2.8,-0.85){$A_{0}=A_{0}'=(0,0)$}

\put(11.85,-0.85){$(n,0)$}
}
\end{picture}
\end{center}
\bigskip

{\em {\small 
{\bf Figure 5.} The Newton polygons of $f$ with respect to $p$ and $q$, one consisting of an edge situated on the $x$-axis, followed by a segment with positive slope, and the other one consisting of a single edge
with positive slope.}}
\medskip

Using the notations in Theorem A, we deduce in this case that
\[
\mathcal{S}_{p}=\{ 0,1,\dots ,j \}\cap  \left( 0,\left\lfloor \frac{n}{2}\right\rfloor \right]  \quad  \makebox{and} \quad \mathcal{S}_{q}=\left\lbrace  i\cdot \frac{n}{\gcd(\nu _{q}(a_{n}),n)}:i\in \mathbb{N} \right\rbrace \cap  \left( 0,\left\lfloor \frac{n}{2}\right\rfloor \right]  , 
\]
so $\mathcal{S}_{p}\cap \mathcal{S}_{q}=\emptyset$, since $j<\frac{n}{\gcd(\nu _{q}(a_{n}),n)}$.  $\square $
\medskip

{\it Proof of Theorem \ref{teorema12}}. \ Here the condition $j>n-\frac{n}{ \gcd(\nu _{q}(a_{n}),n) }$ prevents $j$ to be equal to $0$, thus avoiding again the trivial case when $p$ divides none of the coefficients of $f$. 
Now we see that if $j=n$, then the Newton polygon of $f$ with respect to $p$ will consist of a single segment, and the irreducibility of $f$ follows by directly using the irreducibility criterion of Dumas for the reciprocal of $f$, so here we may assume that $j\leq n-1$. The conditions in the statement of the theorem show in this case that the Newton polygon of $f$ with respect to $p$ consists of two edges, the first one being
in fact a segment joining the points $(0,\nu _{p}(a_{0}))$ and $(j,0)$, and the second one being composed of $n-j$ segments that lie on the $x$-axis, like in Figure 6 below. 

\begin{center}
\setlength{\unitlength}{4.65mm}
\begin{picture}(12,6)
\linethickness{0.075mm}

\put(12,0){\line(0,1){3.5}}   
\thicklines

\put(0,4){\line(5,-3){6.6666}}   

\put(6.6666,0){\line(1,0){5.3334}}   
\linethickness{0.15mm}

\put(0,0){\circle{0.12}}
\put(0,0){\circle{0.08}}
\put(0.25,0.05){\circle{0.08}}
\put(0.5,0.1){\circle{0.08}}
\put(0.75,0.15){\circle{0.08}}
\put(1,0.2){\circle{0.08}}
\put(1.25,0.25){\circle{0.08}}
\put(1.5,0.3){\circle{0.08}}
\put(1.75,0.35){\circle{0.08}}
\put(2,0.4){\circle{0.08}}
\put(2.25,0.45){\circle{0.08}}
\put(2.5,0.5){\circle{0.08}}
\put(2.75,0.55){\circle{0.08}}
\put(3,0.6){\circle{0.08}}
\put(3.25,0.65){\circle{0.08}}
\put(3.5,0.7){\circle{0.08}}
\put(3.75,0.75){\circle{0.08}}
\put(4,0.8){\circle{0.08}}
\put(4.25,0.85){\circle{0.08}}
\put(4.5,0.9){\circle{0.08}}
\put(4.75,0.95){\circle{0.08}}
\put(5,1){\circle{0.08}}
\put(5.25,1.05){\circle{0.08}}
\put(5.5,1.1){\circle{0.08}}
\put(5.75,1.15){\circle{0.08}}
\put(6,1.2){\circle{0.08}}
\put(6.25,1.25){\circle{0.08}}
\put(6.5,1.3){\circle{0.08}}
\put(6.75,1.35){\circle{0.08}}
\put(7,1.4){\circle{0.08}}
\put(7.25,1.45){\circle{0.08}}
\put(7.5,1.5){\circle{0.08}}
\put(7.75,1.55){\circle{0.08}}
\put(8,1.6){\circle{0.08}}
\put(8.25,1.65){\circle{0.08}}
\put(8.5,1.7){\circle{0.08}}
\put(8.75,1.75){\circle{0.08}}
\put(9,1.8){\circle{0.08}}
\put(9.25,1.85){\circle{0.08}}
\put(9.5,1.9){\circle{0.08}}
\put(9.75,1.95){\circle{0.08}}
\put(10,2){\circle{0.08}}
\put(10.25,2.05){\circle{0.08}}
\put(10.5,2.1){\circle{0.08}}
\put(10.75,2.15){\circle{0.08}}
\put(11,2.2){\circle{0.08}}
\put(11.25,2.25){\circle{0.08}}
\put(11.5,2.3){\circle{0.08}}
\put(11.75,2.35){\circle{0.08}}
\put(12,2.4){\circle{0.08}}
\put(12,2.4){\circle{0.12}}

\put(0,0){\vector(1,0){14}}
\put(0,0){\vector(0,1){5.5}}

\put(6.6666,0){\circle{0.12}}

\put(0,4){\circle{0.12}}
\put(0,4){\circle{0.08}}

\put(12,0){\circle{0.12}}
\put(12,0){\circle{0.08}}

{\small 
\put(-6,4){$A_{0}=(0,\nu _{p}(a_{0}))$}

\put(5,-0.85){$A_{1}=(j,0)$}

\put(12.3,2.4){$A_{1}'=(n,\nu _{q}(a_{n}))$}

\put(-1.8,-0.85){$A_{0}'=(0,0)$}

\put(10.85,-0.85){$A_{2}=(n,0)$}
}

\end{picture}
\end{center}
\bigskip

{\em {\small 
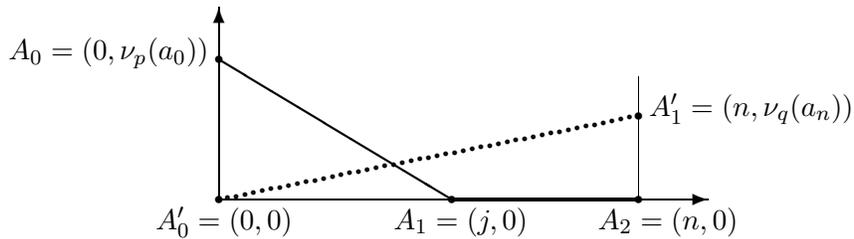
{\bf Figure 6.} The Newton polygons of $f$ with respect to $p$ and $q$, one consisting of a segment with negative slope, followed by an edge situated on the $x$-axis, and the other one consisting of a single edge
with positive slope.}}
\medskip

Using again the notations in Theorem A, we deduce in this case that
\[
\mathcal{S}_{p}=\{ 0,1,\dots ,n-j \}\cap  \left( 0,\left\lfloor \frac{n}{2}\right\rfloor \right]  \quad  \makebox{and} \quad \mathcal{S}_{q}=\left\lbrace  i\cdot \frac{n}{\gcd(\nu _{q}(a_{n}),n)}:i\in \mathbb{N} \right\rbrace \cap  \left( 0,\left\lfloor \frac{n}{2}\right\rfloor \right]  , 
\]
hence $\mathcal{S}_{p}\cap \mathcal{S}_{q}=\emptyset$, since $n-j<\frac{n}{\gcd(\nu _{q}(a_{n}),n)}$.  $\square $
\medskip

{\it Proof of Theorem \ref{teorema13}}. \  If $\gcd(\nu _{q}(a_{n}),n)=1$, the irreducibility of $f$ is a direct consequence of the irreducibility criterion of Dumas, so we may obviously assume that $\gcd(\nu _{q}(a_{n}),n)>1$,
which in view of our condition that $j<\frac{n}{\gcd(\nu _{q}(a_{n}),n)}$ yields $j\leq \lfloor \frac{n}{2}\rfloor $.
The conditions in the statement of the theorem show in this case that
the Newton polygon of $f$ with respect to $p$ consists of two segments with positive different slopes, while the Newton polygon of $f$ with respect to $q$ consists of a single edge (composed of $\gcd(\nu _{q}(a_{n}),n)$ 
segments) with positive slope too, like in Figure 7 below.

\begin{center}
\setlength{\unitlength}{4.65mm}
\begin{picture}(12,6)
\linethickness{0.075mm}

\put(12,0){\line(0,1){5.4}}   
\thicklines

\put(0,0){\line(6,1){5}}   
\put(5,0.83333){\line(2,1){7}}   
\linethickness{0.15mm}

\put(0,0){\circle{0.12}}
\put(0,0){\circle{0.08}}
\put(0.25,0.1){\circle{0.08}}
\put(0.5,0.2){\circle{0.08}}
\put(0.75,0.3){\circle{0.08}}
\put(1,0.4){\circle{0.08}}
\put(1.25,0.5){\circle{0.08}}
\put(1.5,0.6){\circle{0.08}}
\put(1.75,0.7){\circle{0.08}}
\put(2,0.8){\circle{0.08}}
\put(2.25,0.9){\circle{0.08}}
\put(2.5,1){\circle{0.08}}
\put(2.75,1.1){\circle{0.08}}
\put(3,1.2){\circle{0.08}}
\put(3.25,1.3){\circle{0.08}}
\put(3.5,1.4){\circle{0.08}}
\put(3.75,1.5){\circle{0.08}}
\put(4,1.6){\circle{0.08}}
\put(4.25,1.7){\circle{0.08}}
\put(4.5,1.8){\circle{0.08}}
\put(4.75,1.9){\circle{0.08}}
\put(5,2){\circle{0.08}}
\put(5.25,2.1){\circle{0.08}}
\put(5.5,2.2){\circle{0.08}}
\put(5.75,2.3){\circle{0.08}}
\put(6,2.4){\circle{0.08}}
\put(6.25,2.5){\circle{0.08}}
\put(6.5,2.6){\circle{0.08}}
\put(6.75,2.7){\circle{0.08}}
\put(7,2.8){\circle{0.08}}
\put(7.25,2.9){\circle{0.08}}
\put(7.5,3){\circle{0.08}}
\put(7.75,3.1){\circle{0.08}}
\put(8,3.2){\circle{0.08}}
\put(8.25,3.3){\circle{0.08}}
\put(8.5,3.4){\circle{0.08}}
\put(8.75,3.5){\circle{0.08}}
\put(9,3.6){\circle{0.08}}
\put(9.25,3.7){\circle{0.08}}
\put(9.5,3.8){\circle{0.08}}
\put(9.75,3.9){\circle{0.08}}
\put(10,4){\circle{0.08}}
\put(10.25,4.1){\circle{0.08}}
\put(10.5,4.2){\circle{0.08}}
\put(10.75,4.3){\circle{0.08}}
\put(11,4.4){\circle{0.08}}
\put(11.25,4.5){\circle{0.08}}
\put(11.5,4.6){\circle{0.08}}
\put(11.75,4.7){\circle{0.08}}
\put(12,4.8){\circle{0.08}}
\put(12,4.8){\circle{0.12}}

\put(0,0){\vector(1,0){14}}
\put(0,0){\vector(0,1){5.5}}

\put(5,0.8333){\circle{0.12}}
\put(5,0.8333){\circle{0.08}}

\put(12,4.333){\circle{0.12}}
\put(12,4.333){\circle{0.08}}

{\small 

\put(5.8,0.55){$A_{1}=(j,\nu _{p}(a_{j}))$}

\put(12.3,3.9){$A_{2}=(n,\nu _{p}(a_{n}))$}

\put(12.3,4.9){$A_{1}'=(n,\nu _{q}(a_{n}))$}

\put(-2.8,-0.85){$A_{0}=A_{0}'=(0,0)$}

\put(11.85,-0.85){$(n,0)$}
}
\end{picture}
\end{center}
\bigskip

{\em {\small 
{\bf Figure 7.} The Newton polygons of $f$ with respect to $p$ and $q$, one consisting of two segments with positive different slopes, and the other one consisting of a single edge with positive slope.}}
\medskip

The condition $\frac{\nu _{p}(a_{n})}{n}>\frac{\nu _{p}(a_{j})}{j}$ shows that the slope of the segment $A_{1}A_{2}$ exceeds the slope of the segment $A_{0}A_{1}$.
Using again the notations in Theorem A, we deduce that
\[
\mathcal{S}_{p}=\{ j \}  \quad  \makebox{and} \quad \mathcal{S}_{q}=\left\lbrace  i\cdot \frac{n}{\gcd(\nu _{q}(a_{n}),n)}:i\in \mathbb{N} \right\rbrace \cap  \left( 0,\left\lfloor \frac{n}{2}\right\rfloor \right]  , 
\]
so $\mathcal{S}_{p}\cap \mathcal{S}_{q}=\emptyset$, since $j<\frac{n}{\gcd(\nu _{q}(a_{n}),n)}$.  $\square $
\medskip

{\it Proof of Theorem \ref{teorema14}}. \  Again, if $\gcd(\nu _{q}(a_{n}),n)=1$, then the irreducibility of $f$ is an immediate consequence of the irreducibility criterion of Dumas, so we may obviously assume that $\gcd(\nu _{q}(a_{n}),n)>1$,
which in view of our condition that $j<\frac{n}{\gcd(\nu _{q}(a_{n}),n)}$ yields $j\leq \lfloor \frac{n}{2}\rfloor $.
The conditions in the statement of the theorem show in this case that
the Newton polygon of $f$ with respect to $p$ consists of two segments with negative different slopes, while the Newton polygon of $f$ with respect to $q$ consists of a single edge (composed of $\gcd(\nu _{q}(a_{n}),n)$ 
segments) with positive slope, like in Figure 8 below.

\begin{center}

\setlength{\unitlength}{4.65mm}
\begin{picture}(12,6)
\linethickness{0.075mm}

\put(12,0){\line(0,1){3.5}}  
\thicklines

\put(0,4){\line(1,-1){2}}   

\put(2,2){\line(5,-1){10}}   
\linethickness{0.15mm}

\put(0,0){\circle{0.12}}
\put(0,0){\circle{0.08}}
\put(0.25,0.05){\circle{0.08}}
\put(0.5,0.1){\circle{0.08}}
\put(0.75,0.15){\circle{0.08}}
\put(1,0.2){\circle{0.08}}
\put(1.25,0.25){\circle{0.08}}
\put(1.5,0.3){\circle{0.08}}
\put(1.75,0.35){\circle{0.08}}
\put(2,0.4){\circle{0.08}}
\put(2.25,0.45){\circle{0.08}}
\put(2.5,0.5){\circle{0.08}}
\put(2.75,0.55){\circle{0.08}}
\put(3,0.6){\circle{0.08}}
\put(3.25,0.65){\circle{0.08}}
\put(3.5,0.7){\circle{0.08}}
\put(3.75,0.75){\circle{0.08}}
\put(4,0.8){\circle{0.08}}
\put(4.25,0.85){\circle{0.08}}
\put(4.5,0.9){\circle{0.08}}
\put(4.75,0.95){\circle{0.08}}
\put(5,1){\circle{0.08}}
\put(5.25,1.05){\circle{0.08}}
\put(5.5,1.1){\circle{0.08}}
\put(5.75,1.15){\circle{0.08}}
\put(6,1.2){\circle{0.08}}
\put(6.25,1.25){\circle{0.08}}
\put(6.5,1.3){\circle{0.08}}
\put(6.75,1.35){\circle{0.08}}
\put(7,1.4){\circle{0.08}}
\put(7.25,1.45){\circle{0.08}}
\put(7.5,1.5){\circle{0.08}}
\put(7.75,1.55){\circle{0.08}}
\put(8,1.6){\circle{0.08}}
\put(8.25,1.65){\circle{0.08}}
\put(8.5,1.7){\circle{0.08}}
\put(8.75,1.75){\circle{0.08}}
\put(9,1.8){\circle{0.08}}
\put(9.25,1.85){\circle{0.08}}
\put(9.5,1.9){\circle{0.08}}
\put(9.75,1.95){\circle{0.08}}
\put(10,2){\circle{0.08}}
\put(10.25,2.05){\circle{0.08}}
\put(10.5,2.1){\circle{0.08}}
\put(10.75,2.15){\circle{0.08}}
\put(11,2.2){\circle{0.08}}
\put(11.25,2.25){\circle{0.08}}
\put(11.5,2.3){\circle{0.08}}
\put(11.75,2.35){\circle{0.08}}
\put(12,2.4){\circle{0.08}}
\put(12,2.4){\circle{0.12}}

\put(0,0){\vector(1,0){14}}
\put(0,0){\vector(0,1){5.5}}

\put(2,2){\circle{0.12}}
\put(2,2){\circle{0.08}}

\put(0,4){\circle{0.12}}
\put(0,4){\circle{0.08}}

\put(12,0){\circle{0.12}}
\put(12,0){\circle{0.08}}

{\small 
\put(-6,4){$A_{0}=(0,\nu _{p}(a_{0}))$}

\put(2,2.5){$A_{1}=(j,\nu _{p}(a_{j}))$}

\put(12.3,2.4){$A_{1}'=(n,\nu _{q}(a_{n}))$}

\put(-1.8,-0.85){$A_{0}'=(0,0)$}

\put(10.85,-0.85){$A_{2}=(n,0)$}
}

\end{picture}
\end{center}
\bigskip

{\em {\small 
{\bf Figure 8.} The Newton polygons of $f$ with respect to $p$ and $q$, one consisting of two segments with negative different slopes, and the other one consisting of a single edge with positive slope.}}
\medskip

Here the condition $\frac{\nu _{p}(a_{0})}{n}>\frac{\nu _{p}(a_{j})}{n-j}$ shows that the slope of the segment $A_{1}A_{2}$ exceeds the slope of the segment $A_{0}A_{1}$.
In this case we have
\[
\mathcal{S}_{p}=\{ j \}  \quad  \makebox{and} \quad \mathcal{S}_{q}=\left\lbrace  i\cdot \frac{n}{\gcd(\nu _{q}(a_{n}),n)}:i\in \mathbb{N} \right\rbrace \cap  \left( 0,\left\lfloor \frac{n}{2}\right\rfloor \right]  , 
\]
hence $\mathcal{S}_{p}\cap \mathcal{S}_{q}=\emptyset$, since $j<\frac{n}{\gcd(\nu _{q}(a_{n}),n)}$.  $\square $
\medskip

{\it Proof of Theorem \ref{teorema2}.}\ The conditions in the statement of the theorem show that each one of the Newton polygons of $f$ with respect to $p$ and $q$
consists of exactly two segments who intersect on the $x$-axis, and whose slopes have different sign, as in Figure 9 below. 

\begin{center}
\setlength{\unitlength}{4.65mm}
\begin{picture}(12,6)
\linethickness{0.075mm}

\put(12,0){\line(0,1){3.5}}   
\thicklines

\put(0,4){\line(5,-3){6.6666}}   

\put(6.6666,0){\line(3,1){5.3334}}

\put(0,3){\circle{0.12}}
\put(0,3){\circle{0.08}}
\put(0.25,2.75){\circle{0.08}}
\put(0.5,2.5){\circle{0.08}}
\put(0.75,2.25){\circle{0.08}}
\put(1,2){\circle{0.08}}
\put(1.25,1.75){\circle{0.08}}
\put(1.5,1.5){\circle{0.08}}
\put(1.75,1.25){\circle{0.08}}
\put(2,1){\circle{0.08}}
\put(2.25,0.75){\circle{0.08}}
\put(2.5,0.5){\circle{0.08}}
\put(2.75,0.25){\circle{0.08}}
\put(3,0){\circle{0.12}}
\put(3,0){\circle{0.08}}

\put(3,0){\circle{0.12}}
\put(3,0){\circle{0.08}}
\put(3.35,0.1){\circle{0.08}}
\put(3.7,0.2){\circle{0.08}}
\put(4.05,0.3){\circle{0.08}}
\put(4.4,0.4){\circle{0.08}}
\put(4.75,0.5){\circle{0.08}}
\put(5.1,0.6){\circle{0.08}}
\put(5.45,0.7){\circle{0.08}}
\put(5.8,0.8){\circle{0.08}}
\put(6.15,0.9){\circle{0.08}}
\put(6.5,1){\circle{0.08}}
\put(6.85,1.1){\circle{0.08}}
\put(7.2,1.2){\circle{0.08}}
\put(7.55,1.3){\circle{0.08}}
\put(7.9,1.4){\circle{0.08}}
\put(8.25,1.5){\circle{0.08}}
\put(8.6,1.6){\circle{0.08}}
\put(8.95,1.7){\circle{0.08}}
\put(9.3,1.8){\circle{0.08}}
\put(9.65,1.9){\circle{0.08}}
\put(10,2){\circle{0.08}}
\put(10.35,2.1){\circle{0.08}}
\put(10.7,2.2){\circle{0.08}}
\put(11.05,2.3){\circle{0.08}}
\put(11.4,2.4){\circle{0.08}}
\put(11.7,2.48){\circle{0.08}}
\put(12,2.56){\circle{0.08}}
\put(12,2.56){\circle{0.12}}

\linethickness{0.15mm}

\put(0,0){\vector(1,0){14}}
\put(0,0){\vector(0,1){5.5}}

\put(6.6666,0){\circle{0.12}}

\put(0,4){\circle{0.12}}
\put(0,4){\circle{0.08}}

\put(12,1.7778){\circle{0.12}}
\put(12,1.7778){\circle{0.08}}

{\small  
\put(-6,4){$A_{0}=(0,\nu _{p}(a_{0}))$}

\put(-6,3){$A_{0}'=(0,\nu _{q}(a_{0}))$}

\put(6.3,-0.85){$A_{1}=(j_{1},0)$}

\put(1.5,-0.85){$A_{1}'=(j_{2},0)$}

\put(12.3,1.7778){$A_{2}=(n,\nu _{p}(a_{n}))$}

\put(12.3,2.7778){$A_{2}'=(n,\nu _{q}(a_{n}))$}

\put(-1.5,-0.85){$(0,0)$}

\put(11.85,-0.85){$(n,0)$}
}

\end{picture}
\end{center}
\bigskip

{\em {\small 
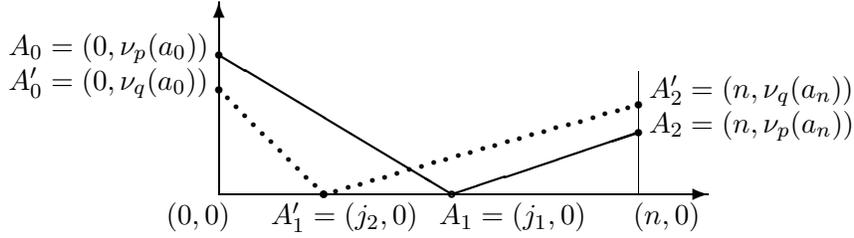
{\bf Figure 9.} The Newton polygons of $f$ with respect to $p$ and $q$, each one consisting of two edges whose slopes have different sign.}}
\medskip

Therefore, we see that 
\[
\mathcal{S}_{p}=\{ \min(j_{1},n-j_{1}) \} \quad  \makebox{and} \quad \mathcal{S}_{q}=\{ \min(j_{2},n-j_{2}) \} , 
\]
so $\mathcal{S}_{p}\cap \mathcal{S}_{q}=\emptyset$, since $j_{1}\neq j_{2}$ and $j_{1}+j_{2}\neq n$.  $\square $
\medskip

{\it Proof of Corollary \ref{corolarul4}.} Here all we have to do is to observe that among any three different indices in the set $\{ 1,\dots ,n-1\} $, there exist two of them whose sum is not equal to $n$.  $\square $

\medskip

{\it Proof of Theorem \ref{teorema3}.} The proof of Theorem \ref{teorema2} applies here as well. The point in this case is that the left-most edge of the Newton polygon of $f$ with respect to $p$ joins the 
points $A_{0}=(0,\nu _{p}(a_{0}))$ and $A_{1}=(1,0)$, and hence
$A_{0}A_{1}$ must be a segment in the Newton polygon, since it contains no points with integer coordinates other than $A_{0}$ and $A_{1}$, no matter what value $\nu _{p}(a_{0})$ takes.
$\square $

\medskip

{\it Proof of Theorem \ref{teorema15}.} \  If $j_{1}=0$, the irreducibility of $f$ follows from the irreducibility criterion of Dumas, so we may assume that $j_{1}>0$. In this case the Newton polygon of $f$ with respect to $p$ consists of two edges, the first one being composed of $j_{1}$ segments situated on the $x$-axis, and the second one being in fact a segment with positive slope,
while the Newton polygon of $f$ with respect to $q$ consists of two segments whose slopes have different sign, as in Figure 10 below.

\begin{center}

\setlength{\unitlength}{4.65mm}
\begin{picture}(12,6)
\linethickness{0.075mm}

\put(12,0){\line(0,1){4.5}}   
\thicklines

\put(0,0){\line(1,0){4}}   

\put(4,0){\line(2,1){8}}   

\linethickness{0.15mm}

\put(0,4){\circle{0.12}}
\put(0,4){\circle{0.08}}
\put(0.25,3.875){\circle{0.08}}
\put(0.5,3.75){\circle{0.08}}
\put(0.75,3.625){\circle{0.08}}
\put(1,3.5){\circle{0.08}}
\put(1.25,3.375){\circle{0.08}}
\put(1.5,3.25){\circle{0.08}}
\put(1.75,3.125){\circle{0.08}}
\put(2,3){\circle{0.08}}
\put(2.25,2.875){\circle{0.08}}
\put(2.5,2.75){\circle{0.08}}
\put(2.75,2.625){\circle{0.08}}
\put(3,2.5){\circle{0.08}}
\put(3.25,2.375){\circle{0.08}}

\put(3.5,2.25){\circle{0.08}}
\put(3.75,2.125){\circle{0.08}}
\put(4,2){\circle{0.08}}
\put(4.25,1.875){\circle{0.08}}
\put(4.5,1.75){\circle{0.08}}
\put(4.75,1.625){\circle{0.08}}
\put(5,1.5){\circle{0.08}}
\put(5.25,1.375){\circle{0.08}}
\put(5.5,1.25){\circle{0.08}}
\put(5.75,1.125){\circle{0.08}}
\put(6,1){\circle{0.08}}
\put(6.25,0.875){\circle{0.08}}
\put(6.5,0.75){\circle{0.08}}
\put(6.75,0.625){\circle{0.08}}
\put(7,0.5){\circle{0.08}}
\put(7.25,0.375){\circle{0.08}}
\put(7.5,0.25){\circle{0.08}}
\put(7.75,0.125){\circle{0.08}}
\put(8,0){\circle{0.08}}
\put(8,0){\circle{0.12}}

\put(8.25,0.125){\circle{0.08}}
\put(8.5,0.25){\circle{0.08}}
\put(8.75,0.375){\circle{0.08}}
\put(9,0.5){\circle{0.08}}
\put(9.25,0.625){\circle{0.08}}
\put(9.5,0.75){\circle{0.08}}
\put(9.75,0.875){\circle{0.08}}
\put(10,1){\circle{0.08}}
\put(10.25,1.125){\circle{0.08}}
\put(10.5,1.25){\circle{0.08}}
\put(10.75,1.375){\circle{0.08}}
\put(11,1.5){\circle{0.08}}
\put(11.25,1.625){\circle{0.08}}
\put(11.5,1.75){\circle{0.08}}
\put(11.75,1.875){\circle{0.08}}
\put(12,2){\circle{0.08}}
\put(12,2){\circle{0.12}}

\put(0,0){\vector(1,0){14}}
\put(0,0){\vector(0,1){5.5}}

\put(0,0){\circle{0.12}}
\put(0,0){\circle{0.08}}

\put(4,0){\circle{0.12}}
\put(4,0){\circle{0.08}}

\put(0,4){\circle{0.12}}
\put(0,4){\circle{0.08}}

\put(12,4){\circle{0.12}}
\put(12,4){\circle{0.08}}

{\small  

\put(-6,4){$A_{0}'=(0,\nu _{q}(a_{0}))$}

\put(6.3,-0.85){$A_{1}'=(j_{2},0)$}

\put(1.5,-0.85){$A_{1}=(j_{1},0)$}

\put(12.3,4){$A_{2}=(n,\nu _{p}(a_{n}))$}

\put(12.3,2){$A_{2}'=(n,\nu _{q}(a_{n}))$}

\put(-3,-0.85){$A_{0}=(0,0)$}

\put(11.85,-0.85){$(n,0)$}
}

\end{picture}
\end{center}
\bigskip

{\em {\small 
{\bf Figure 10.} The Newton polygons of $f$ with respect to $p$ and $q$, one consisting of an edge situated on the $x$-axis, followed by a segment with positive slope, and the other one consisting of two segments whose slopes have different sign.}}
\medskip

Therefore, we see in this case that 
\[
\mathcal{S}_{p}=\{ 1,\dots ,j_{1} \} \quad  \makebox{and} \quad \mathcal{S}_{q}=\{ \min(j_{2},n-j_{2}) \} , 
\]
so $\mathcal{S}_{p}\cap \mathcal{S}_{q}=\emptyset$, since $j_{1}<\min(j_{2},n-j_{2})$.  $\square $
\medskip

{\it Proof of Theorem \ref{teorema16}.} \ In this case the Newton polygon of $f$ with respect to $p$ consists of two segments with positive different slopes, while the Newton polygon of $f$ with respect to $q$ consists
of two segments whose slopes have different sign, as in Figure 11 below.

\begin{center}
\setlength{\unitlength}{4.65mm}
\begin{picture}(12,6)
\linethickness{0.075mm}

\put(12,0){\line(0,1){4.8}}   
\thicklines

\put(0,0){\line(5,1){5}}   

\put(5,1){\line(2,1){7}}   

\linethickness{0.15mm}

\put(0,4){\circle{0.12}}
\put(0,4){\circle{0.08}}
\put(0.25,3.875){\circle{0.08}}
\put(0.5,3.75){\circle{0.08}}
\put(0.75,3.625){\circle{0.08}}
\put(1,3.5){\circle{0.08}}
\put(1.25,3.375){\circle{0.08}}
\put(1.5,3.25){\circle{0.08}}
\put(1.75,3.125){\circle{0.08}}
\put(2,3){\circle{0.08}}
\put(2.25,2.875){\circle{0.08}}
\put(2.5,2.75){\circle{0.08}}
\put(2.75,2.625){\circle{0.08}}
\put(3,2.5){\circle{0.08}}
\put(3.25,2.375){\circle{0.08}}

\put(3.5,2.25){\circle{0.08}}
\put(3.75,2.125){\circle{0.08}}
\put(4,2){\circle{0.08}}
\put(4.25,1.875){\circle{0.08}}
\put(4.5,1.75){\circle{0.08}}
\put(4.75,1.625){\circle{0.08}}
\put(5,1.5){\circle{0.08}}
\put(5.25,1.375){\circle{0.08}}
\put(5.5,1.25){\circle{0.08}}
\put(5.75,1.125){\circle{0.08}}
\put(6,1){\circle{0.08}}
\put(6.25,0.875){\circle{0.08}}
\put(6.5,0.75){\circle{0.08}}
\put(6.75,0.625){\circle{0.08}}
\put(7,0.5){\circle{0.08}}
\put(7.25,0.375){\circle{0.08}}
\put(7.5,0.25){\circle{0.08}}
\put(7.75,0.125){\circle{0.08}}
\put(8,0){\circle{0.08}}
\put(8,0){\circle{0.12}}

\put(8.25,0.125){\circle{0.08}}
\put(8.5,0.25){\circle{0.08}}
\put(8.75,0.375){\circle{0.08}}
\put(9,0.5){\circle{0.08}}
\put(9.25,0.625){\circle{0.08}}
\put(9.5,0.75){\circle{0.08}}
\put(9.75,0.875){\circle{0.08}}
\put(10,1){\circle{0.08}}
\put(10.25,1.125){\circle{0.08}}
\put(10.5,1.25){\circle{0.08}}
\put(10.75,1.375){\circle{0.08}}
\put(11,1.5){\circle{0.08}}
\put(11.25,1.625){\circle{0.08}}
\put(11.5,1.75){\circle{0.08}}
\put(11.75,1.875){\circle{0.08}}
\put(12,2){\circle{0.08}}
\put(12,2){\circle{0.12}}

\put(0,0){\vector(1,0){14}}
\put(0,0){\vector(0,1){5.5}}

\put(0,0){\circle{0.12}}
\put(0,0){\circle{0.08}}

\put(5,1){\circle{0.12}}
\put(5,1){\circle{0.08}}

\put(0,4){\circle{0.12}}
\put(0,4){\circle{0.08}}

\put(12,4.5){\circle{0.12}}
\put(12,4.5){\circle{0.08}}

{\small  

\put(-6,4){$A_{0}'=(0,\nu _{q}(a_{0}))$}

\put(6.3,-0.85){$A_{1}'=(j_{2},0)$}

\put(-1.9,1.15){$A_{1}=(j_{1},\nu _{p}(a_{j_{1}}))$}

\put(12.3,4){$A_{2}=(n,\nu _{p}(a_{n}))$}

\put(12.3,2){$A_{2}'=(n,\nu _{q}(a_{n}))$}

\put(-3,-0.85){$A_{0}=(0,0)$}

\put(11.85,-0.85){$(n,0)$}
}

\end{picture}
\end{center}
\bigskip

{\em  {\small 
{\bf Figure 11.} The Newton polygons of $f$ with respect to $p$ and $q$, one consisting of two segments with positive different slopes, and the other one consisting of two segments whose slopes have different sign.}}
\medskip

Now, using again the notation in the statement of Theorem A, we see that 
\[
\mathcal{S}_{p}=\{ \min(j_{1},n-j_{1}) \}  \quad  \makebox{and} \quad \mathcal{S}_{q}=\{ \min(j_{2},n-j_{2}) \} , 
\]
so $\mathcal{S}_{p}\cap \mathcal{S}_{q}=\emptyset$, since $j_{1}\neq j_{2}$ and $j_{1}+j_{2}\neq n$.  $\square $
\medskip

{\it Proof of Theorem \ref{teorema17}.}\ In this case the Newton polygon of $f$ with respect to $p$ consists of an edge situated on the $x$-axis, followed by a segment  with positive slope, while the Newton polygon of $f$ with respect to $q$ consists of two segments with positive different slopes, as in Figure 12 below.

\begin{center}
\setlength{\unitlength}{4.65mm}
\begin{picture}(12,6)
\linethickness{0.075mm}

\put(12,0){\line(0,1){5}}   
\thicklines

\put(0,0){\line(1,0){4}}   

\put(4,0){\line(2,1){8}}   

\linethickness{0.15mm}

\put(4,0){\circle{0.08}}
\put(4,0){\circle{0.12}}

\put(12,4){\circle{0.12}}
\put(12,4){\circle{0.08}}

\put(0,0){\circle{0.12}}
\put(0,0){\circle{0.08}}
\put(0.25,0.025){\circle{0.08}}
\put(0.5,0.05){\circle{0.08}}
\put(0.75,0.075){\circle{0.08}}
\put(1,0.1){\circle{0.08}}
\put(1.25,0.125){\circle{0.08}}
\put(1.5,0.15){\circle{0.08}}
\put(1.75,0.175){\circle{0.08}}
\put(2,0.2){\circle{0.08}}
\put(2.25,0.225){\circle{0.08}}
\put(2.5,0.25){\circle{0.08}}
\put(2.75,0.275){\circle{0.08}}
\put(3,0.3){\circle{0.08}}
\put(3.25,0.325){\circle{0.08}}
\put(3.5,0.35){\circle{0.08}}
\put(3.75,0.375){\circle{0.08}}
\put(4,0.4){\circle{0.08}}
\put(4.25,0.425){\circle{0.08}}
\put(4.5,0.45){\circle{0.08}}
\put(4.75,0.475){\circle{0.08}}
\put(5,0.5){\circle{0.08}}
\put(5.25,0.525){\circle{0.08}}
\put(5.5,0.55){\circle{0.08}}
\put(5.75,0.575){\circle{0.08}}
\put(6,0.6){\circle{0.08}}
\put(6.25,0.625){\circle{0.08}}
\put(6.5,0.65){\circle{0.08}}
\put(6.75,0.675){\circle{0.08}}
\put(7,0.7){\circle{0.08}}
\put(7.25,0.725){\circle{0.08}}
\put(7.5,0.75){\circle{0.08}}
\put(7.75,0.775){\circle{0.08}}
\put(8,0.8){\circle{0.08}}
\put(8,0.8){\circle{0.12}}
\put(8.25,0.9){\circle{0.08}}
\put(8.5,1){\circle{0.08}}
\put(8.75,1.1){\circle{0.08}}
\put(9,1.2){\circle{0.08}}
\put(9.25,1.3){\circle{0.08}}
\put(9.5,1.4){\circle{0.08}}
\put(9.75,1.5){\circle{0.08}}
\put(10,1.6){\circle{0.08}}
\put(10.25,1.7){\circle{0.08}}
\put(10.5,1.8){\circle{0.08}}
\put(10.75,1.9){\circle{0.08}}
\put(11,2){\circle{0.08}}
\put(11.25,2.1){\circle{0.08}}
\put(11.5,2.2){\circle{0.08}}
\put(11.75,2.3){\circle{0.08}}
\put(12,2.4){\circle{0.08}}
\put(12,2.4){\circle{0.12}}

\put(0,0){\vector(1,0){14.9}}
\put(0,0){\vector(0,1){5.5}}

{\small  
\put(3.5,-0.85){$A_{1}=(j_{1},0)$}

\put(8.7,0.35){$A_{1}'=(j_{2},\nu _{q}(a_{j_{2}}))$}

\put(12.3,2.4){$A_{2}'=(n,\nu _{q}(a_{n}))$}

\put(12.3,4){$A_{2}=(n,\nu _{p}(a_{n}))$}

\put(-4,-0.85){$A_{0}=A_{0}'=(0,0)$}

\put(11.85,-0.85){$(n,0)$}
}

\end{picture}
\end{center}
\bigskip

{\em {\small
{\bf Figure 12.} The Newton polygons of $f$ with respect to $p$ and $q$, one consisting of an edge situated on the $x$-axis followed by a segment with positive slope, and the other one consisting of two segments with positive different slopes.}}
\medskip

In this case we see that 
\[
\mathcal{S}_{p}=\{ 1,\dots ,j_{1} \} \quad  \makebox{and} \quad \mathcal{S}_{q}=\{ \min(j_{2},n-j_{2}) \} , 
\]
so $\mathcal{S}_{p}\cap \mathcal{S}_{q}=\emptyset$, since $j_{1}<\min(j_{2},n-j_{2})$.  $\square $
\medskip

{\it Proof of Theorem \ref{teorema18}.} \ In this case the Newton polygon of $f$ with respect to $p$ consists of an edge situated on the $x$-axis, followed by a segment  with positive slope, while the Newton polygon of $f$ with respect to $q$ consists of two segments with negative different slopes, as in Figure 13 below.

\begin{center}
\setlength{\unitlength}{4.65mm}
\begin{picture}(12,6)
\linethickness{0.075mm}

\put(12,0){\line(0,1){4.5}}   
\thicklines

\put(0,0){\line(1,0){4}}   

\put(4,0){\line(2,1){8}}   

\linethickness{0.15mm}

\put(0,4){\circle{0.12}}
\put(0,4){\circle{0.08}}
\put(0.25,3.875){\circle{0.08}}
\put(0.5,3.75){\circle{0.08}}
\put(0.75,3.625){\circle{0.08}}
\put(1,3.5){\circle{0.08}}
\put(1.25,3.375){\circle{0.08}}
\put(1.5,3.25){\circle{0.08}}
\put(1.75,3.125){\circle{0.08}}
\put(2,3){\circle{0.08}}
\put(2.25,2.875){\circle{0.08}}
\put(2.5,2.75){\circle{0.08}}
\put(2.75,2.625){\circle{0.08}}
\put(3,2.5){\circle{0.08}}
\put(3.25,2.375){\circle{0.08}}

\put(3.5,2.25){\circle{0.08}}
\put(3.75,2.125){\circle{0.08}}
\put(4,2){\circle{0.08}}
\put(4.25,1.875){\circle{0.08}}
\put(4.5,1.75){\circle{0.08}}
\put(4.75,1.625){\circle{0.08}}
\put(5,1.5){\circle{0.08}}
\put(5.25,1.375){\circle{0.08}}
\put(5.5,1.25){\circle{0.08}}
\put(5.75,1.125){\circle{0.08}}
\put(6,1){\circle{0.08}}
\put(6.25,0.875){\circle{0.08}}
\put(6.5,0.75){\circle{0.08}}
\put(6.75,0.625){\circle{0.08}}

\put(7,0.5){\circle{0.12}}
\put(7,0.5){\circle{0.08}}
\put(7.25,0.475){\circle{0.08}}
\put(7.5,0.45){\circle{0.08}}
\put(7.75,0.425){\circle{0.08}}
\put(8,0.4){\circle{0.08}}

\put(8.25,0.375){\circle{0.08}}
\put(8.5,0.35){\circle{0.08}}
\put(8.75,0.325){\circle{0.08}}
\put(9,0.3){\circle{0.08}}
\put(9.25,0.275){\circle{0.08}}
\put(9.5,0.25){\circle{0.08}}
\put(9.75,0.225){\circle{0.08}}
\put(10,0.2){\circle{0.08}}
\put(10.25,0.175){\circle{0.08}}
\put(10.5,0.15){\circle{0.08}}
\put(10.75,0.125){\circle{0.08}}
\put(11,0.1){\circle{0.08}}
\put(11.25,0.075){\circle{0.08}}
\put(11.5,0.05){\circle{0.08}}
\put(11.75,0.025){\circle{0.08}}
\put(12,0){\circle{0.08}}
\put(12,0){\circle{0.12}}

\put(0,0){\vector(1,0){14}}
\put(0,0){\vector(0,1){5.5}}

\put(0,0){\circle{0.12}}
\put(0,0){\circle{0.08}}

\put(4,0){\circle{0.12}}
\put(4,0){\circle{0.08}}

\put(0,4){\circle{0.12}}
\put(0,4){\circle{0.08}}

\put(12,4){\circle{0.12}}
\put(12,4){\circle{0.08}}

{\small  

\put(-6,4){$A_{0}'=(0,\nu _{q}(a_{0}))$}

\put(7,0.85){$A_{1}'=(j_{2},\nu _{q}(a_{j_{2}}))$}

\put(4,-0.85){$A_{1}=(j_{1},0)$}

\put(12.3,4){$A_{2}=(n,\nu _{p}(a_{n}))$}

\put(-3,-0.85){$A_{0}=(0,0)$}

\put(11.85,-0.85){$A_{2}'=(n,0)$}
}

\end{picture}
\end{center}
\bigskip

{\em {\small 
{\bf Figure 13.} The Newton polygons of $f$ with respect to $p$ and $q$, one consisting of an edge situated on the $x$-axis followed by a segment with positive slope, and the other one consisting of two segments with negative different slopes.}}
\medskip

Again, we deduce that
\[
\mathcal{S}_{p}=\{ 1,\dots ,j_{1} \} \quad  \makebox{and} \quad \mathcal{S}_{q}=\{ \min(j_{2},n-j_{2}) \} , 
\]
so $\mathcal{S}_{p}\cap \mathcal{S}_{q}=\emptyset$, since $j_{1}<\min(j_{2},n-j_{2})$.  $\square $
\medskip

{\it Proof of Theorem \ref{teorema5}.} The conditions in the statement of the theorem show that each one of the Newton polygons of $f$ with respect to $p$ and $q$ consists of only two segments,
whose intersection point lies above the $x$-axis, as in Figure 14 below.

\begin{center}
\setlength{\unitlength}{4.65mm}
\begin{picture}(12,6)
\linethickness{0.075mm}

\put(12,0){\line(0,1){5}}   
\thicklines

\put(0,0){\line(5,1){6}}   

\put(6,1.2){\line(2,1){6}}   
\put(6,1.2){\circle{0.08}}
\put(6,1.2){\circle{0.12}}
\put(12,4.2){\circle{0.12}}

\put(0,0){\circle{0.12}}
\put(0,0){\circle{0.08}}
\put(0.25,0.025){\circle{0.08}}
\put(0.5,0.05){\circle{0.08}}
\put(0.75,0.075){\circle{0.08}}
\put(1,0.1){\circle{0.08}}
\put(1.25,0.125){\circle{0.08}}
\put(1.5,0.15){\circle{0.08}}
\put(1.75,0.175){\circle{0.08}}
\put(2,0.2){\circle{0.08}}
\put(2.25,0.225){\circle{0.08}}
\put(2.5,0.25){\circle{0.08}}
\put(2.75,0.275){\circle{0.08}}
\put(3,0.3){\circle{0.08}}
\put(3.25,0.325){\circle{0.08}}
\put(3.5,0.35){\circle{0.08}}
\put(3.75,0.375){\circle{0.08}}
\put(4,0.4){\circle{0.08}}
\put(4.25,0.425){\circle{0.08}}
\put(4.5,0.45){\circle{0.08}}
\put(4.75,0.475){\circle{0.08}}
\put(5,0.5){\circle{0.08}}
\put(5.25,0.525){\circle{0.08}}
\put(5.5,0.55){\circle{0.08}}
\put(5.75,0.575){\circle{0.08}}
\put(6,0.6){\circle{0.08}}
\put(6.25,0.625){\circle{0.08}}
\put(6.5,0.65){\circle{0.08}}
\put(6.75,0.675){\circle{0.08}}
\put(7,0.7){\circle{0.08}}
\put(7.25,0.725){\circle{0.08}}
\put(7.5,0.75){\circle{0.08}}
\put(7.75,0.775){\circle{0.08}}
\put(8,0.8){\circle{0.08}}
\put(8,0.8){\circle{0.12}}
\put(8.25,0.9){\circle{0.08}}
\put(8.5,1){\circle{0.08}}
\put(8.75,1.1){\circle{0.08}}
\put(9,1.2){\circle{0.08}}
\put(9.25,1.3){\circle{0.08}}
\put(9.5,1.4){\circle{0.08}}
\put(9.75,1.5){\circle{0.08}}
\put(10,1.6){\circle{0.08}}
\put(10.25,1.7){\circle{0.08}}
\put(10.5,1.8){\circle{0.08}}
\put(10.75,1.9){\circle{0.08}}
\put(11,2){\circle{0.08}}
\put(11.25,2.1){\circle{0.08}}
\put(11.5,2.2){\circle{0.08}}
\put(11.75,2.3){\circle{0.08}}
\put(12,2.4){\circle{0.08}}
\put(12,2.4){\circle{0.12}}

\linethickness{0.15mm}

\put(0,0){\vector(1,0){14.9}}
\put(0,0){\vector(0,1){5.5}}

{\small  
\put(0.5,1.7){$A_{1}=(j_{1},\nu _{p}(a_{j_{1}}))$}

\put(8.7,0.35){$A_{1}'=(j_{2},\nu _{q}(a_{j_{2}}))$}

\put(12.3,2.4){$A_{2}'=(n,\nu _{q}(a_{n}))$}

\put(5.8,4.2){$A_{2}=(n,\nu _{p}(a_{n}))$}

\put(-4,-0.85){$A_{0}=A_{0}'=(0,0)$}

\put(11.85,-0.85){$(n,0)$}
}

\end{picture}
\end{center}
\bigskip

{\em {\small
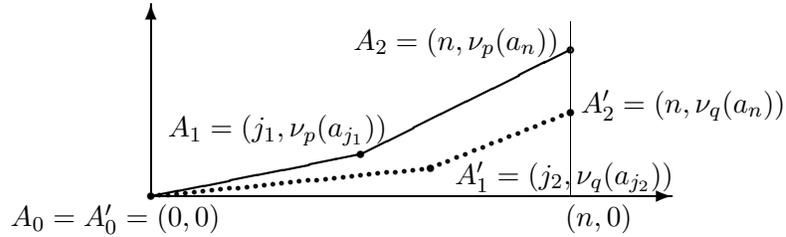
{\bf Figure 14.} The Newton polygons of $f$ with respect to $p$ and $q$, each one consisting of two segments with positive, different slopes.}}
\medskip

Note that the conditions $\frac{\nu _{p}(a_{n})}{n}>\frac{\nu _{p}(a_{j_{1}})}{j_{1}}$   and  $\frac{\nu _{q}(a_{n})}{n}>\frac{\nu _{q}(a_{j_{2}})}{j_{2}}$ show that the slopes of the segments $A_{1}A_{2}$ and $A_{1}'A_{2}'$ exceed
the slopes of the segments $A_{0}A_{1}$ and $A_{0}'A_{1}'$, respectively.
Using again the notation in the statement of Theorem A, we see that 
\[
\mathcal{S}_{p}=\{ \min(j_{1},n-j_{1}) \} \quad  \makebox{and} \quad \mathcal{S}_{q}=\{ \min(j_{2},n-j_{2}) \} , 
\]
hence $\mathcal{S}_{p}\cap \mathcal{S}_{q}=\emptyset$, since $j_{1}\neq j_{2}$ and $j_{1}+j_{2}\neq n$.  $\square $
\medskip

{\it Proof of Theorem \ref{teorema6}.} The conditions in the statement of the theorem show that the Newton polygon of $f$ with respect to $p$ consists of only two segments with negative different slopes,
while the Newton polygon of $f$ with respect to $q$ consists of only two segments with positive different slopes, as in Figure 15 below. Here the conditions $\frac{\nu _{p}(a_{0})}{n}>\frac{\nu _{p}(a_{j_{1}})}{n-j_{1}}$
and  $\frac{\nu _{q}(a_{n})}{n}>\frac{\nu _{q}(a_{j_{2}})}{j_{2}}$ show that the slopes of the segments $A_{1}A_{2}$ and $A_{1}'A_{2}'$ exceed the slopes of the segments
$A_{0}A_{1}$ and $A_{0}'A_{1}'$, respectively.

\begin{center}

\setlength{\unitlength}{4.65mm}
\begin{picture}(12,6)
\linethickness{0.075mm}

\put(12,0){\line(0,1){5}}   
\thicklines

\put(0,4){\line(2,-1){6}}   
\put(6,1){\line(6,-1){6}}   
\linethickness{0.15mm}

\put(0,4){\circle{0.08}}
\put(0,4){\circle{0.12}}
\put(6,1){\circle{0.08}}
\put(6,1){\circle{0.12}}
\put(12,0){\circle{0.08}}
\put(12,0){\circle{0.12}}

\put(0,0){\circle{0.12}}
\put(0,0){\circle{0.08}}
\put(0.25,0.025){\circle{0.08}}
\put(0.5,0.05){\circle{0.08}}
\put(0.75,0.075){\circle{0.08}}
\put(1,0.1){\circle{0.08}}
\put(1.25,0.125){\circle{0.08}}
\put(1.5,0.15){\circle{0.08}}
\put(1.75,0.175){\circle{0.08}}
\put(2,0.2){\circle{0.08}}
\put(2.25,0.225){\circle{0.08}}
\put(2.5,0.25){\circle{0.08}}
\put(2.75,0.275){\circle{0.08}}
\put(3,0.3){\circle{0.08}}
\put(3.25,0.325){\circle{0.08}}
\put(3.5,0.35){\circle{0.08}}
\put(3.75,0.375){\circle{0.08}}
\put(4,0.4){\circle{0.08}}
\put(4.25,0.425){\circle{0.08}}
\put(4.5,0.45){\circle{0.08}}
\put(4.75,0.475){\circle{0.08}}
\put(5,0.5){\circle{0.08}}
\put(5.25,0.525){\circle{0.08}}
\put(5.5,0.55){\circle{0.08}}
\put(5.75,0.575){\circle{0.08}}
\put(6,0.6){\circle{0.08}}
\put(6.25,0.625){\circle{0.08}}
\put(6.5,0.65){\circle{0.08}}
\put(6.75,0.675){\circle{0.08}}
\put(7,0.7){\circle{0.08}}
\put(7.25,0.725){\circle{0.08}}
\put(7.5,0.75){\circle{0.08}}
\put(7.75,0.775){\circle{0.08}}
\put(8,0.8){\circle{0.08}}
\put(8.25,0.825){\circle{0.08}}
\put(8.5,0.85){\circle{0.08}}
\put(8.75,0.875){\circle{0.08}}
\put(9,0.9){\circle{0.08}}
\put(9,0.9){\circle{0.12}}
\put(9.25,1){\circle{0.08}}
\put(9.5,1.2){\circle{0.08}}
\put(9.75,1.4){\circle{0.08}}
\put(10,1.6){\circle{0.08}}
\put(10.25,1.8){\circle{0.08}}
\put(10.5,2){\circle{0.08}}
\put(10.75,2.2){\circle{0.08}}
\put(11,2.4){\circle{0.08}}
\put(11.25,2.6){\circle{0.08}}
\put(11.5,2.8){\circle{0.08}}
\put(11.75,3){\circle{0.08}}
\put(12,3.2){\circle{0.08}}
\put(12,3.2){\circle{0.12}}

\put(0,0){\vector(1,0){14.5}}
\put(0,0){\vector(0,1){5.5}}

{\small  
\put(-1,1){$A_{1}=(j_{1},\nu _{p}(a_{j_{1}}))$}

\put(10,0.65){$A_{1}'=(j_{2},\nu _{q}(a_{j_{2}}))$}

\put(12.3,3.2){$A_{2}'=(n,\nu _{q}(a_{n}))$}

\put(-5.8,4){$A_{0}=(0,\nu _{p}(a_{0}))$}

\put(-3,-0.85){$A_{0}'=(0,0)$}

\put(11.85,-0.85){$A_{2}=(n,0)$}
}

\end{picture}
\end{center}
\bigskip

{\em {\small
{\bf Figure 15.} The Newton polygons of $f$ with respect to $p$ and $q$, one consisting of two segments with negative, different slopes, and the other one consisting of two segments with positive, different slopes.}}
\medskip

As in the proof of Theorem \ref{teorema5}, we see that 
\[
\mathcal{S}_{p}=\{ \min(j_{1},n-j_{1}) \} \quad  \makebox{and} \quad \mathcal{S}_{q}=\{ \min(j_{2},n-j_{2}) \} , 
\]
so $\mathcal{S}_{p}\cap \mathcal{S}_{q}=\emptyset$, since $j_{1}\neq j_{2}$ and $j_{1}+j_{2}\neq n$.  $\square $
\medskip

{\it Proof of Theorem \ref{teorema7}.} The conditions in the statement of the theorem show that each one of the Newton polygons of $f$ with respect to $p_{1},\dots ,p_{k}$ consists of only two segments whose slopes
have the same sign, as in Figure 15 above. One applies now Theorem B.  $\square $

\section{Examples} \label{se3}

We end with several examples of polynomials whose irreducibility can not be proved by a direct use of Sch\" onemann-Eisenstein criterion, or of Dumas' criterion. 

\medskip

1) Let $f(X)=q^3+p^2q^3X+p^2q^3X^2+p^2q^3X^3+p^2q^3X^4+p^2q^3X^5+p^2X^6$, where $p,q$ are distinct prime numbers. One can easily see that $f$ is irreducible by Corollary \ref{corolarul1}.

2) Let $f(X)=p+pq^2X+pq^2X^2+pq^2X^3+q^2X^4+pq^2X^5+pq^2X^6$, where $p,q$ are distinct prime numbers. In this case we deduce that $f$ is irreducible by using Corollary \ref{corolarul13} with  $k=2$ and $j=4$.

3) Let $f(X)=p^m+q^2X+pq^2X^2+pq^2X^3+pq^2X^4+pq^2X^5+pq^2X^6$, where $p,q$ are distinct prime numbers and $m$ is an arbitrary non-negative integer. Here we deduce that $f$ is irreducible by using Corollary \ref{corolarul15} with $k=2$.

4) Let $f(X)=p+pq^2X+pq^2X^2+pq^2X^3+pq^2X^4+q^2X^5+p^mq^2X^6$, where $p,q$ are distinct prime numbers and $m$ is an arbitrary non-negative integer. According to Corollary \ref{corolarul17} with $k=2$, $f$ must be irreducible.

5) Let $f(X)=1+q^2X+q^2X^2+pq^2X^3+pq^2X^4+pq^2X^5+pq^2X^6$, where $p,q$ are distinct prime numbers. Here we deduce that $f$ is irreducible by Corollary \ref{corolarul19} with $j=k=2$.

6) Let $f(X)=p+pq^2X+pq^2X^2+pq^2X^3+q^2X^4+q^2X^5+q^2X^6$, where $p,q$ are distinct prime numbers. By Corollary \ref{corolarul21} with $k=2$ and $j=4$ we deduce that $f$ is irreducible.

7) Let $f(X)=pq+pqX+pqX^2+qX^3+pX^4+pqX^5+pqX^6$, where $p,q$ are distinct prime numbers. In this case we see that $f$ is irreducible by Corollary \ref{corolarul3} with $j_{1}=3$ and $j_{2}=4$.

8) Let $f(X)=q+qX+qX^2+pX^3+pqX^4+pqX^5+pqX^6$, where $p,q$ are distinct prime numbers. Here we see that $f$ is irreducible by Corollary \ref{corolarul27} with $j_{1}=2$ and $j_{2}=3$.

\end{document}